\documentclass{m2an}
\usepackage{graphicx, amssymb, amsmath, epsfig}
\usepackage{color}
\usepackage{hyperref}
\usepackage{booktabs}
\usepackage{boxedminipage}
\usepackage{algorithmic,algorithm}
\usepackage[utf8]{inputenc}
\usepackage{epstopdf}       
\usepackage{float}

\graphicspath{{./}{figs/}}

\newcommand{\R}{\mathbb R}

\newtheorem{teo}{Theorem}[section]

\begin{document}
\title{A-posteriori snapshot location for POD \\in optimal control of linear parabolic equations}\thanks{We would like to thank Z J. Zhou from Shandong Normal University, China for providing the data and code of the space-time approximation in \cite{GHZ12}. First author also acknowledges the support of US Department of Energy (grant number DE-SC0009324). }
\author{Alessandro Alla}\address{Florida State University, Department of Scientific Computing, FL-32306 Tallahassee, USA  ({\texttt aalla@fsu.edu})}
\author{Carmen Gr\"assle}\address{University of Hamburg, Department of Mathematics, D-20146 Hamburg, Germany ({\texttt carmen.graessle@uni-hamburg.de})}
\author{Michael Hinze}\address{University of Hamburg, Department of Mathematics, D-20146 Hamburg, Germany ({\texttt michael.hinze@uni-hamburg.de})}
\date{30 August 2016}
\begin{abstract} 
In this paper we study the approximation of an optimal control problem for linear para\-bolic PDEs with model order reduction based on Proper Orthogonal Decomposition (POD-MOR). POD-MOR is a Galerkin approach where the basis functions are obtained upon information contained in time snapshots of the parabolic PDE related to given input data. In the present work we show that for POD-MOR in optimal control of parabolic equations it is important to have knowledge about the controlled system at the right time instances. We propose to determine the time instances (snapshot locations) by an a-posteriori error control concept. This method is based on a reformulation of the optimality system of the underlying optimal control problem as a second order in time and fourth order in space elliptic system which is approximated by a space-time finite element method.  
Finally, we present numerical tests to illustrate our approach and to show the effectiveness of the method in comparison to existing approaches.
\end{abstract}
%
%
\maketitle
\section{Introduction}

Optimization with PDE constraints is nowadays a well-studied topic motivated by its relevance in industrial applications. We are interested in the numerical approximation of such optimization problems in an efficient and reliable way using surrogate models obtained with POD-MOR. The surrogate models are built upon {\em snapshots} of the system to provide information about the underlying problem. This stage is usually called the offline stage. For the snapshot POD approach we refer the reader to \cite{Sir87}.\\
\noindent
Several works focus their attention on the choice of the snapshots, in order to approximate either dynamical systems or optimal control problems by suitable surrogate models. In \cite{KV10}, it is proposed to optimize the choice of the time instances such that the error between POD and the trajectory of the dynamical system is minimized. A recent approach proposes to choose the snapshots by an a-posteriori error estimator in order to equidistribute the state error on the time grid related to the snapshot locations (see \cite{HL14}). We also mention an adaptive method, proposed in \cite{OKAC15}, where the aim is to reduce expensive offline costs selecting the snapshots according to an adaptive time-stepping algorithm using time error-control. For further references we refer the interested reader to \cite{OKAC15}.\\
In optimal control problems the reduced model is usually built upon a forecast on the control. This approach does not guarantee a proper construction of the surrogate model since we do not know how far away the optimal solution is from the reference control. More sophisticated approaches select snapshots by solving an optimization problem in order to improve the selection of the snapshots according to the desired controlled dynamics. For this purpose optimality system for POD (OS-POD) is introduced in \cite{KV08}. In OS-POD, the computation of the basis functions is performed by means of the solution of an enlarged optimal control problem which involves the full problem, the reduced equation and the eigenvalue problem for the POD modes.\\
The reduction of optimal control problems with particular focus on adaptive adjustments of the surrogate models can be found in \cite{AH01,AFS02}. We should also mention another adaptive method for feedback control problems by means of the Hamilton-Jacobi-Bellman equation, introduced in \cite{AF12}.\\
Recently, an a-posteriori error estimator was introduced in \cite{TV09,KTV13} for optimal control problems. In these works the error between the unknown optimal and the computed POD suboptimal control is estimated for linear and nonlinear problems, and it is shown that increasing the number of basis functions leads to the desired convergence. OS-POD and a-posteriori error estimation is combined in \cite{V11}.\\
All these works have in common that they compute basis functions for optimal control problems. In our paper we address the question of an efficient and suitable selection of snapshot locations by means of an a-posteriori error control approach proposed in \cite{GHZ12}. We rewrite the optimality conditions as a second order in time and fourth order in space elliptic equation for the adjoint variable and we generalize this approach to control constraints. In particular, a time adaptive concept is used to build the snapshot grid which should be used to construct the POD surrogate model for the approximate solution of the optimal control problem. Here the novelty for the reduced control problem is twofold: we directly obtain snapshots related to an approximation of the 
optimal control and, at the same time, we get information about the time grid. 

We have proposed a similar approach based on a reformulation of the optimality system 
with respect to the state variable in \cite{AGH15}. Now, we focus our approach on the adjoint variable 
and generalize the idea presented in \cite{GHZ12} to time dependent control intensities with control 
shape functions including control constraints. Furthermore, we 
certify our approach by means of several error bounds for the state, adjoint state and control 
variable.\\
The outline of this paper is as follows. In Section 2 we present the 
optimal control problem together with the optimality conditions. In Section 
3 we recall the main results of \cite{GHZ12}. Proper Orthogonal 
Decomposition and its application to optimal control problems is 
presented in Section 4. The focus of Section 5 lies in 
investigating our snapshot location strategy. Finally, numerical tests are discussed in Section 6 and conclusions are driven in Section 7.

\section{Optimal Control Problem}

In this section we describe the optimal control problem. The governing equation is given by a 
linear parabolic PDE:
\begin{equation}\label{heat}
\left.
\begin{array}{rcll}
y_t-\Delta y & = & f+\mathcal{B}u &\text{ in } \Omega_T,\\
y(\cdot,0) & = & y_0 & \text{ in } \Omega,\\
y &= & 0 &\text{ on } \Sigma_T,
\end{array}
\right\}
\end{equation}
\noindent
where $\Omega\subset\R^q, q \in \{1,2,3\}$ is an open bounded domain with smooth 
boundary, $T>0$, $\Omega_T:=\Omega\times (0,T]$ is 
the space-time cylinder, $\Sigma_T:=\partial\Omega\times (0,T]$, and the state 
is denoted by  $ y:\Omega_T\rightarrow\R$. As control space we use 
$\left(L^2(0,T;\mathbb{R}^m), \langle \cdot,\cdot \rangle_U\right)$, where 
$\langle u,v\rangle_U:=\sum_{i=1}^m \langle u_i,v_i\rangle_{L^2(0,T)}$, and define the control 
operator as $\mathcal{B}:U \rightarrow L^2(0,T;H^{-1}(\Omega))$, $(\mathcal{B}u)(t) = 
\sum_{i=1}^m u_i(t) \chi_i$, where $\chi_i \in H^{-1}(\Omega) (1 \leq i \leq m)$ denote 
specified control actions. Thus $\mathcal{B}$ is linear and bounded. For the control 
variable we require
$$u\in U_{ad} := \{u\in U \; | \; u_a(t) \leq u(t) \leq u_b(t) \text{ in } 
\mathbb{R}^m \text{ a.e. in } [0,T] \} \subset L^\infty(0,T;\mathbb{R}^m)$$ 
with $u_a, u_b \in L^\infty(0,T;\mathbb{R}^m), u_a(t) \leq u_b(t)$ almost everywhere in $(0,T)$. 
It is well-known (see \cite{L71}, for example) that for a given initial condition $y_0\in L^2(\Omega)$ and a 
forcing term $f\in L^2(0,T;H^{-1}(\Omega))$ the equation \eqref{heat} 
admits a unique solution $y=y(u)\in W(0,T)$, where
$$W(0,T):=\left\{v\in L^2\left(0,T;H^1_0(\Omega)\right), \dfrac{\partial v}{\partial t}\in L^2\left(0,T;H^{-1}(\Omega)\right)\right\}.$$
If $y_0 \in H_0^1(\Omega)$, higher regularity results can be derived 
according to \cite{E10}. We also note that the unconstrained case is related to $u_a\equiv-\infty, u_b\equiv+\infty$.

\noindent
 The weak formulation of \eqref{heat} is given by: find $y \in W(0,T)$ with $y(0)=y_0$ and 
\begin{equation}\label{weak:heat}
  \int_\Omega y_t(t)v dx + \int_\Omega \nabla y(t) \cdot \nabla v dx
=\int_\Omega  (f+\mathcal{B}u)(t) v dx \quad \forall v \in H_0^1(\Omega).
\end{equation}


\noindent The cost functional we want to minimize is given by
\begin{equation}
J(y,u):=\dfrac{1}{2} \|y-y_d\|^2_{L^2(\Omega_T)}+\dfrac{\alpha}{2} \|u\|^2_U,
\end{equation}
where $y_d\in L^2(\Omega_T)$ is the desired state and the regularization parameter $\alpha$ is a real positive constant.
The optimal control problem then reads
\begin{equation}\label{ocp}
\min_{u\in U_{ad}} \hat{J}(u):=J(y(u),u) \mbox{, where } y(u) \mbox{ satisfies } \eqref{heat}.
\end{equation}
Note that $U_{ad}$ is a non-empty, bounded, convex and closed subset of $L^\infty(0,T;\mathbb{R}^m)$.
Hence, it is easy to argue that \eqref{ocp} admits a unique solution $u \in U$ with associated state $y(u)\in W(0,T)$, see e.g. \cite{L71}.

\noindent
The first order optimality system of the optimal control problem \eqref{ocp} is given by the state equation \eqref{heat}, together with the adjoint equation
\begin{equation}\label{adj}
\left.
\begin{array}{rcll}
-p_t-\Delta p & = & y-y_d &\text{ in } \Omega_T,\\
p(\cdot,T) & = & 0 &\text{ in } \Omega,\\
p & = & 0 &\text{ on } \Sigma_T,
\end{array}
\right\}
\end{equation}
and the variational inequality
\begin{equation}\label{opt_con}
\langle\alpha u + \mathcal{B}^*p,v-u \rangle_U \geq 0\quad \mbox{ for all }v \in U_{ad},
\end{equation}
\noindent where $\mathcal{B}^* : L^2(0,T;H^{-1}(\Omega))^* \to U^*$ is the dual 
operator of $\mathcal{B}$. In \eqref{opt_con} we have identified $L^2(0,T;H^{-1}(\Omega))^*$ 
with $L^2(0,T; H^1_0(\Omega))$ and $U^*$ with $U$, where we use that Hilbert spaces are reflexive. 
The variational inequality \eqref{opt_con} is equivalent to the projection formula
\begin{equation}\label{Proj}
 u(t)=\mathcal{P}_{U_{ad}}\left\lbrace -\dfrac{1}{\alpha}(\mathcal{B}^*p)(t) \right\rbrace \mbox{ for almost all } t\in[0,T],
\end{equation}
where $\mathcal{P}_{U_{ad}}:U\rightarrow U_{ad}$ denotes the orthogonal projection onto $U_{ad}$. It 
follows from the reflexivity of the involved spaces that the action of the adjoint operator 
$\mathcal{B}^*$ is given as $$(\mathcal{B}^*v)(t)=\left(\langle \chi_1,v\rangle_{H^{-1},H^1_0},\ldots,
\langle \chi_m,v\rangle_{H^{-1},H^1_0}\right)$$ and
\begin{equation*}
 \mathcal{P}_{U_{ad}}\left\lbrace-\frac{1}{\alpha} \mathcal{B}^*p\right\rbrace_i = \max \left\{u_a , \min \{ 
u_b , -\frac{1}{\alpha} \langle \chi_i , p \rangle_{H^{-1},H_0^1}\} \right\}.
\end{equation*}

Since our domain is smooth, the regularities of the optimal state, the optimal control and the 
associated adjoint state are limited through the regularities of the initial state $y_0$, the right 
hand side $f$, the control $\mathcal{B}u$ and the desired state $y_d.$\\
\noindent
The numerical approximation of the optimality 
system \eqref{heat}-\eqref{adj}-\eqref{opt_con} with a standard Finite Element 
Method (FEM) in the spatial variable leads to a high-dimensional system of ordinary 
differential equations:

\begin{equation}\label{opt_disc}
\left.
\begin{array}{rclrcl}
M\dot{y}^N- A y^N & = & f^N+ \mathcal{B}^Nu, & \quad y^N(0) & = & y_0^N,\\
-M\dot{p}^N-A p^N & = & y^N-y^N_d, & \quad p^N(T) & = & 0,\\
 \langle \alpha u + (\mathcal{B}^*)^Np^N, v-u \rangle_\mathcal{U} & \geq & 0. & & & 
\end{array}
\right\}
\end{equation}
\noindent
Here $y^N, p^N:[0,T]\rightarrow\R^N$ are the 
semi-discrete state and adjoint, respectively,  $\dot{y}^N, \dot{p}^N$ are the time derivatives,
$M\in\R^{N\times N}$ denotes the mass matrix and 
$A\in\R^{N\times N}$ the stiffness matrix. Note that the 
dimension $N$ of each equation in the semi-discrete system \eqref{opt_disc} is related to the number of element nodes chosen in the FEM approach.\\



\section{Space-Time approximation}

In this section, we consider the reformulation of the 
optimality system \eqref{heat}-\eqref{adj}-\eqref{opt_con} as an elliptic 
equation of fourth order in space and second order in time for the adjoint 
variable $p$. This is carried out for the unconstrained control 
problem in \cite{GHZ12} and generalized to control constrained optimal control 
problems in \cite{NPS11}. Following these works, we include control constraints. Here, we aim to derive an 
a-posteriori error estimate for the time discretization as suggested in \cite{GHZ12}, which 
then turns out to be the basis for our model reduction approach to solve \eqref{ocp}.\\
\noindent
We define 
$$H^{2,1}_0(\Omega_T):=\left\{v\in H^{2,1}(\Omega_T): v(T)=0 \mbox { in }\Omega\right\},$$
where 
$$H^{2,1}(\Omega_T)=L^2\left(0,T;H^2(\Omega\right)\cap H^1_0\left(\Omega)\right)\cap H^1\left(0,T; L^2(\Omega)\right)$$ is equipped with the norm
$$\|w\|_{H^{2,1}(\Omega_T)}^2:=\left(\|w\|^2_{L^2(0,T;H^2(\Omega))}+\|w\|^2_{H^1(0,T;L^2(\Omega))}\right).$$
\noindent Under the assumptions 
$y_0 \in H_0^1(\Omega)$, $\chi_i \in L^2(\Omega)$ for $i=1, \dotsc, m$ and $y_d \in H^{2,1}(\Omega_T)$, the regularity of
$y,p \in H^{2,1}(\Omega_T)$ is ensured, see \cite{E10} for 
the details. Then, the first order optimality 
conditions \eqref{heat}-\eqref{adj}-\eqref{opt_con} can be transformed into 
an initial boundary value problem for $p$ in space-time:

\begin{equation}\label{2ordp}
\left.
\begin{array}{rcll}
-p_{tt}+\Delta^2 p-\mathcal{B}\mathcal{P}_{U_{ad}}\left(-\dfrac{1}{\alpha}\mathcal{B}^*p\right) & = &-(y_d)_t+\Delta y_d & \text{ in } \Omega_T,\\
p(\cdot, T) & = & 0 &\text{ in } \Omega,\\
p & = & 0 &\text{ on } \Sigma_T,\\
\Delta p & = & y_d &\text{ on } \Sigma_T,\\
\left(p_t+\Delta p\right)(0) & = & y_d(0)-y_0 &\text{ in }\Omega,
\end{array}
\right\}
\end{equation}

\noindent where, without loss of generality, we have set $f\equiv0$. We note that the quantity 
$$\mathcal{B}\mathcal{P}_{U_{ad}}\left(-\dfrac{1}{\alpha}\mathcal{B}^*p\right)$$
is nondifferentiable and nonlinear in $p$ and thus \eqref{2ordp} 
becomes a semilinear second order in time and fourth order in space elliptic 
problem with a monotone nonlinearity. Existence of a unique weak solution for 
\eqref{2ordp} can be proved analogously to \cite{NPS11} and follows from the 
fact that the optimal control problem \eqref{ocp} in the case of control constraints 
with closed and convex $U_{ad} \subset U$ admits a unique solution.\\
In order to provide the weak formulation of \eqref{2ordp}, we define 
the operator $A_0$ and the linear form $L_0$ as
$$A_0:H^{2,1}_0(\Omega_T)\times H^{2,1}_0(\Omega_T)\rightarrow\R,\qquad L_0:H^{2,1}_0(\Omega_T)\rightarrow\R,$$
$$A_0(v,w):=\int_{\Omega_T} \left( v_tw_t-\mathcal{B}\mathcal{P}_{U_{ad}} \left(-\dfrac{1}{\alpha}\mathcal{B}^* v \right)w\right)+\int_{\Omega_T}\Delta v\Delta w +\int_\Omega \nabla v(0) \nabla w(0),$$
$$ L_0(v):=\int_{\Omega_T} \langle -\dfrac{\partial y_d}{\partial t}+\Delta y_d,v\rangle_{H^{-1}(\Omega)\times H^1_0(\Omega)}-\int_\Omega (y_d(0)-y_0)v(0)+\int_{\Sigma_T} y_d \nabla v \cdot \hat{n},$$
where $\hat{n}$ denotes the outer normal to the boundary $\partial\Omega$.
The weak formulation of equation \eqref{2ordp} for
given $y_d\in H^{2,1}(\Omega_T),\, y_0\in H^1_0(\Omega),$ reads:
\begin{equation}\label{weak_2ord}
\mbox{ find } p\in H_0^{2,1}(\Omega_T) \mbox{ with }A_0(p,v)=L_0(v)\quad \forall v\in H^{2,1}_0(\Omega_T).
\end{equation}

\noindent
It follows from the monotonicity of the orthogonal projection that 
\eqref{weak_2ord} admits a unique solution $p$, compare e.g. \cite[Th. 1.25]{HPUU09}.
We put our attention on the semi-discrete approximation of \eqref{2ordp}
and investigate a-priori and a-posteriori error estimates for the time discrete problem, 
where the space is kept continuous. Let us consider the time 
discretization $0=t_0<t_1<\ldots<t_n=T$ with $\Delta t_j=t_j-t_{j-1}$ 
and $\Delta t:=\max_j \Delta t_j$. Let $I_j:=[t_{j-1},t_j]$. 
We define the time discrete space
$$V_t^k:=\left\{v\in H^{2,1}(\Omega_T): \; v(\cdot)|_{I_j} \in P_1(I_j)\right\}, 
\qquad \bar{V}_t^k:=V_t^k\cap  H^{2,1}_0(\Omega_T),$$
where the notation $P_1(I_j)$ stands for the polynomials of first order on the interval 
$I_j$. Then, we consider the semi-discrete problem:
\begin{equation}\label{weak_dis}
\mbox{find }p_k\in \bar{V}_t^k \mbox{ with } A_0(p_k,v_k)=L_0(v_k),\quad \forall v_k\in \bar{V}_t^k.
\end{equation}
Using the arguments of e.g. \cite[Th. 1.25]{HPUU09} one can show that problem \eqref{weak_dis} 
admits a unique solution $p_k \in \bar{V}_t^k$.\\ 
We note that with \eqref{weak_2ord} and \eqref{weak_dis} we have the Galerkin orthogonality 
\begin{equation}\label{GO}
 A_0(p,v_k) - A_0(p_k,v_k) = 0 \quad \forall v_k \in \bar{V}_t^k. 
\end{equation}
Thus, for 
$v \in H_0^{2,1}(\Omega_T)$ it holds true
\begin{equation*}
 \begin{array}{r c l}
 A_0(p,v) - A_0(p_k,v) & = & A_0(p,v-v_k) - A_0(p_k,v-v_k) \quad 
 \forall v_k \in \bar{V}_t^k.\\
 \end{array}
\end{equation*}
The following Theorem states a temporal residual type a-posteriori error 
estimate for $p$, which transfers the estimation of \cite[Theorem 3.5]{GHZ12}
to the control constrained optimal control problem \eqref{ocp}:
\begin{teo}\label{thm:apost}
Let $p\in H_0^{2,1}(\Omega_T)$ and $p_k\in \bar{V}^k_t$ denote the 
solutions to \eqref{weak_2ord} and \eqref{weak_dis}, respectively. Then we obtain
\begin{equation}\label{est-thm31}
\|p-p_k\|_{H^{2,1}(\Omega_T)}^2\leq C_1\eta^2,
\end{equation}
where $C_1>0$ and
$$\eta^2=\sum_j \Delta t_j^2 \int_{I_j} \left\| -\dfrac{\partial y_d}{\partial t}+\Delta y_d+  \dfrac{\partial^2 p_k}{\partial t^2} +\mathcal{B}\mathcal{P}_{U_{ad}}\left(-\dfrac{1}{\alpha}\mathcal{B}^*p_k\right)-\Delta^2 p_k \right\|^2_{L^2(\Omega)} + \sum_j \int_{I_j} \|y_d-\Delta p_k\|_{L^2(\partial \Omega)}^2.$$
\end{teo}
\textit{Proof.} 
We start the proof showing a consequence of the monotonicity of the projector operator $-\mathcal{P}_{U_{ad}}\{-\mathcal{B}^*p\}$. We find that
$$\left\langle - \mathcal{P}_{U_{ad}}\left\{-\frac{1}{\alpha} \mathcal{B}^* p_1
\right\} + \mathcal{P}_{U_{ad}}\left\{-\frac{1}{\alpha} \mathcal{B}^* p_2
\right\}, \mathcal{B}^*p_1 - \mathcal{B}^*p_2 \right\rangle_U \geq 0,\quad \forall p_1,p_2\in H^{2,1}_0(\Omega_T),$$
and hence
\begin{equation}\label{monotonicity}
  \int_{\Omega_T} \left( - \mathcal{B} \mathcal{P}_{U_{ad}} \left\{ 
-\frac{1}{\alpha} \mathcal{B}^* p_1 \right\} + \mathcal{B}\mathcal{P}_{U_{ad}}
\left\{ -\frac{1}{\alpha} \mathcal{B}^* p_2 \right\} \right)(p_1 - p_2) 
\geq 0 .
\end{equation}

\noindent For easier notation, we set $N(p):= -\mathcal{B}\mathcal{P}_{U_{ad}} 
\left\{ -\frac{1}{\alpha} \mathcal{B}^*p \right\}$.\\
Let $e^p:= p - p_k$ and let $\pi_k e^p$ denote the standard Lagrange type 
temporal interpolation of $e^p$. Using the inequality
 $$ \| v \|_{H^{2,1}(\Omega_T)}^2 \leq C \left(\| \frac{\partial v}{\partial t} \|^2_{L^2(\Omega_T)} + \| \Delta v \|_{L^2(\Omega_T)}^2\right)$$
for $v \in H_0^{2,1}(\Omega_T)$ and $C > 0$ from \cite[Lemma 2.5]{GHZ12},
the monotonicity \eqref{monotonicity} and the Galerkin 
orthogonality \eqref{GO}, we 
can estimate:
\begin{equation*}
 \begin{array}{r c l}
 & &  c \| p - p_k \|_{H^{2,1}(\Omega_T)}^2  \\[2ex]
  & \leq &  \left\| \displaystyle\frac{\partial (p-p_k)}{\partial t} \right\|_{L^2(\Omega_T)}^2 + \| \Delta (p-p_k) \|_{L^2(\Omega_T)}^2 \\
 & \leq & \left\| \displaystyle\frac{\partial (p-p_k)}{\partial t} 
 \right\|_{L^2(\Omega_T)}^2 + \| \Delta (p-p_k) \|_{L^2(\Omega_T)}^2
 +  \displaystyle\int_{\Omega_T} (N(p)-N(p_k))(p-p_k) \\[2ex]
& = & \displaystyle\int_{\Omega_T} \displaystyle\frac{\partial (p-p_k)}{\partial t} \frac{\partial e^p}{\partial t} + \int_{\Omega_T} \Delta (p-p_k) \Delta e^p + \int_{\Omega_T}(N(p) - N(p_k))e^p \\[2ex]
 & = & \displaystyle\int_{\Omega_T} \displaystyle\frac{\partial (p-p_k)}{\partial t} \frac{\partial (e^p-\pi_k e^p)}{\partial t} + \int_{\Omega_T} \Delta (p-p_k) \Delta (e^p-\pi_k e^p) + \int_{\Omega_T}(N(p) - N(p_k))(e^p - \pi_k e^p) \\[2ex]
 & = & \displaystyle\int_{\Omega_T} (-\frac{\partial y_d}{\partial t} + \Delta y_d)(e^p - \pi_k e^p) + \int_{\Sigma_T} y_d \nabla(e^p - \pi_k e^p) \cdot \hat{n} - \int_{\Omega_T} \frac{\partial p_k}{\partial t} \frac{\partial (e^p - \pi_k e^p)}{\partial t}\\[2ex]
 & & - \displaystyle\int_{\Omega_T} \Delta p_k \Delta (e^p - \pi_k e^p)- \int_{\Omega_T}N(p_k)(e^p - \pi_k e^p)
 \end{array}
\end{equation*}
Integration by parts on each time interval and Green's formula lead to
\begin{equation*}
 \begin{array}{r c l}
 & &  c \| p - p_k \|_{H^{2,1}(\Omega_T)}^2  \\[2ex]
 & \leq & \displaystyle\sum_j \int_{I_j} \int_\Omega (- \frac{\partial y_d}{\partial t} + \Delta y_d + \frac{\partial^2 p_k}{\partial t^2} - \Delta^2 p_k - N(p_k))(e^p - \pi_k e^p) + \sum_j \int_{I_j} \int_{\partial \Omega} (y_d - \Delta p_k) \nabla (e^p - \pi_k e^p) \cdot \hat{n}.
 \end{array}
\end{equation*}
Utilizing error estimates of the Lagrange interpolation $\pi_k$, the trace 
inequality and Young's inequality, we find
\begin{equation*}
 \begin{array}{r c l}
 \hspace*{1cm} & &  \| p - p_k \|_{H^{2,1}(\Omega_T)}^2  \\[2ex]
 & \leq & C_1 \displaystyle\sum_j \Delta t_j^2 \int_{I_j} \left\| - 
 \frac{\partial y_d}{\partial t} + \Delta y_d + 
 \frac{\partial^2 p_k}{\partial t^2} - \Delta^2 p_k + \mathcal{B} 
 \mathcal{P}_{U_{ad}} \left\{-\frac{1}{\alpha} \mathcal{B}^* p_k \right\} \right\|_{L^2(\Omega)}^2 \\[2ex]
 & & + C_1 \displaystyle\sum_j \int_{I_j} \| y_d - \Delta p_k \|_{L^2(\partial \Omega)}^2. \hspace{9.5cm} \square
 \end{array}
\end{equation*}

Theorem \ref{thm:apost} provides a tool to refine the time grid by means of the residual of the system \eqref{2ordp}. Due to \eqref{Proj}, the time instances of 
this grid may be regarded as ideal snapshot locations for POD-MOR applied to problem \eqref{ocp}.



\section{POD for optimal control problems}

In this section, we recall the POD method which we use in order to replace the 
original problem \eqref{ocp} by a surrogate model. The main interest when applying the POD method is to reduce computation times and storage capacity while retaining a 
satisfying approximation quality. This is possible due to the key fact that 
POD basis functions (unlike typical 
finite element ansatz functions) contain information about the 
underlying model, since the POD modes are derived from snapshots of 
a solution data set. For this reason it is important to use rich snapshot ensembles 
reflecting 
the dynamics of the modeled system. Usually, we are able to improve 
the accuracy of a POD suboptimal solution by enlarging the number of utilized POD basis functions or enriching the 
snapshot ensemble, for instance. The snapshot form of POD proposed 
by Sirovich in \cite{Sir87} works in the continuous version as follows.\\
\noindent
Let us suppose that the continuous solution $y(t)$ of \eqref{heat} and $p(t)$ of \eqref{adj} 
belongs to a real 
separable Hilbert space $V$, where $V=H_0^1(\Omega)$ or $L^2(\Omega)$, 
equipped with its inner product $\langle\cdot,\cdot\rangle$ and associated 
norm $\|\cdot\|^2=\langle\cdot,\cdot\rangle$. We set   
\mbox{$\mathcal{V}:=\mbox{span}\{z^k(t) \; | \; t \in [0,T] \text{ and } 1 \leq k \leq 3 \} \subseteq V$}, 
where $z^1(t) := y(t)$, $z^2(t) := p(t)$, $z^3(t) := \dot{p}(t)$. Note that the initial condition 
$y(0) = y_0$ is included in $\mathcal{V}$. The aim is to determine 
a POD basis $\{\psi_1,\ldots,\psi_\ell\} \subset V$ of 
rank $\ell \in \{1, ..., d\}$ with $d= \text{dim}(\mathcal{V}) \leq \infty$, by solving 
the following constrained minimization problem:
\begin{eqnarray}\label{prb:pod}
\min_{\psi_1,\ldots,\psi_\ell} \sum_{k=1}^3\int_0^T \left\|z^k(t)-\sum_{i=1}^\ell \langle z^k(t),\psi_i\rangle \; \psi_i\right\|^2 dt
\quad \mbox{ s.t. } \langle\psi_j,\psi_i\rangle=\delta_{ij}\quad\mbox{for }1\leq i,j\leq \ell,
\end{eqnarray}
where $\delta_{ij}$ denotes the Kronecker symbol, i.e. $\delta_{ij}=0$ for $i \neq j$ 
and $\delta_{ii} = 1$. \\
\noindent
It is well-known (see \cite{GV13}) that a solution to problem 
(\ref{prb:pod}) is given by the first $\ell$ eigenvectors $\{\psi_1, \ldots, 
\psi_\ell\}$ corresponding to the $\ell$ largest eigenvalues
$\lambda_i > 0$ of the self-adjoint linear operator 
$\mathcal R:V \rightarrow V,$ i.e. \mbox{$\mathcal{R}\psi_i=\lambda_i\psi_i$,} $i=1, \dotsc, \ell$, 
where $\mathcal R$ is defined as follows:
$$ \mathcal{R}\psi=\sum_{k=1}^3 \int_0^T \langle z^k(t),\psi\rangle \; z^k(t) dt \quad \mbox{for } \psi\in V.$$ 
\noindent Moreover, we can quantify the POD approximation error by the neglected 
eigenvalues (more details in \cite{GV13}) as follows:
\begin{equation}\label{err-POD}
\sum_{k=1}^3 \int_0^T \left\|z^k(t)-\sum_{i=1}^\ell \langle z^k(t),\psi_i\rangle \; 
 \psi_i\right\|^2 dt =\sum_{i=\ell+1}^d \lambda_i.
 \end{equation}

\noindent
Let us assume that we have computed POD basis functions $\{ \psi_i\}_{i=1}^\ell$.
Then, we define the POD Galerkin ansatz of order $\ell$ for the state $y$ as:
\begin{equation}\label{pod_ans}
y^\ell(t)=\sum_{i=1}^\ell w_i(t)\psi_i,
\end{equation}

\noindent
where $y^\ell\in V^\ell:= \text{span}\{\psi_1,\ldots,\psi_\ell\}$ and the 
unknown coefficients are denoted by $\{w_i\}_{i=1}^\ell$.
If we plug this ansatz into the weak formulation of the state 
equation \eqref{weak:heat} and use $V^\ell$ as the test space, we get the 
following reduced order model for \eqref{weak:heat} of low 
dimension:

\begin{equation}\label{weakpod:heat}
\begin{array}{r c l}
\displaystyle\int_\Omega y_t^\ell(t) \psi dx + \int_\Omega \nabla y^\ell(t) \cdot \nabla \psi dx & = & \displaystyle\int_\Omega (f + \mathcal{B} u)(t) \psi dx 
\quad \forall \psi \in V^\ell \text{ and } t \in (0,T] \text{ a.e.},\\[2ex]
\displaystyle\int_\Omega y^\ell(0) \psi dx & = & \displaystyle\int_\Omega y_0 \psi dx\\[2ex]
\end{array}
\end{equation}

\noindent Choosing $\psi = \psi_i$ 
for $i=1, \dotsc, \ell$ and utilizing \eqref{pod_ans}, we infer from \eqref{weakpod:heat} 
that the coefficients \mbox{$(w_1(t), \dotsc, w_\ell(t)) =: w(t)$} satisfy
$$M^\ell \dot{w}(t) + A^\ell w(t) = F^\ell(t) \quad \text{a.e. in} (0,T], \quad M^\ell w(0) = y_0^\ell, $$
where $(M^\ell)_{ij}  = \int_\Omega \psi_j \psi_i dx$, 
$(A^\ell)_{ij}  = \int_\Omega \nabla \psi_j \cdot 
\nabla \psi_i dx$, $(F^\ell (t))_j = \int_\Omega (f + \mathcal{B}u)(t) \psi_j dx$ and 
$(y_0^\ell)_j = \int_\Omega y_0 \psi_j dx$. Note that $M^\ell$ is the 
identity matrix, if we choose as inner product $\langle \cdot , \cdot \rangle :=
\langle \cdot , \cdot \rangle_{L^2(\Omega)}$.
\noindent
The reduced order model surrogate (ROM) for the optimal control
problem is given by
\begin{equation}\label{ocp_pod}
\min_{u \in U_{ad}} \hat{J}^\ell(u) \mbox{ s.t. } y^\ell(u) \mbox{ satisfies } \eqref{weakpod:heat},
\end{equation}
where $\hat{J}^\ell$ is the reduced cost functional, i.e. $\hat{J}^\ell 
(u):= \hat{J} (y^\ell(u),u) $. We recall that the discretization of the optimal solution $\bar{u}^\ell$ to \eqref{ocp_pod} is determined 
by the relation between the adjoint state and control and refer 
to \cite{H05} for more details about the variational discretization concept.\\
\noindent In order to solve the reduced optimal control problem \eqref{ocp_pod}, we consider 
the well-known first order optimality condition given 
by the variational inequality
$$\langle \nabla \hat{J}^\ell(\bar{u}^\ell), 
u-\bar{u}^\ell \rangle_U \geq 0 \quad \forall u \in U_{ad},$$ 
which is sufficient since the underlying problem is convex.\\
The first order optimality conditions of \eqref{ocp_pod} also deliver that the adjoint POD scheme 
for the approximation of $p$ is given by: find 
$p^\ell(t) \in V^\ell$ with $p^\ell(T) = 0$ satisfying
\begin{equation}\label{weakpod:adjoint}
-\displaystyle\int_\Omega p_t^\ell(t) \psi dx + \int_\Omega \nabla p^\ell(t) \cdot \nabla \psi dx  =  \displaystyle\int_\Omega (y^\ell - y_d)(t) \psi dx 
\quad \forall \psi \in V^\ell \text{ and } t \in (0,T) \text{ a.e.}
\end{equation}

\section{The snapshot location strategy}

In Section 4, the POD method in the continuous framework is recalled, where the 
POD basis functions are computed in such a way that the error between the trajectories $y(t)$ of 
\eqref{heat} and $p(t)$ of \eqref{adj} and its POD Galerkin approximation is minimized in \eqref{prb:pod}. In practice, we 
do not have the whole solution trajectories $\{z^k(t)\}_{t\in[0,T]}, 1 \leq k \leq 4$ at hand. 
But we have snapshots available, which are the solutions $\{y(t_j)\}_{j=0}^n$ to \eqref{heat} and 
the solutions $\{p(t_j)\}_{j=0}^n$ to \eqref{adj} at times $\{t_j\}_{j=0}^n$. 
This motivates to replace the time integration in \eqref{prb:pod} by an appropriate quadrature 
rule based on $t_0,\ldots,t_n$, i.e. $\int_0^T g(t)dt\approx \sum_{j=0}^n \alpha_j g(t_j)$ for 
$g\in C^0([0,T])$ with quadrature weights $\beta_0,\ldots,\beta_n\in\R.$ We later choose the 
weights for the trapezoidal rule, compare \eqref{post_p}. In the present work, we neglect the error 
introduced by quadrature weights.

%
The minimization problem related to \eqref{prb:pod} then becomes
$$\min \; \sum_{k=1}^3\sum_{j=0}^n \beta_j \left\| z^k(t_j) - \sum_{i=1}^\ell \langle z^k(t_j),\psi_i \rangle \psi_i 
\right\|^2, \text{ s.t. } \langle \psi_j , \psi_i \rangle = \delta_{ij} \quad \text{ for } 1 \leq i,j \leq \ell$$
 \noindent and obviously constitutes a strong dependence of the POD basis functions 
 on the chosen snapshot locations $t_0,\ldots, t_n$. The related snapshots shall have the 
 property to capture the main features of the dynamics of the truth solution
 as good as possible. Here it is important to select suitable time instances at which 
 characteristic dynamical properties of the optimal state are located. A natural question is:
 
 \medskip

  \begin{center}
  {\em How to pick time instances that represent 
 good locations for snapshots in POD-MOR for \eqref{ocp_pod}?}\\
  \end{center}
 
 \medskip
 
Moreover, we face some difficulties since the reduction of optimal control 
 problems is usually initialized with snapshots computed from a 
 given input control $u_\circ \in U_{ad}$. This problem is usually addressed in the 
 offline stage for POD, which is the phase needed for snapshot generation, POD basis 
 computation and building the reduced order model. Mostly, we 
 do not have any information about the optimal control, such that 
 in POD-MOR the input control $u_\circ$ is often chosen as $u_\circ \equiv 0$. This circumstance raises the 
 question about the quality of the POD basis and the quality of the POD suboptimal solution.
The a-posteriori error estimator \eqref{est-thm31} in Section 3 motivates a suitable location 
of time instances for the POD adjoint
state and at the same time we get an approximation of the optimal control which can 
be used as an input control $u_\circ$ in order to generate the snapshots. 

\noindent The use, in the offline-stage, of a time adaptive mesh refinement process allows to overcome the choice of 
an input control $u_\circ$ and the choice of the snapshot locations by 
solving equation \eqref{2ordp}. Then, we take advantage of the a-posteriori error 
estimation presented in Theorem \ref{thm:apost}. Equation 
\eqref{2ordp} provides the optimal adjoint state associated with \eqref{ocp}, which
does not require the explicit knowledge of a control input $u_\circ$. We note that the ellipticity of equation \eqref{2ordp} play a crucial role in this approach. The same approach would not work, if one solves the 
optimality conditions directly. The numerical approximation of $p$ provides important information 
about the control input. In fact, thanks to the variational inequality \eqref{opt_con} we are 
first able to build an approximate control $u$ and finally compute the associated state 
$y(u)$. In this way our snapshot set will contain information about the state corresponding to an 
approximation of the optimal control. Thanks to this numerical approximation 
of the optimal control problem we can build the snapshot matrix and compute the POD basis 
functions where the number $\ell$ is chosen such that $\sum_{i=\ell+1}^d \lambda_i\approx0$.\\ 
The approximation of equation \eqref{2ordp} is very useful in model order reduction since we overcome the choice of the initial input control to generate the snapshot set. Moreover, we also gain information about a temporal grid, which allows us to better 
resolve $p$ with respect to time.
\noindent
The a-posteriori error estimation \eqref{est-thm31}
guarantees that the finite element approximation of \eqref{2ordp} in the time variable is below a certain tolerance. Therefore, the 
reduced optimal control problem \eqref{ocp_pod} is set up and solved on the resulting adaptive time grid. Now the 
question is:

\medskip

  \begin{center}
\textit{How good is the quality of the computed time grid in terms of the error between\\ the optimal solution
and the POD surrogate solution? } \\
  \end{center}

\medskip

\subsection{Error Analysis for the adjoint variable}
Let us motivate our approach by analyzing the error 
$\|p(u)-p_k^\ell(u_k^\ell)\|_{L^2(0,T,V)}$  between the optimal adjoint solution $p(u)$
of \eqref{adj} associated with the optimal control $u$ for 
\eqref{ocp}, i.e. $u=\mathcal{P}_{U_{ad}}(-\dfrac{1}{\alpha}\mathcal{B}^*p)$ 
and the POD reduced approximation $p_k^\ell(u_k^\ell)$, which is 
the time discrete solution to the POD-ROM for \eqref{adj} associated with the 
time discrete optimal control $u_k^\ell$ for \eqref{ocp_pod}, i.e. $y=y(u_k^\ell)$ in \eqref{adj}. We denote by $V$ the space 
$V=H_0^1(\Omega)$ and by $H$ the space $L^2(\Omega)$. By the triangular inequality we get the following estimates for the $L^2(0,T;V)$-norm:
\begin{align}\label{err:est_new}
\|p(u)-p_k^\ell(u_k^\ell)\| &\leq \underbrace{\|p(u)-p_k(u_k)\|}_{(\ref{err:est_new}.1)}+\underbrace{\|p_k(u_k)- \mathcal{P}^\ell p_k(u_k)\|}_{(\ref{err:est_new}.2)}+ \underbrace{\| \mathcal{P}^\ell p_k(u_k) - \mathcal{P}^\ell p_k(u_k^\ell) \|}_{(\ref{err:est_new}.3)} + \underbrace{\| \mathcal{P}^\ell p_k(u_k^\ell) - p_k^\ell(u_k^\ell) \|}_{(\ref{err:est_new}.4)}
\end{align}
where $p_k(u_k)$ is the time discrete adjoint solution of \eqref{weak_dis} associated with 
the control $u_k=\mathcal{P}_{U_{ad}}(-\dfrac{1}{\alpha}\mathcal{B}^*p_k)$ and 
$p_k(u_k^\ell)$ is the time discrete adjoint solution to \eqref{adj} with respect to the 
suboptimal control $u_k^\ell$, i.e. $y=y(u_k^\ell)$ in \eqref{adj}. By $\mathcal{P}^\ell:V\rightarrow V^\ell$ we 
denote the orthogonal POD projection operator as follows:  
$$\mathcal{P}^\ell y:=\sum_{i=1}^\ell \langle y,\psi_i\rangle_V \psi_i\quad \mbox{ for } y\in V.$$
The term (\ref{err:est_new}.1) can be estimated by \eqref{est-thm31} and concerns 
the snapshot generation. Thus, we can 
decide on a certain tolerance in order to have a prescribed error. The second term 
(\ref{err:est_new}.2) in \eqref{err:est_new} is the POD projection error and can be estimated 
by the sum of the neglected eigenvalues. Then, we note that the third term 
(\ref{err:est_new}.3) can be estimated as follows: 
\begin{equation}\label{err:est_third}
 \| \mathcal{P}^\ell p_k(u_k) - \mathcal{P}^\ell p_k(u_k^\ell) \| \leq \| \mathcal{P}^\ell \| \,  \|p_k(u_k) - p_k(u_k^\ell) \| \leq C_2 \|u_k-u_k^\ell \|_U,
 \end{equation}
where $\| \mathcal{P}^\ell \| \leq 1$ and $C_2>0$ is the constant referring to the Lipschitz continuity of $p_k$ independent of $k$ as in \cite{NV11}.

In order to control the quantity $\| u_k-u_k^\ell \|_U \leq \| u_k - u \|_U + \| u - u_k^\ell \|_U$ we make use of the 
a-posteriori error estimation of \cite{TV09}, which provides an upper bound for the error 
between the (unknown) optimal control and any arbitrary control $u_p$ (here $u_p = u_k$ and 
$u_p = u_k^\ell$) by 
$$ \| u - u_p \|_U \leq \frac{1}{\alpha} \| \zeta_p \|_U, $$ where  $\alpha$ is the 
regularization parameter in the cost functional and $\zeta_p \in 
L^2(0,T;\mathbb{R}^m)$ is chosen such that 
$$ \langle \alpha u_p - \mathcal{B}^* p(u_p) + \zeta_p, u - u_p
\rangle_U  \geq 0 \quad \forall u \in U_{ad} $$
is satisfied. Finally, last term (\ref{err:est_new}.4) can be estimated according to \cite{HV08} and involves 
the sum of the eigenvalues not considered, the first derivative of the time discrete adjoint 
variable and the difference between the state and the POD state:
\begin{equation}\label{err:est_fourth}
 \| \mathcal{P}^\ell p_k(u_k^\ell) - p_k^\ell(u_k^\ell) \|^2  \leq C_3 \left(\sum_{i=\ell+1}^d \lambda_i^k 
+  \| \dot{p}_k(u_k^\ell)-\mathcal{P}^\ell \dot{p}_k(u_k^\ell)  \|_{L^2(0,T,V')}^2 + 
\|y_k(u_k^\ell) - y_k^\ell(u_k^\ell) \|_{L^2(0,T,H)}^2\right),
\end{equation}
for a constant $C_3>0$.
We note that the sum of the neglected eigenvalues is sufficiently small provided that 
$\ell$ is large enough. Furthermore, error estimation \eqref{err:est_fourth} depends on the 
time derivative $\dot{p}_k$. To avoid this dependence, we include time derivative 
information concerning the adjoint variable into the snapshot set, see \cite{KV02}. 


To summarize, the error estimation reads:

\begin{equation}\label{err:summarized}
 \| p(u) - p_k^\ell(u_k^\ell) \|_{L^2(0,T,V)}  \leq \sqrt{C}_1 \eta + 
 \frac{C_2}{\alpha}( \|\zeta_k \|_U + \|\zeta_k^\ell \|_U )
+ \sqrt{  C_3 \left( \sum_{i=\ell+1}^d \lambda_i^k +
  \|y_k - y_k^\ell \|_{L^2(0,T,H)}^2 \right)}.
\end{equation}

Finally, we note that estimation \eqref{err:est_fourth} 
involves the state variable which is estimated in the following Section 5.2.

\subsection{Error Analysis for the state variable}

In this section we address the problem of the certification of the quality for POD approximation 
for the state variable.
It may happen that the time grid selected for the adjoint $p$ will not be accurate enough for 
the state variable $y$. Therefore a further refinement of the time grid might be 
useful in order to reduce the error between the POD state and the true state below
a given threshold. This is not guarenteed if we use the time grid, which 
results from the use of the estimate \eqref{est-thm31}. 
Here, we consider the error between the full solution $y(u_k^\ell)$ corresponding to the suboptimal control $u_k^\ell$ and the time discrete POD 
solution $y_k^\ell(u_k^\ell)$, where we assume to have the same temporal grid for snapshots and 
the solution of our POD reduced order problem. In this situation, the following estimate is proved in \cite{KV02}:
%

\begin{subequations}\label{post_p}
\begin{align}
\displaystyle\sum_{j=0}^n \beta_j \|y(t_j;u_k^\ell)-y_j^\ell(u_k^\ell) \|_H^2 &\leq \quad \displaystyle\sum_{j=1}^n 
\left( \Delta t_j^2 C_y ((1+c_p^2)\|y_{tt}(u_k^\ell)\|^2_{L^2(I_j,H)}+\|y_t(u_k^\ell)\|_{L^2(I_j;V)})\right) \label{post_p_a}\\[0.4cm]
& \hspace{1cm} +\displaystyle\sum_{j=1}^n C_y\left(\sum_{i=\ell+1}^d\left(|\langle \psi_i,y_0\rangle_V|^2+\lambda_i\right)\right) \label{post_p_b}\\[0.4cm]
& \hspace{1cm} +\displaystyle\sum_{j=1}^n \sum_{i=\ell+1}^d C_y \dfrac{\lambda_i}{\Delta{t_j^2}} \label{post_p_c}
\end{align}
\end{subequations}


\noindent where $C_y>0$ is a constant depending on $T$, but independent of the time grid 
$\{t_j\}_{j=0}^n$. We note that $y(t_j;u_k^\ell)$ is the continuous solution of \eqref{heat} at given time instances related to the suboptimal control $u_k^\ell$. The 
temporal step size in the subinterval $[t_{j-1},t_j]$ is denoted by $\Delta t_j$. The positive 
weights $\beta_j$ are given by
$$ \beta_0 = \frac{\Delta t_1}{2}, \quad \beta_j = \frac{\Delta t_j + 
\Delta t_{j+1}}{2} \text{ for } j = 1, \dotsc, n-1, \quad \text{and } \beta_n = \frac{\Delta t_n}{2}. $$
The constant $c_p$ is an upper bound of the projection 
operator. A similar estimate can be carried out for the $V-$norm. We refer the interested reader to \cite{KV02}.

Estimate \eqref{post_p} provides now a recipe for further refinement of the time 
grid in order to approximate the state $y$ within a prescribed tolerance. One option here 
consists in equidistributing the error contributions of the term \eqref{post_p_a}, while the 
number of modes has to be adapted to the time grid size according to the term \eqref{post_p_c}. 
Finally, the number $\ell$ of modes should be chosen such that the term in \eqref{post_p_b} 
remains within the prescribed tolerance.

\subsection{The algorithm}

The a-posteriori error control concept for \eqref{2ordp} now offers the possibility to select 
snapshot locations by a time adaptive procedure. For this purpose, \eqref{2ordp} is solved 
adaptively in time, where the spatial resolution ($\Delta x$ in Algorithm \ref{Alg:OPTPOD}) is 
chosen to be very coarse in order to keep the computational costs low. This is possible due to the fact that spatial and temporal 
discretization decouple when using the solution technique of \cite{GHZ12} as we will see 
in Section 6, compare Figure \ref{fig3:mesh}.
The resulting time grid points now serve as snapshot locations, on which our POD reduced order model 
for the optimization is based. 
The snapshots are now obtained from a simulation of \eqref{heat} with high spatial resolution h, which is used in \eqref{adj} to obtain highly resolved snapshots of p, which are accomplished with time finite differences of those adjoint snapshots. The right-hand side $u$ in the simulation of \eqref{heat} is obtained from \eqref{opt_con} with $p$ from \eqref{adj} computed with spatially coarse resolution $\Delta x$. 
The certification of the state variable is then performed according to \eqref{post_p} as a post-processing procedure.
This strategy might not deliver the optimal time instances, but it is a practical and efficient 
strategy, which turns out to deliver good approximation results (compare Section 6) at low 
costs. 

The algorithm is summarized below in Algorithm \ref{Alg:OPTPOD}. 

\begin{algorithm}[htbp]
\caption{Adaptive snapshot selection for optimal control problems.}
\label{Alg:OPTPOD}
\begin{algorithmic}[1]
\REQUIRE coarse spatial grid size $\Delta x$, fine spatial grid size $h$, maximal number of degrees of freedom (dof)\\ for the adaptive time discretization, $T>0$.
\STATE Solve \eqref{2ordp} adaptively w.r.t. time with spatial resolution $\Delta x$ and obtain the time grid $\mathcal{T}$ with solution $p_{\Delta x}$. 
\STATE Set $u_{\Delta x}=\mathcal{P}_{U_{ad}}\left(-\dfrac{1}{\alpha}\mathcal{B}^*p_{\Delta x}\right).$
\STATE Solve \eqref{heat} on $\mathcal{T}$ with spatial resolution $\Delta x$ corresponding to the control $u_{\Delta x}$.
\STATE Refine the time interval $\mathcal{T}$ according to \eqref{post_p} and construct the time grid $\mathcal{T}_{new}$.
\STATE Generate state and adjoint snapshots by solving \eqref{heat} with r.h.s. $u_{\Delta x}$ and \eqref{adj}, respectively, on $\mathcal{T}_{new}$ with spatial resolution $h$. Generate time derivative adjoint snapshots with time finite differences on those adjoint snapshots. 
\STATE Compute a POD basis of order $\ell$ and build the POD reduced order model \eqref{ocp_pod} based on the state, \\adjoint state and time derivative adjoint state snapshots.
\STATE Solve \eqref{ocp_pod} with the time grid $\mathcal{T}_{new}$ 

%
%
\end{algorithmic}
\end{algorithm}

%

\section{Numerical Tests}

In our numerical computations we use a one-dimensional spatial domain and a 
finite element discretization in space by means of conformal piecewise linear 
polynomials. We use the implicit 
Euler method for time integration. The solution of the optimal control problem \eqref{ocp_pod}
is done by a gradient method with stopping criteria $\|\hat{J}'(u^k)\| \leq \tau_r \|\hat{J}'(u^k)\|_U +\tau_a$ and an Armijo linesearch. In the following numerical examples, we apply Algorithm \ref{Alg:OPTPOD} in 
order to validate this strategy by numerical results.\\ 
The numerical tests illustrate that utilizing a time adaptive grid for snapshot location and for solving the
POD reduced order optimal control problem delivers more accurate approximation results than utilizing 
a uniform time grid. We show three different numerical tests.
The first example presents a steep gradient at the end of the time interval in the adjoint variable. 
 In the second example the adjoint state develops an interior layer in the middle of the 
 time interval and finally we introduce control contraints in the third example. Moreover we 
 also show the benefits of the post processing for the state variable 
 (step 4 in Algorithm \ref{Alg:OPTPOD}) to achieve more accurate approximation results for both state and adjoint state. \\
 All coding is done in \textsc{Matlab R2015}a and the computations are performed on a 2.50GHz computer.\\

%

\subsection{Test 1: Solution with steep gradient towards final time}

\begin{figure}[htbp]
\includegraphics[scale=0.33]{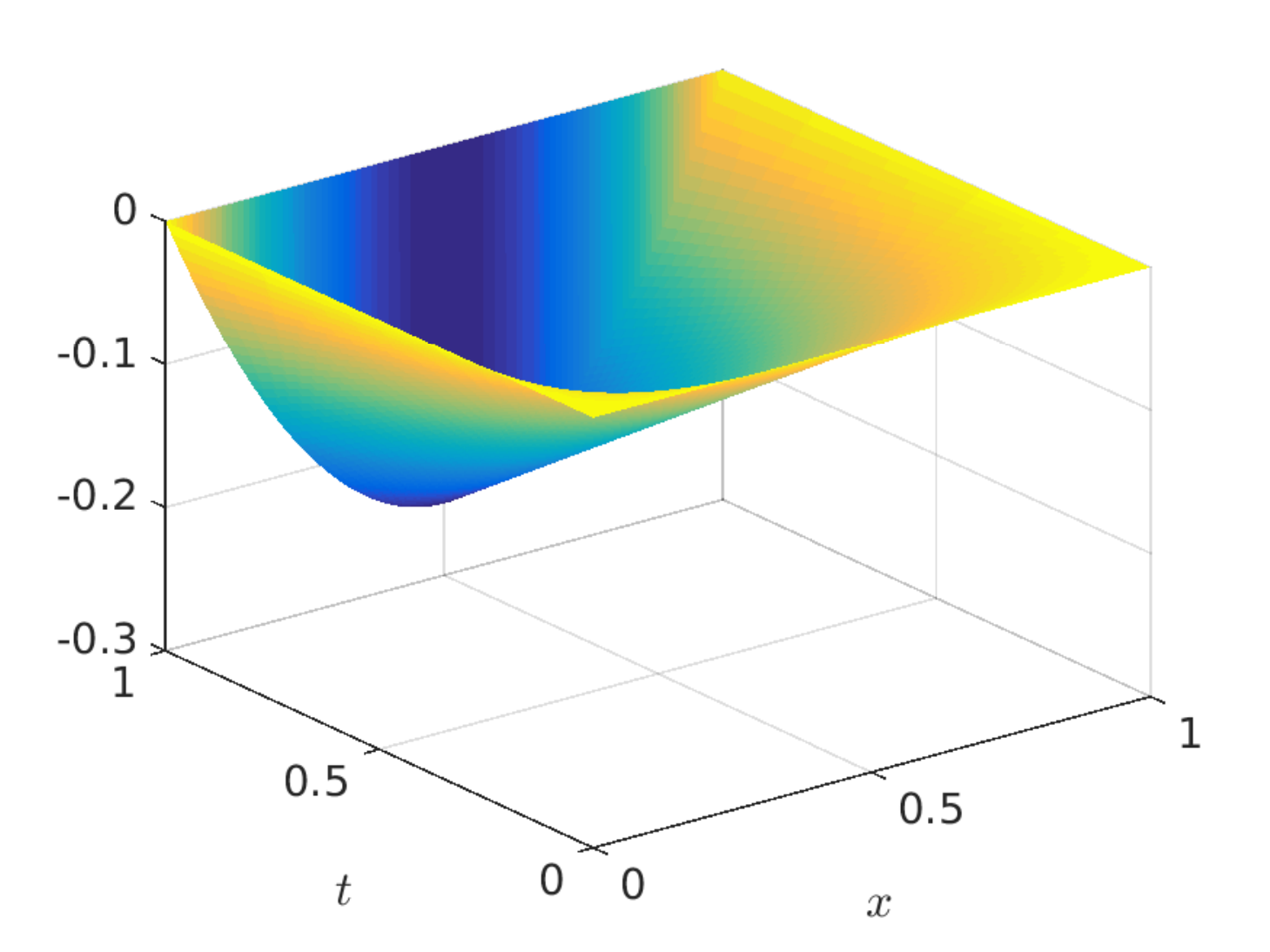} \includegraphics[scale=0.33]{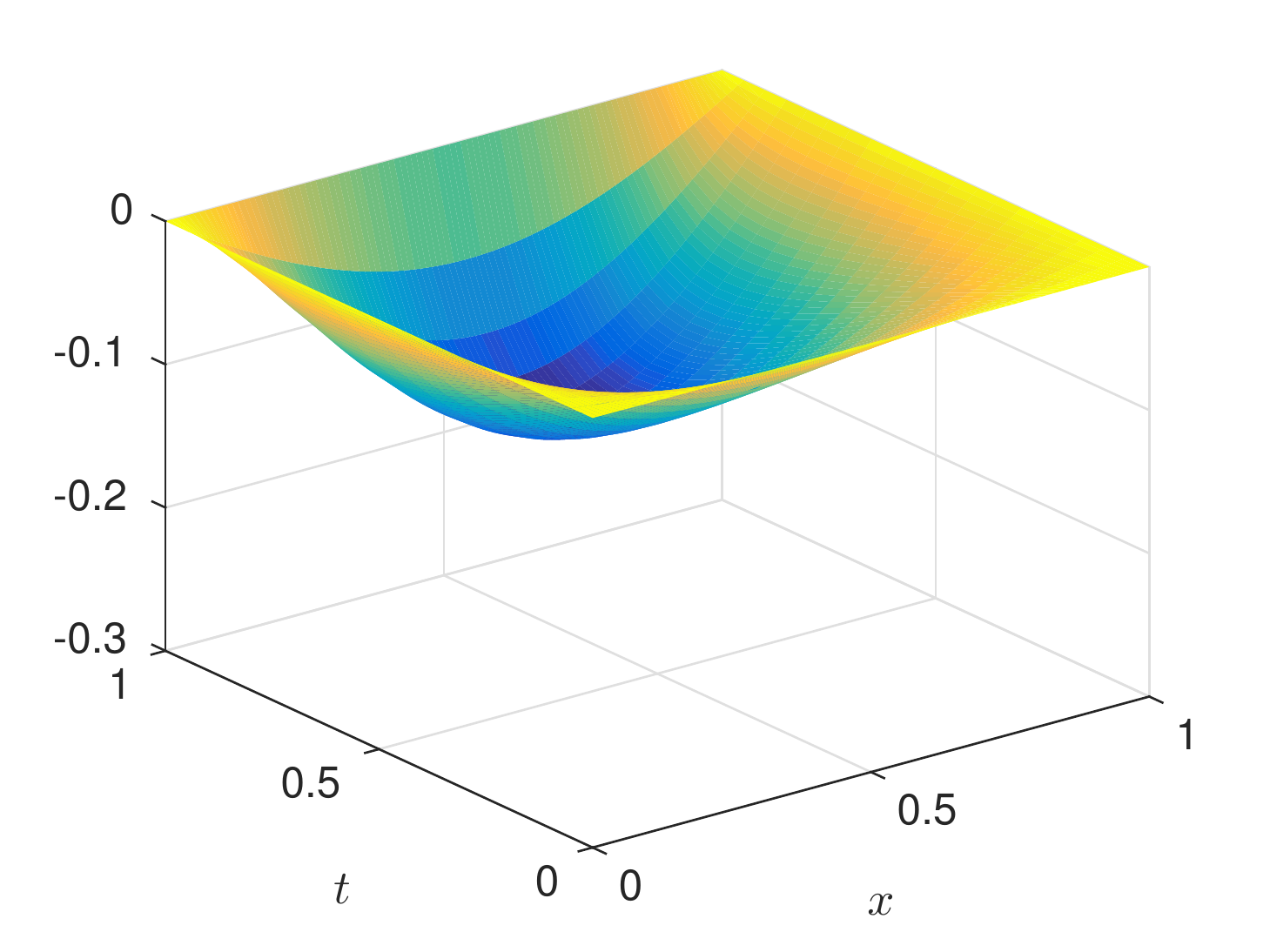} \includegraphics[scale=0.33]{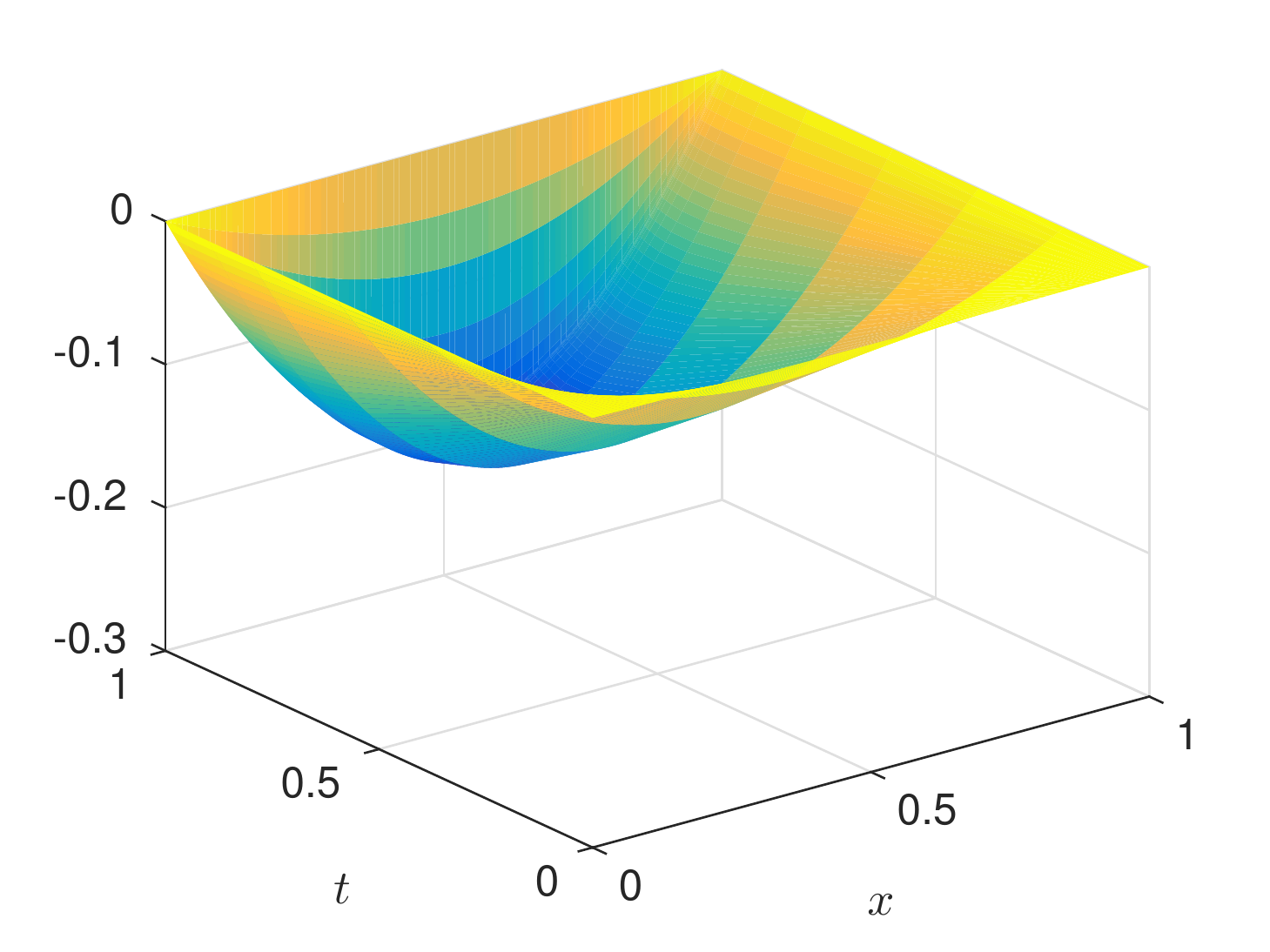}
 \caption{Test 1: Analytical optimal adjoint state $\bar{p}$ (left), POD adjoint 
 solution $p^\ell$ utilizing an equidistant time grid with $ \Delta t = 1/20$ (middle), 
 POD adjoint solution $p^\ell$ utilizing an adaptive time grid with dof=21 (right).}
 \label{fig1:sol} 
 \end{figure}

\begin{figure}[htbp]
 \includegraphics[scale=0.33]{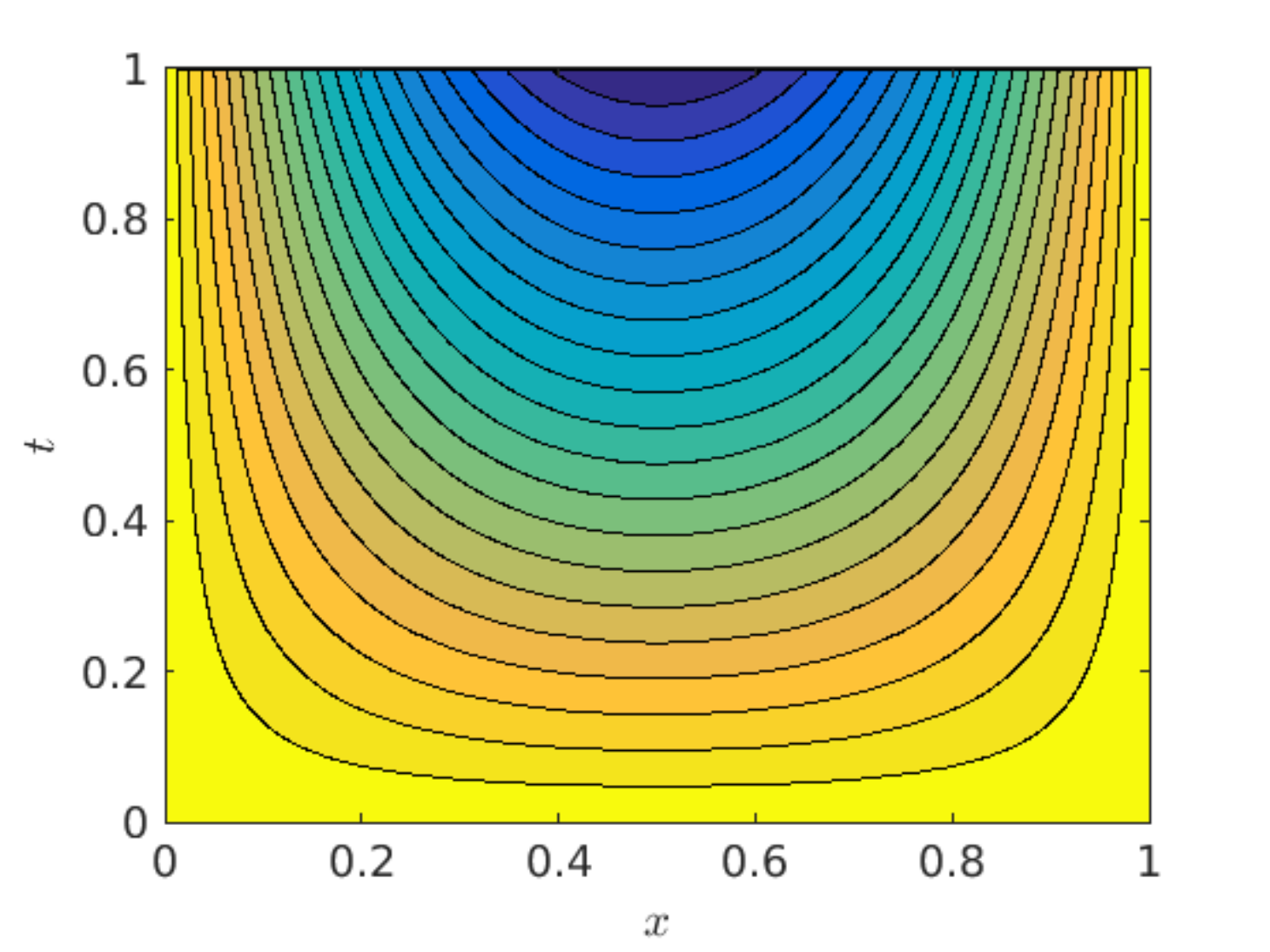}  \includegraphics[scale=0.33]{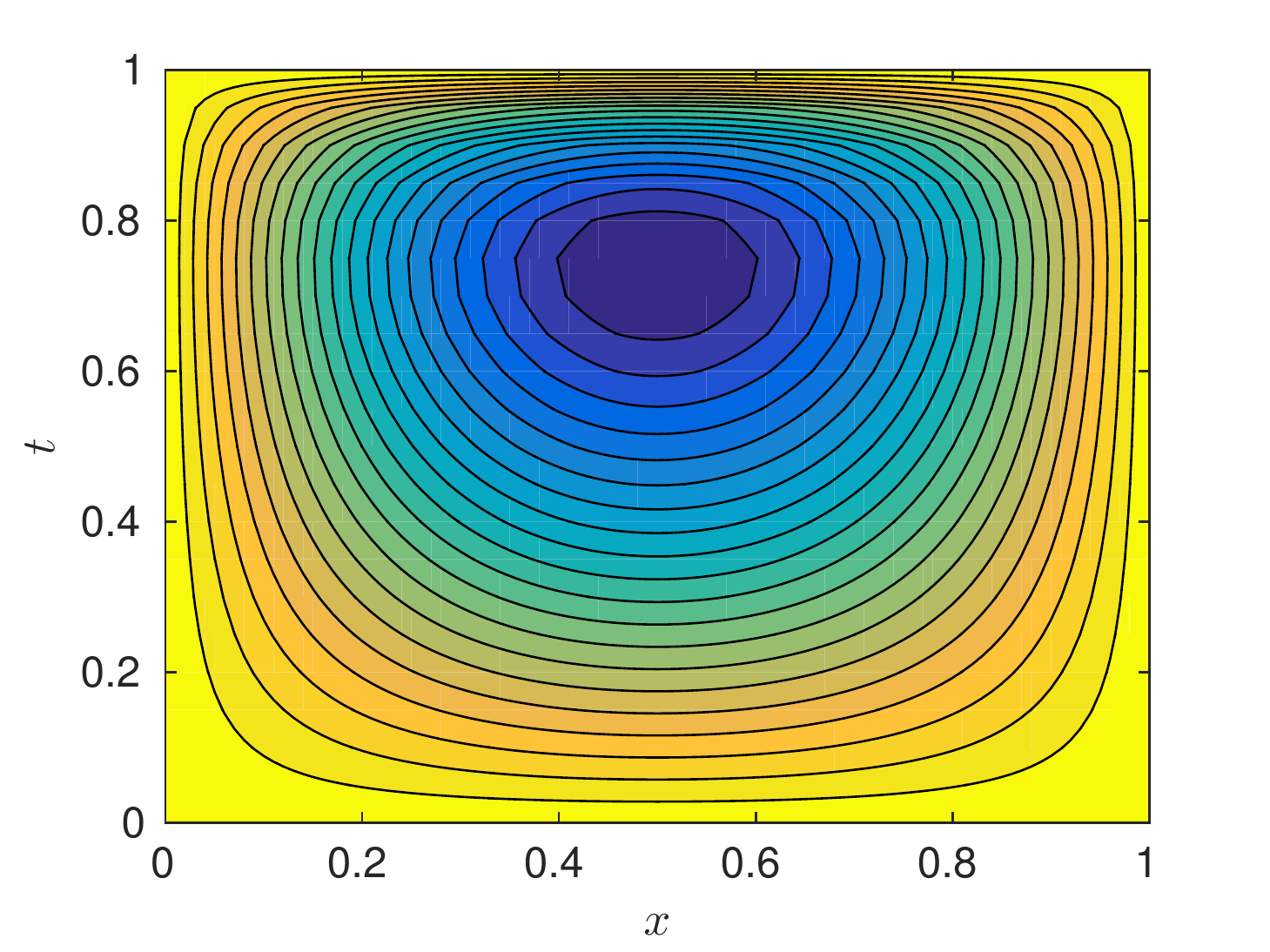}  \includegraphics[scale=0.33]{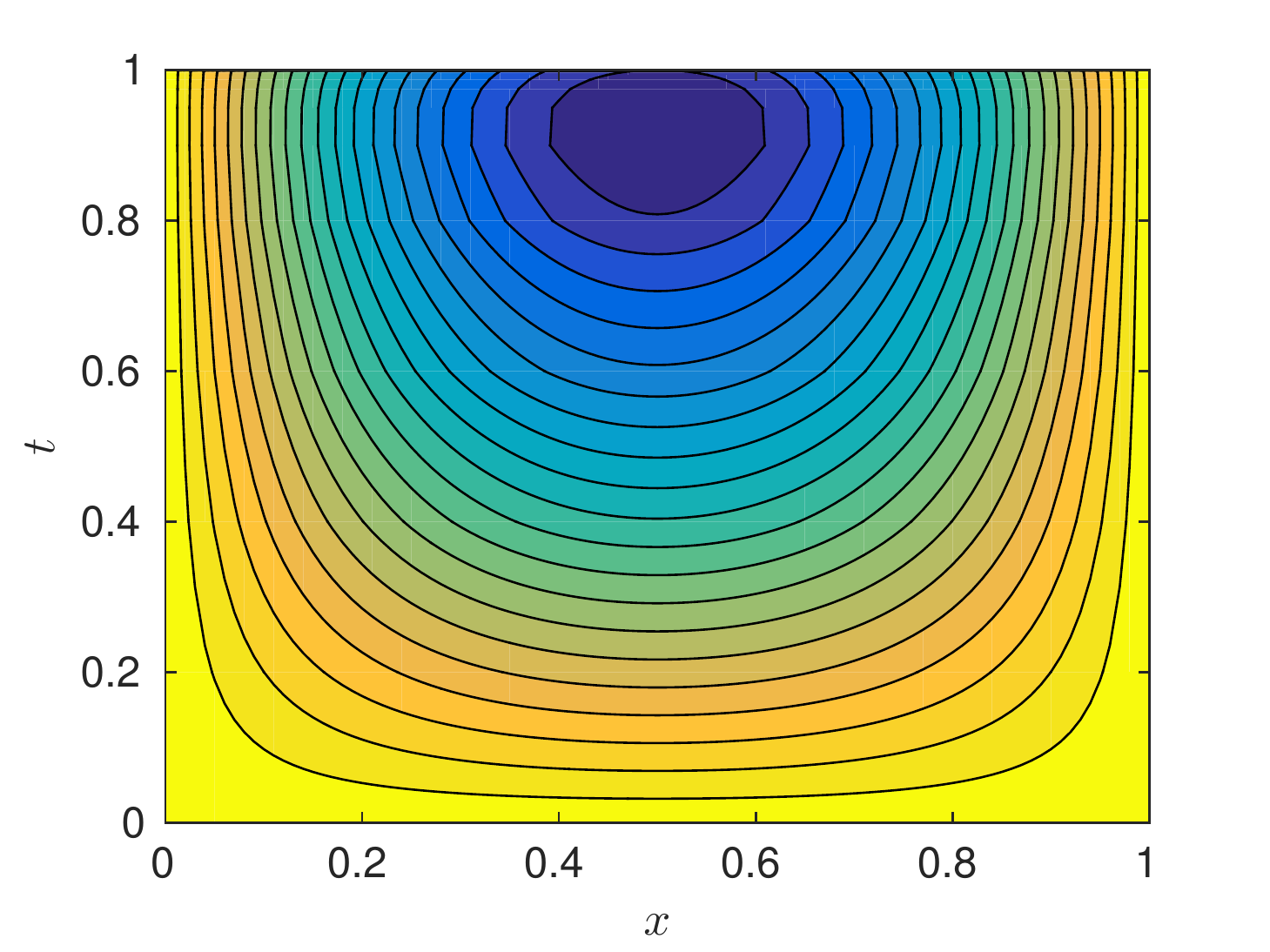}
 \caption{Test 1: Contour lines of the analytical optimal adjoint state 
 $\bar{p}$ (left), POD adjoint solution $p^\ell$ utilizing an equidistant 
 time grid with $ \Delta t = 1/20$ (middle), POD adjoint solution $p^\ell$ utilizing an adaptive 
 time grid with dof=21 (right).}
 \label{fig2:con}
 \end{figure}

The data for this test example is inspired from Example 5.3 in \cite{GHZ12}, with the 
following choices: $\Omega = (0,1)$ and 
$[0,T] = [0,1]$. We set $U_{ad} = L^\infty(0,T;
\mathbb{R}^m)$. The example is built in such a way that the exact optimal 
solution $(\bar{y},\bar{u})$ of problem \eqref{ocp} with associated optimal 
adjoint state $\bar{p}$ is known: $$\bar{y}(x,t) = \sin (\pi x) \sin(\pi t), \quad 
\bar{p}(x,t) = x(x-1)\left(t-\frac{e^{(t-1)/\varepsilon}-
e^{-1/\varepsilon}}{1-e^{-1/\varepsilon}}\right), \quad
\bar{u}(t) = -\frac{1}{\alpha} \mathcal{B}^*\bar{p}(x,t) = 
-t+\frac{e^{(t-1)/\varepsilon}-e^{-1/\varepsilon}}{1-e^{-1/\varepsilon}}$$ 
with $m=1$ and the control shape function $\chi(x) = x(x-1)$ for the operator $\mathcal{B}$. This leads to 
the right hand side 
$$f(x,t) = \pi \sin(\pi x)(\cos(\pi t) + \pi \sin(\pi t)) + x(x-1)\left(t-
\frac{e^{(t-1)/\varepsilon}-e^{-1/\varepsilon}}{1-e^{-1/\varepsilon}}\right),$$ 
\noindent the desired state $$y_d (x,t) = \sin(\pi x)\sin(\pi t) + x(x-1)\left(1-\frac{e^{(t-1)/\varepsilon}\cdot 1/\varepsilon}{1-e^{-1/\varepsilon}}\right) 
+ 2\left(t-\frac{e^{(t-1)/\varepsilon}-e^{-1/\varepsilon}}{1-e^{-1/\varepsilon}}\right)$$ and the initial 
condition $y_0 (x) = 0$. We choose the regularization parameter to be $\alpha = 1/30$. For 
small values of $\varepsilon$ (we use $\varepsilon = 10^{-4}$), 
the adjoint state $\bar{p}$ develops a layer towards $t = 1$, which can be 
seen in the left plots of Figure \ref{fig1:sol} and Figure \ref{fig2:con}.

\begin{figure}[htbp]
\includegraphics[scale=0.33]{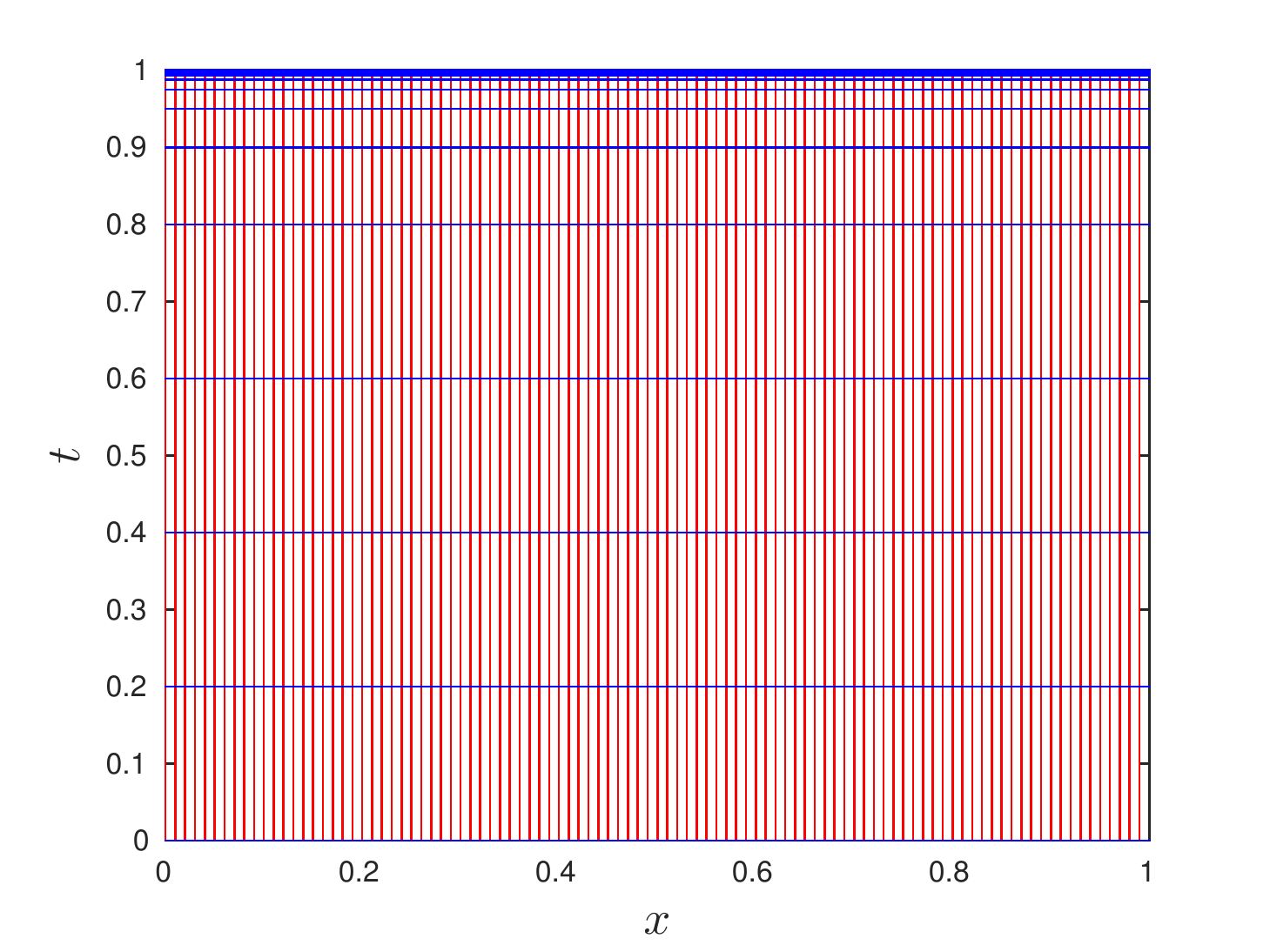}
\includegraphics[scale=0.33]{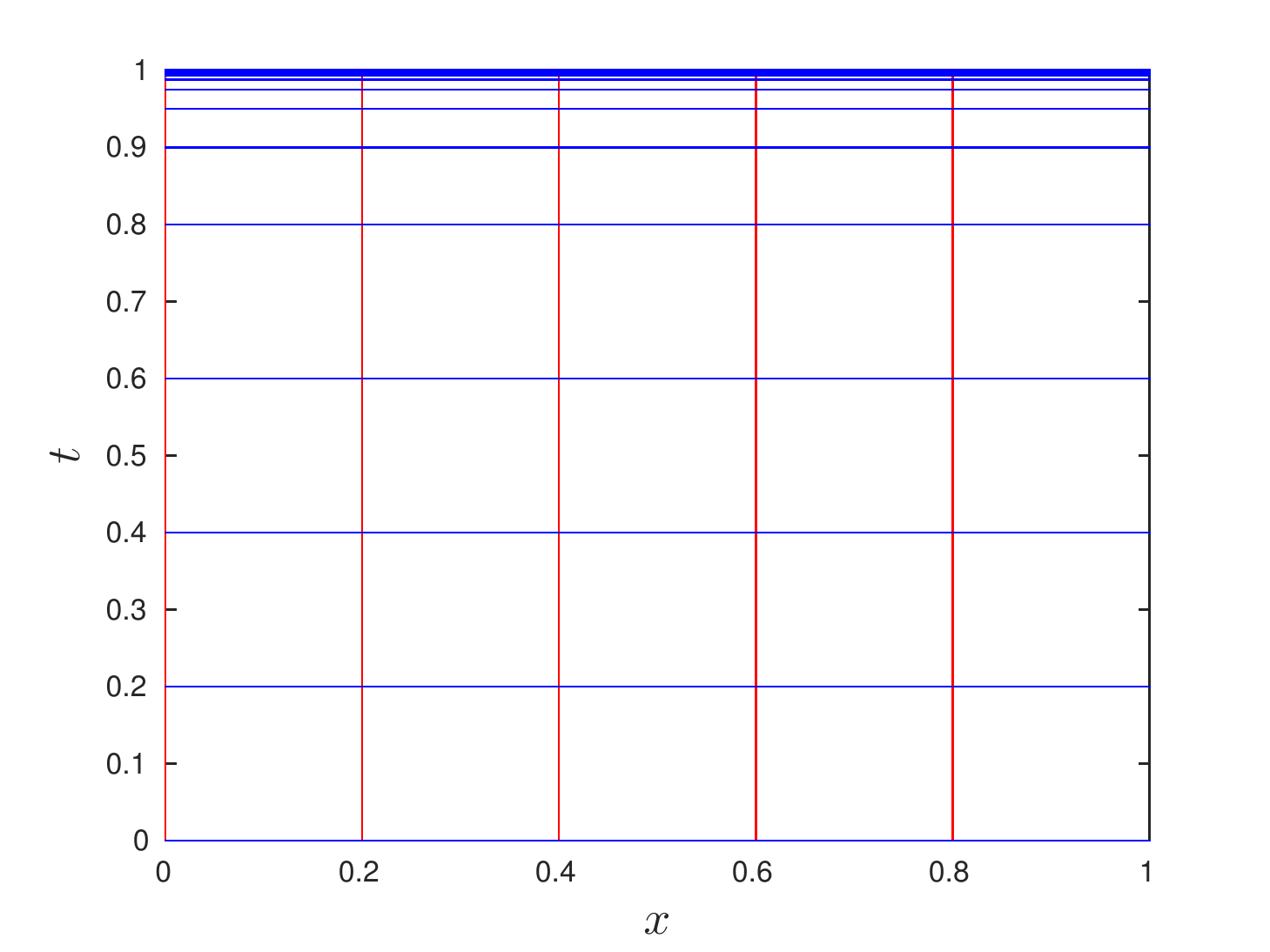}
\includegraphics[scale=0.33]{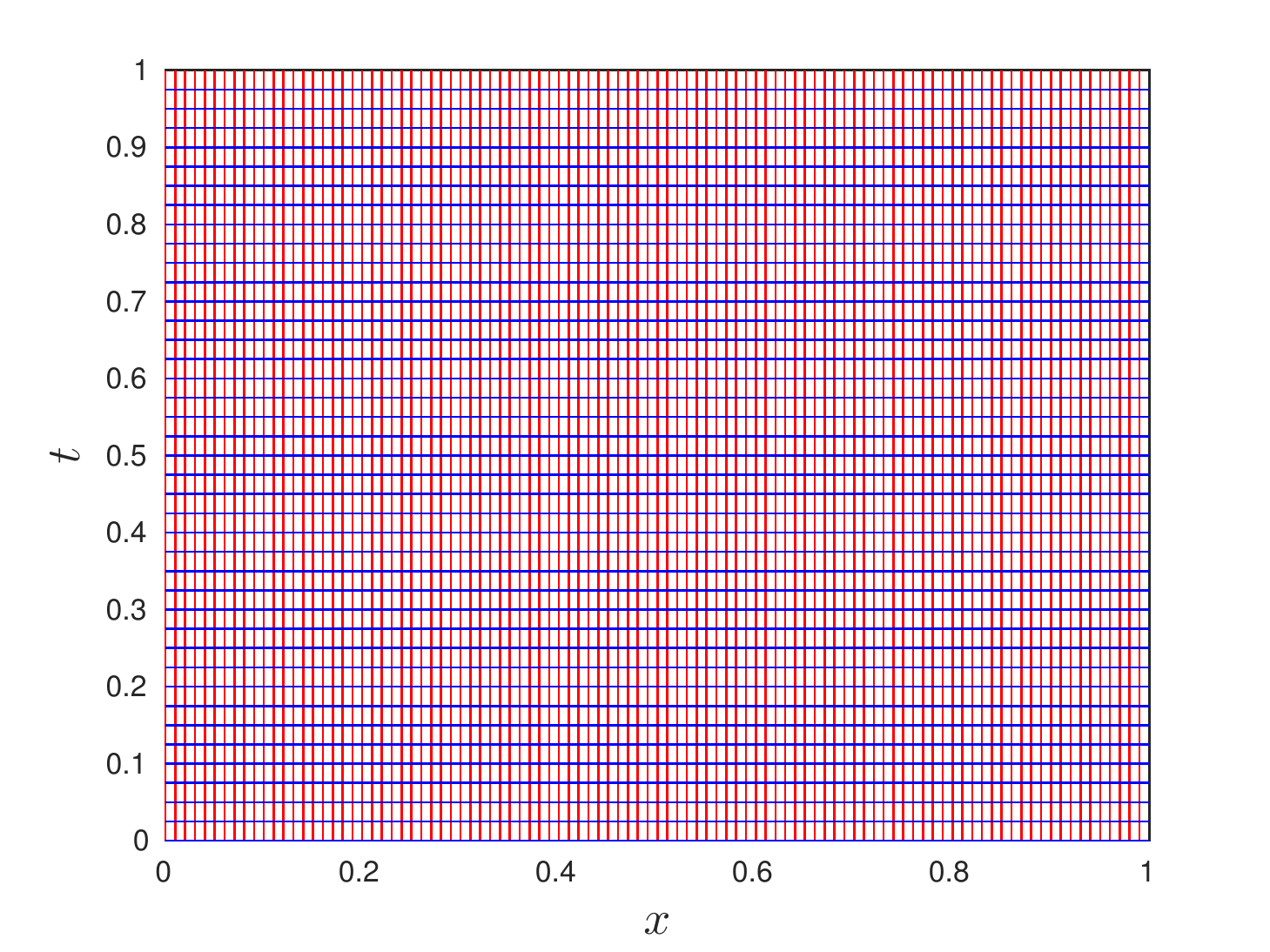}
    \caption{Test 1: Adaptive space-time grids with dof $= 21$
    according to the strategy in \cite{GHZ12} and 
    $\Delta x = 1/100$ (left) and $\Delta x = 1/5$ (middle), 
    respectively, and the equidistant grid with $\Delta t = 1/20$ (right)}
    \label{fig3:mesh}
\end{figure}

 In this test run we focus on the influence of the 
 time grid to approximate of the POD solution. Therefore, we 
 compare the use of two different types of time grids: 
 an equidistant time grid characterized by the time increment $\Delta t = 1/n$ 
 and a non-equidistant (adaptive) time grid characterized by $n+1$ degrees of 
 freedom (dof). We build the POD-ROM from the uncontrolled problem; we create the snapshot ensemble 
 by determining the associated state $y(u_\circ)$ and adjoint state $p(u_\circ)$ corresponding to 
 the control function $u_\circ \equiv 0$ and we also include the initial condition $y_0$ and the time 
 derivatives of the adjoint $p_t(u_\circ)$ into our snapshot set, which is accomplished with time 
 finite differences of the adjoint snapshots. We use $\ell = 1$ POD basis 
 function. Although we would also have the possibility to use 
 suboptimal snapshots corresponding to an approximation $u_{\Delta x}$ 
 of the optimal control, here, we want to emphasize the importance 
 of the time grid. Nevertheless in this example, the quality of the 
 POD solution does not really differ, if we consider 
 suboptimal or uncontrolled snapshots. First, we leave out 
 the post-processing step 4 of Algorithm \ref{Alg:OPTPOD} and 
 discuss the inclusion of this part later.\\
 Figure \ref{fig3:mesh} visualizes the space-time mesh of the numerical solution 
 of (\ref{2ordp}) utilizing the temporal residual type 
 a-posteriori error estimate (\ref{est-thm31}). The first grid in Figure \ref{fig3:mesh} 
 corresponds to the choice of dof=21 and $\Delta x = 1/100$, whereas 
 the grid in the middle refers to using dof = 21 and $\Delta x = 1/5$. 
 Both choices for spatial discretization lead to the exact same time grid, 
 which displays fine time steps towards the end of the time horizon 
 (where the layer in the optimal adjoint state is located), whereas at the beginning and in the middle of the 
 time interval the time steps are larger. This clearly indicates that the 
 resulting time adaptive grid is very insensitive against changes in the spatial resolution. 
For the sake of completeness, the equidistant grid with the same number of degrees of freedom is 
shown in the right plot of Figure \ref{fig3:mesh}.\\
Since the generation of the time adaptive grid as well as the approximation of the 
optimal solution is done in the offline computation part of POD-MOR, this process 
shall be perfomred quickly, which is why we pick $\Delta x = 1/5$ for step 1 in 
Algorithm 1.\\

 \begin{figure}[htbp]
\includegraphics[scale=0.33]{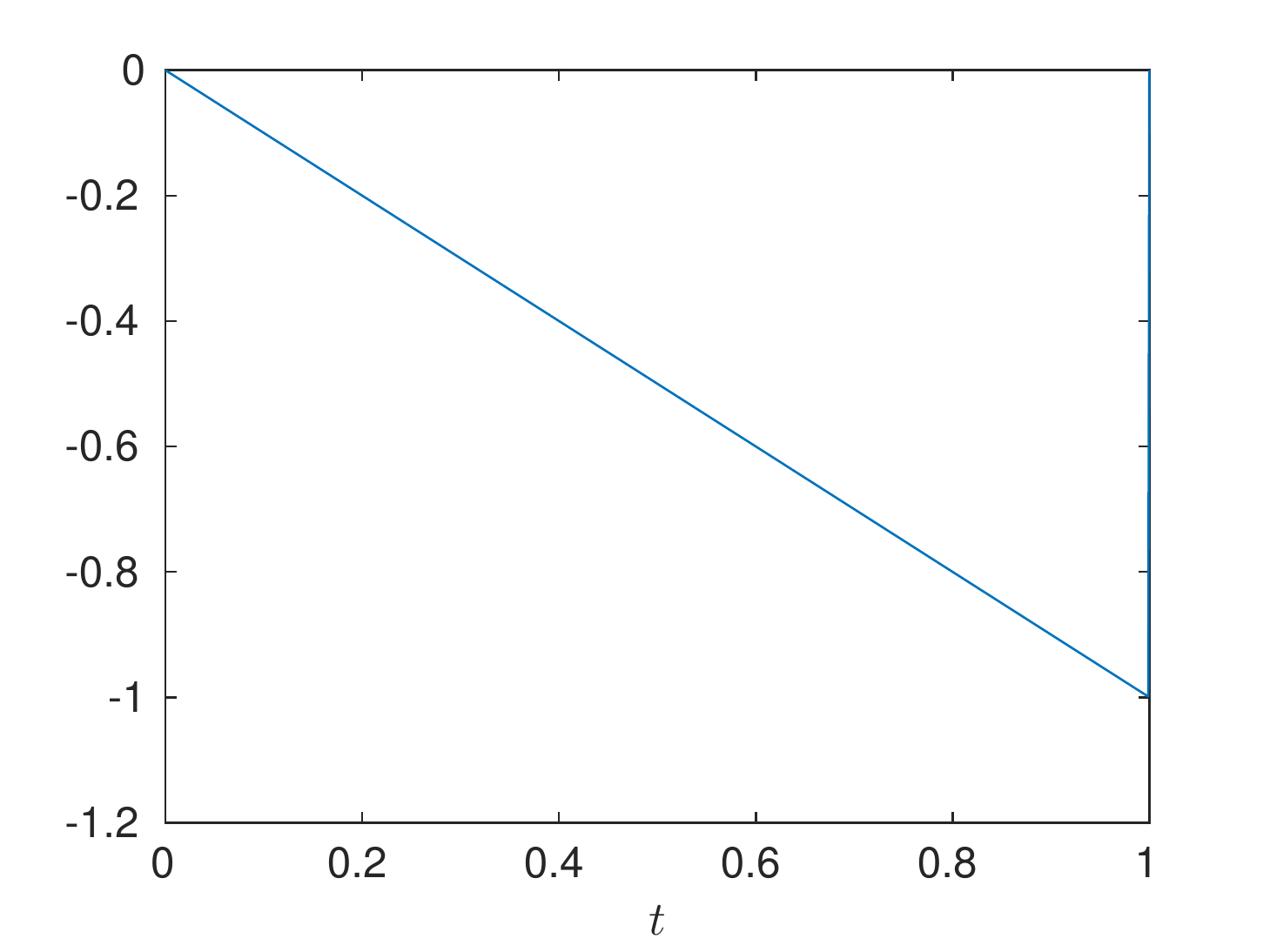}
\includegraphics[scale=0.33]{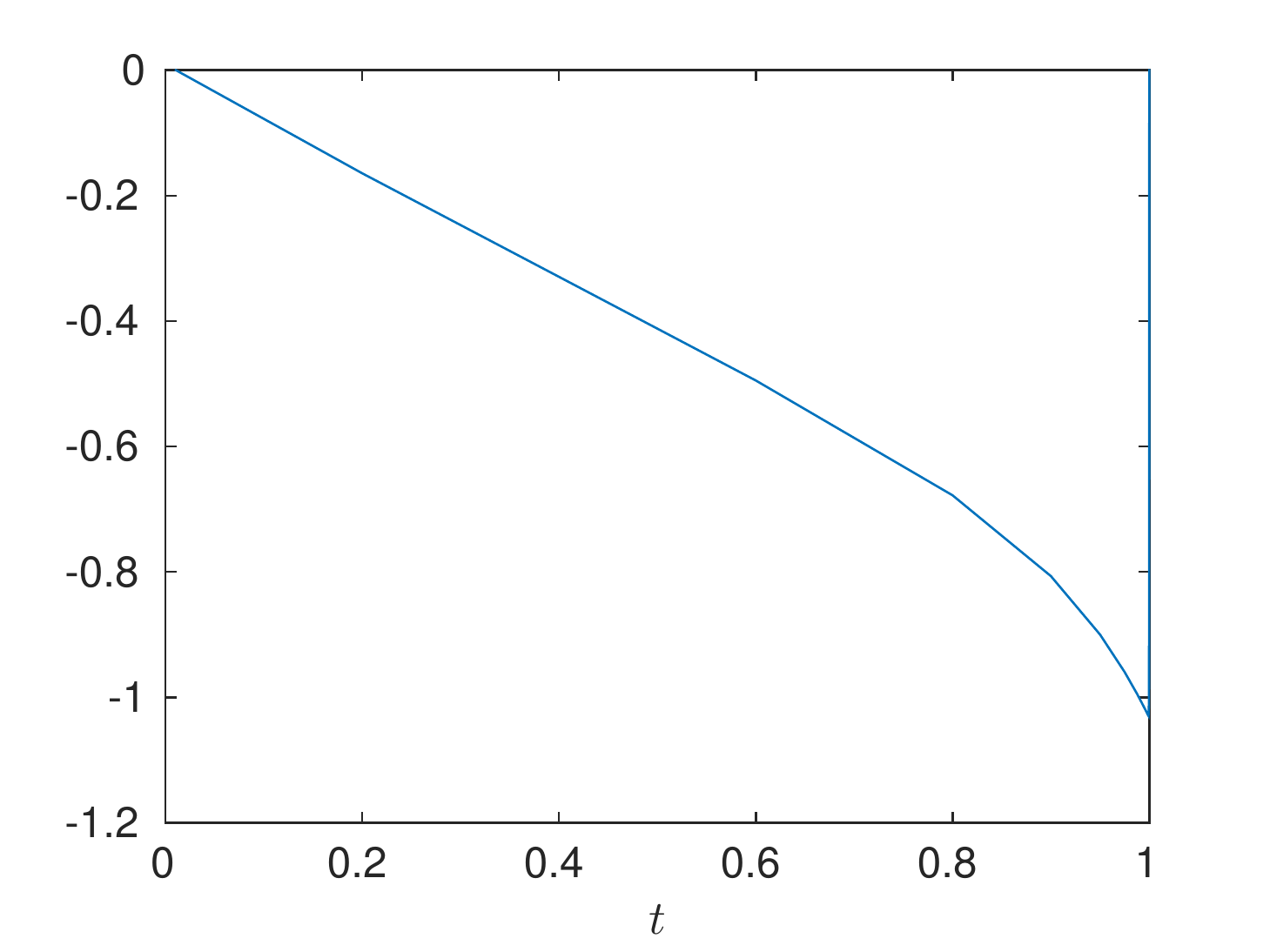}\\
\includegraphics[scale=0.33]{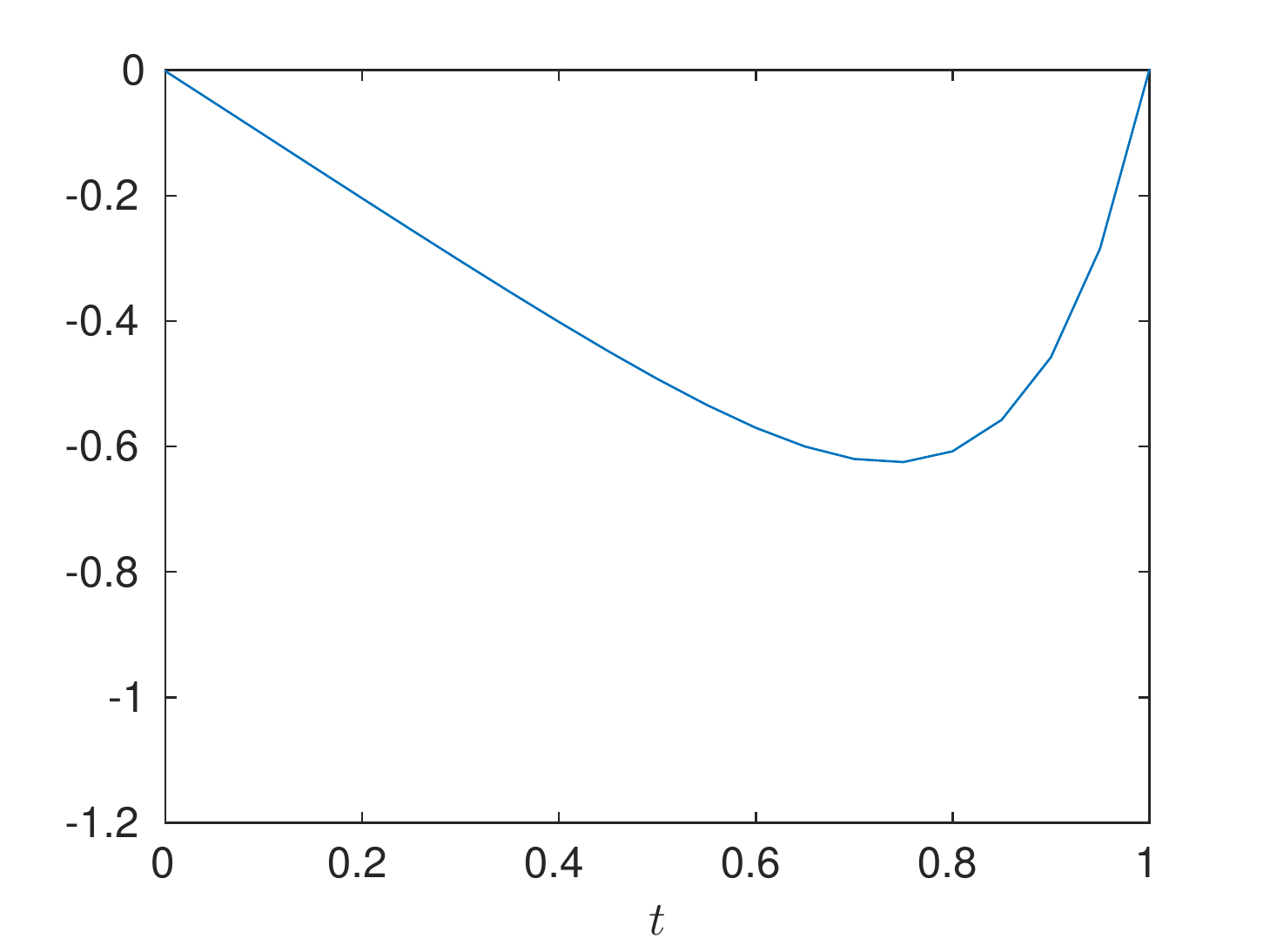}
\includegraphics[scale=0.33]{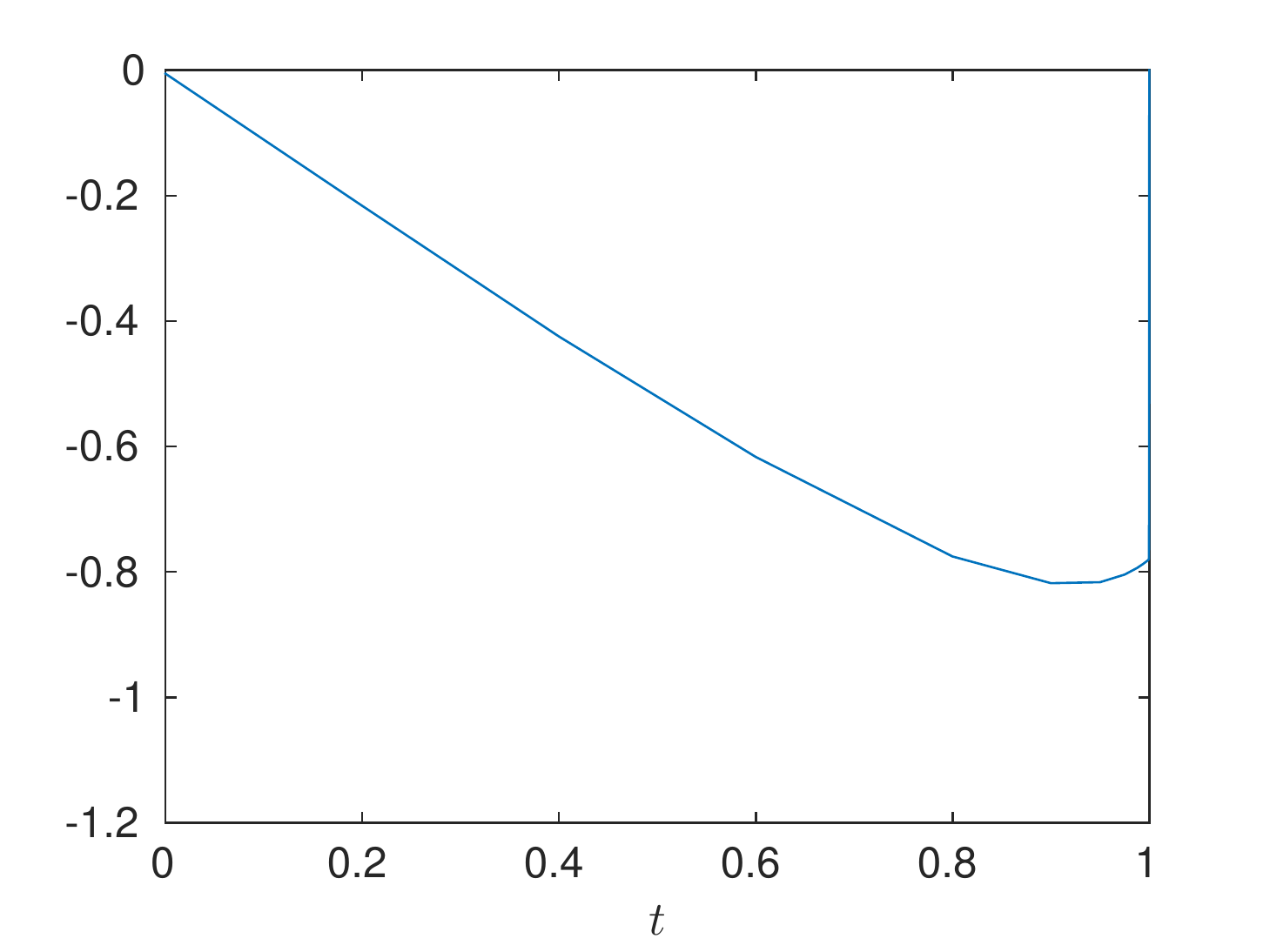}

 \caption{Test 1: Analytical optimal control $\bar{u}$ (top left), 
 approximation $u_{\Delta x}$ of the optimal control gained by step 1
 of Algorithm \ref{Alg:OPTPOD} (top right);
 POD control utilizing a uniform time grid with $\Delta t = 1/20$
 (bottom left), POD control utilizing an adaptive time grid with dof=21
 (bottom right)}
    \label{fig2b:control}
\end{figure}

Figures \ref{fig1:sol} and 
\ref{fig2:con} (middle and right plots) show the surface and contour lines of the 
POD adjoint state utilizing an equidistant time grid and 
utilizing the time adaptive grid, respectively. The analytical 
control intensity $\bar{u}(t)$, the approximation $u_{\Delta x}$ 
of the optimal control computed in step 1 of Algorithm \ref{Alg:OPTPOD} as well 
as the POD controls utilizing a uniform and time adaptive grid, respectively, 
are shown in Figure \ref{fig2b:control}.\\ Table \ref{tab:1} summarizes the approximation quality of the POD solution 
depending on different time discretizations. The fineness of the time discretization 
(characterized by $\Delta t$ and dof, respectively) is chosen in 
such a way that the results of uniform and 
adaptive temporal discretization are comparable. The absolute errors between the 
analytical optimal state $\bar{y}$ and the POD solution $y^\ell$, defined by 
$\varepsilon_{\text{abs}}^y := \|\bar{y}-y^\ell \|_{L^2(\Omega_T)}$, are listed in 
columns 2 and 6; same applies for the errors in the control 
$\varepsilon_{\text{abs}}^u := \|\bar{u}-u^\ell \|_{\mathcal{U}}$ (columns 3 and 7)
and adjoint state $\varepsilon_{\text{abs}}^p 
:= \|\bar{p}-p^\ell \|_{L^2(\Omega_T)}$ (columns 4 and 8). If 
we compare the results, we note that we gain one order of accuracy for the adjoint and control variable with the time adaptive grid. In order to achieve an 
accuracy in the control variable of order $10^{-2}$ utilizing an 
equidistant time grid, we need about $n = 10000$ time steps (not listed 
in Table \ref{tab:1}). This emphasizes that using an appropriate 
(non-equidistant) time grid for the adjoint variable is of particular 
importance in order to efficiently achieve POD controls of good quality.

\begin{table}[htbp]
\centering
 \begin{tabular}{ c | c | c | c || c | c | c | c}
 \toprule
 $\Delta t$ & $\varepsilon_{\text{abs}}^y$ 
 & $\varepsilon_{\text{abs}}^u$ & $\varepsilon_{\text{abs}}^p$ & dof & 
 $\varepsilon_{\text{abs}}^y$ & $\varepsilon_{\text{abs}}^u$ &  $\varepsilon_{\text{abs}}^p$ \\
 \hline
 1/20 & $1.5120 \cdot 10^{-02}$  & $1.9837 \cdot 10^{-01}$  & $3.6247 \cdot 10^{-02}$  & 21 &  $5.1874 \cdot 10^{-02}$ &  $5.3428 \cdot 10^{-02}$ &  $9.6343 \cdot 10^{-03}$ \\
 1/42 & $1.1186 \cdot 10^{-02}$  & $2.1071 \cdot 10^{-01}$  & $3.8490 \cdot 10^{-02}$  & 43 &  $5.1634 \cdot 10^{-02}$ &  $2.4868 \cdot 10^{-02}$ &  $4.3611 \cdot 10^{-03}$ \\
 1/61 & $1.0774 \cdot 10^{-02}$  & $2.1447 \cdot 10^{-01}$  & $3.9173 \cdot 10^{-02}$  & 62 &  $5.1599 \cdot 10^{-02}$ &  $2.3275 \cdot 10^{-02}$ &  $4.0691 \cdot 10^{-03}$ \\
 1/114 & $1.1157 \cdot 10^{-02}$  & $2.1846 \cdot 10^{-01}$  & $3.9893 \cdot 10^{-02}$  & 115 &  $5.1568 \cdot 10^{-02}$ &  $2.3027 \cdot 10^{-02}$ & $4.0340 \cdot 10^{-03}$\\
 \bottomrule 
 \end{tabular}
 \vspace{0.4cm} \caption{Test 1: Absolute errors between the analytical optimal 
 solution and the POD solution depending on the time discretization (equidistant: 
 columns 1-4, adaptive: columns 5-8)}
 \label{tab:1}
  \end{table}

Table \ref{tab:2} contains the evaluations of each term in \eqref{err:summarized}. 
The value $\eta_p^i$ ($\eta_p^b$) refers to the first (second) part in (\ref{est-thm31}). 
For this test example, we note that the term $\eta_p^i$ influences the estimation. However, we observe that the better the semi-discrete adjoint state $p_{\Delta x}$
from step 1 of Algorithm \ref{Alg:OPTPOD} is, the better 
will be the POD adjoint solution. Since all summands of (\ref{err:summarized}) 
can be estimated, Table \ref{tab:2} allows us to control the approximation of 
the POD adjoint state. The estimation \eqref{post_p} concerning the state variable 
will be investigated later on.

\begin{table}[htbp]
\centering
 \begin{tabular}{ c | c | c | c | c | c }
 \toprule
 dof & $\varepsilon_{\text{abs}}^p$ 
  & $\eta_p^i$ & $\eta_p^b$   & $\| \zeta_k \|_U + \| \zeta_k^\ell \|_U $   & $\sum_{i=\ell+1}^d \lambda_i$ \\
 \hline
 21 & $9.6343 \cdot 10^{-03}$  & $4.9518 \cdot 10^{+00}$  & $4.8031 \cdot 10^{-04}$  & $1.6033 \cdot 10^{-02}$    & $3.3938 \cdot 10^{-04}$   \\
 43 & $4.3611 \cdot 10^{-03}$  & $1.1976 \cdot 10^{+00}$  & $5.0087 \cdot 10^{-05}$   & $1.9200 \cdot 10^{-02}$   &  $2.9454 \cdot 10^{-04}$    \\
 62 & $4.0691 \cdot 10^{-03}$  & $7.2852 \cdot 10^{-01}$ & $2.9835 \cdot 10^{-05}$   & $1.9707 \cdot 10^{-02}$   &  $2.9212 \cdot 10^{-04}$ \\
 115 & $4.0340 \cdot 10^{-03}$  & $3.4966 \cdot 10^{-01}$ & $1.4845 \cdot 10^{-05}$  & $2.0191 \cdot 10^{-02}$  &   $2.9090 \cdot 10^{-04}$ \\
 \bottomrule 
 \end{tabular}
 \vspace{0.4cm} \caption{Test 1: Evaluation of each summand of the error 
 estimation \ref{err:summarized}}
 \label{tab:2}
  \end{table}

\begin{figure}[htbp]
\includegraphics[scale=0.33]{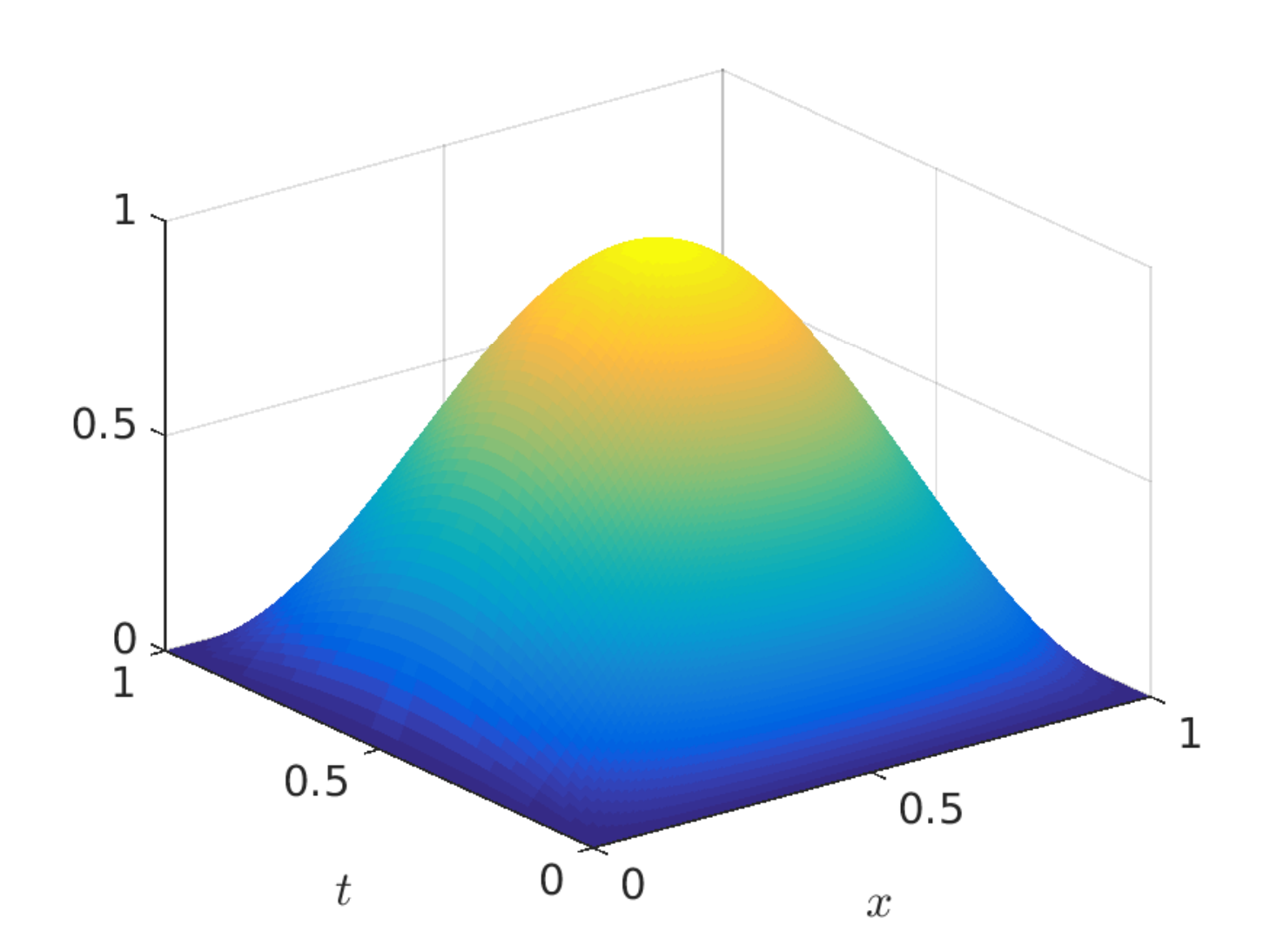}
\includegraphics[scale=0.33]{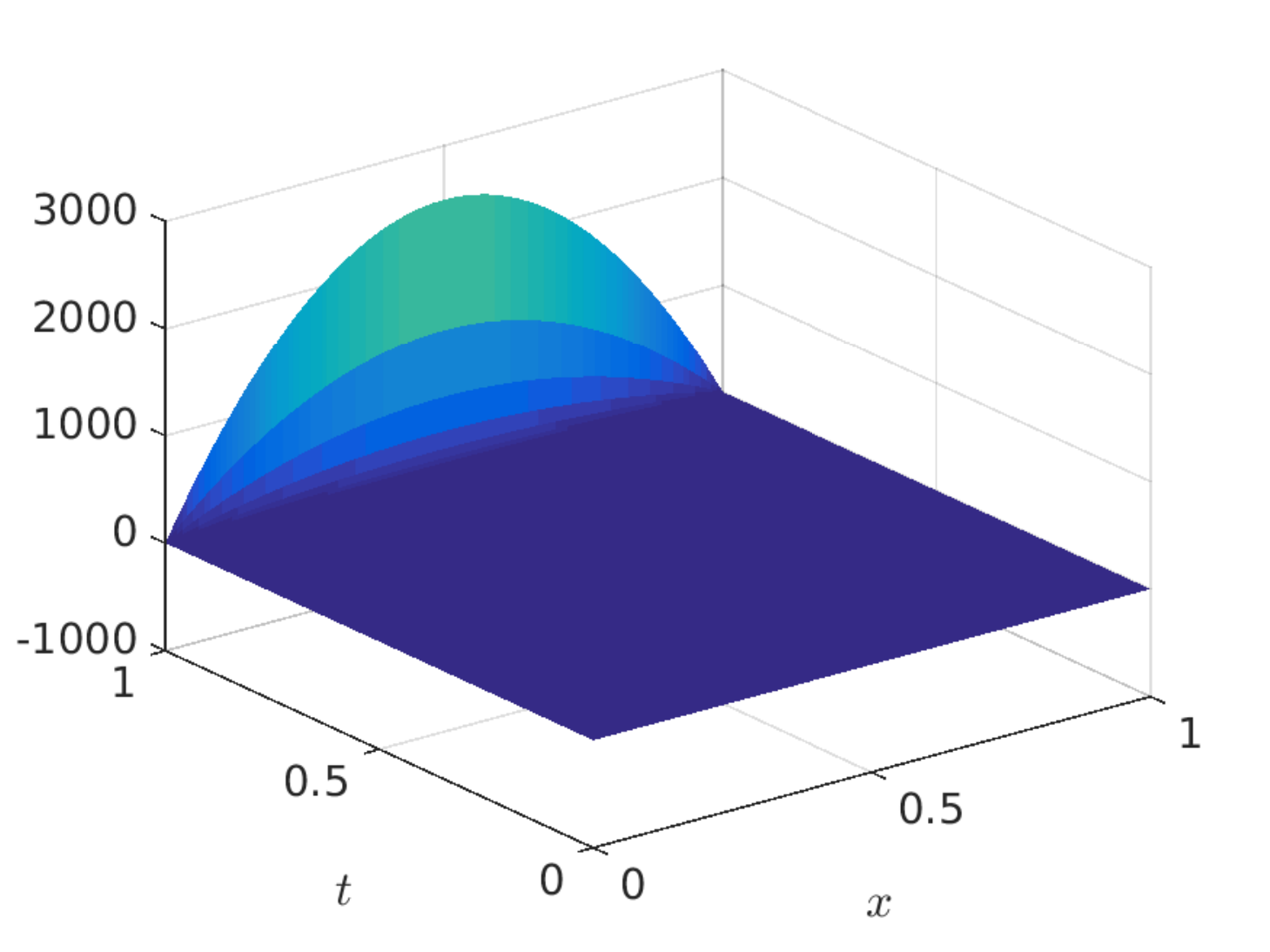}\\
\includegraphics[scale=0.33]{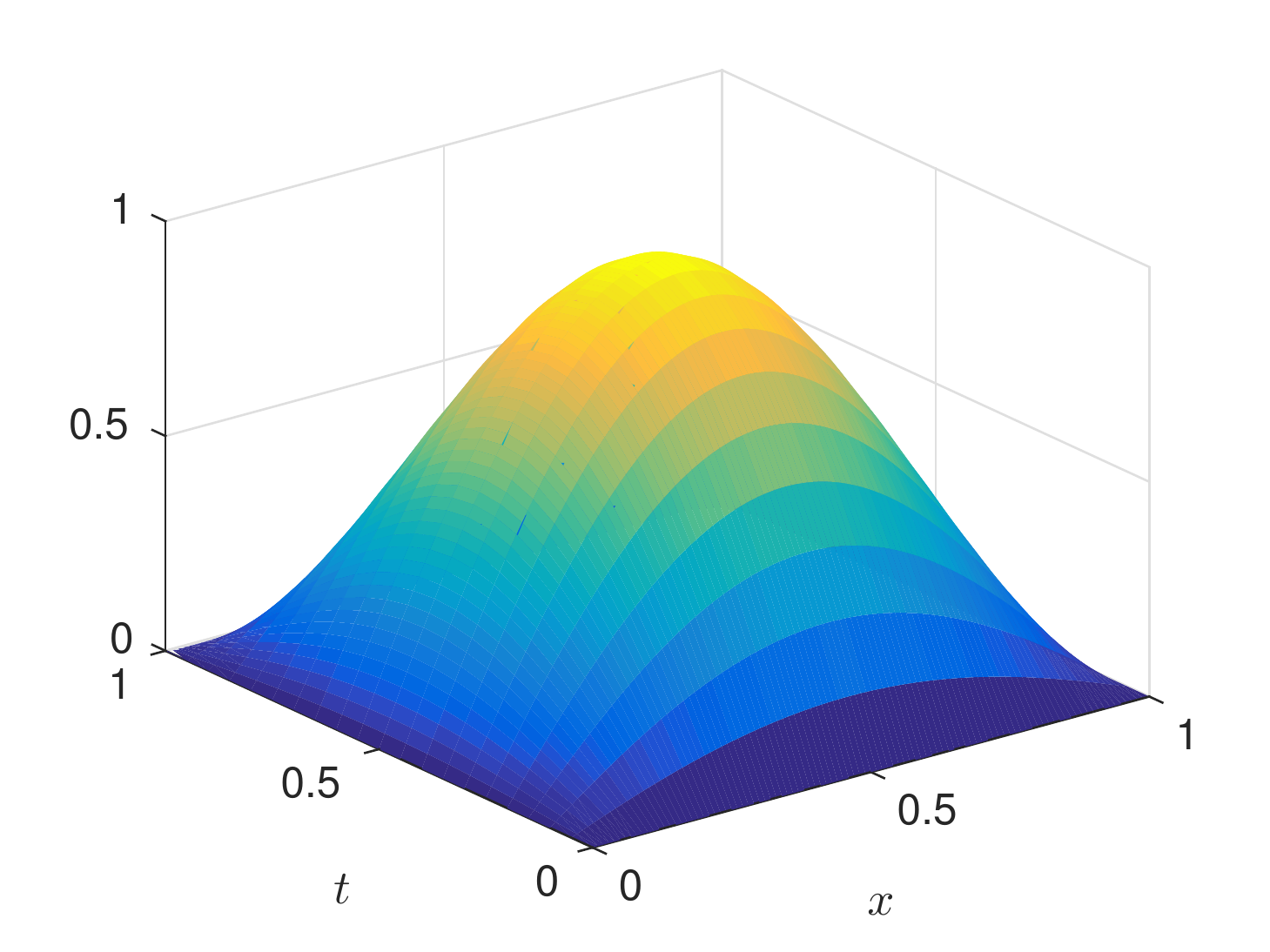}
\includegraphics[scale=0.33]{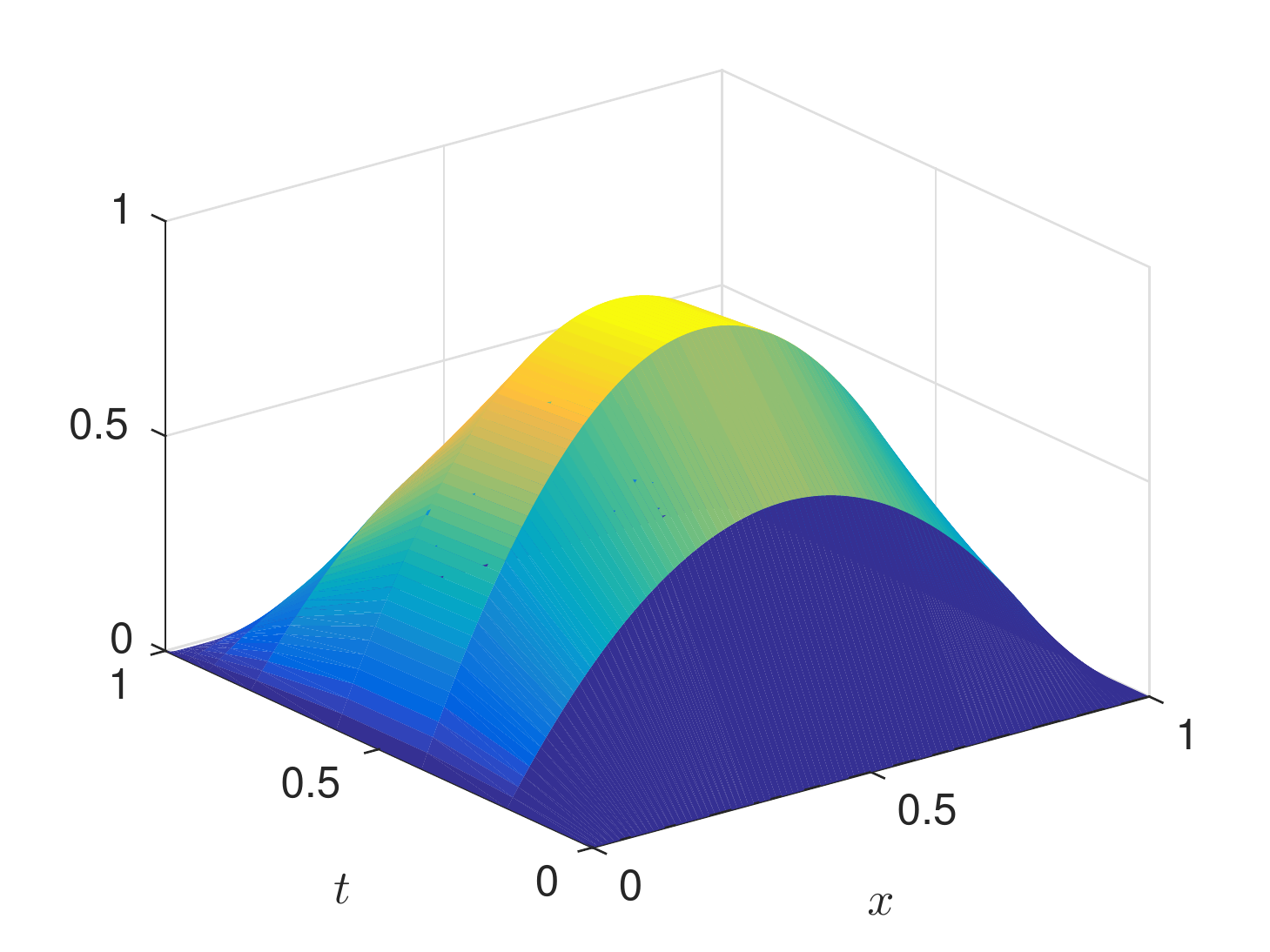}

 \caption{Test 1: Analytical optimal state $\bar{y}$ (top left), desired state $y_d$ 
 (top right); POD state $y^\ell$ utilizing a uniform time grid with $\Delta t = 1/20$ 
 (bottom left), POD state $y^\ell$ utilizing an adaptive time grid with dof = 21 (bottom right) }
    \label{fig4:states}
\end{figure}

  Moreover, a comparison of the value of the cost functional is given in 
  Table \ref{tab:3}. The aim of the optimization problem (\ref{ocp}) is to 
  minimize the quantity of interest $J(y,u)$. The analytical value of the 
  cost functional at the optimal solution is $J(\bar{y},\bar{u}) \approx 8.3988 
  \cdot 10^{+01}$. Table \ref{tab:3} clearly points out that the use of 
  a time adaptive grid is fundamental for solving the optimal control 
  problem (\ref{ocp}). The huge differences in the values of the cost functional 
  is due to the great increase of the desired state $y_d$ at the end of the time 
  interval (see Figure \ref{fig4:states}). Small time steps at the end of the time interval, as it is the 
  case in the time adaptive grid, lead to much more accurate results.

\begin{table}[htbp]
\centering
 \begin{tabular}{ c | c || c | c  }
 \toprule
 $\Delta t$ & $J(y^\ell,u)$ & dof & $J(y^\ell,u)$  \\
 \hline
 1/20 & $4.1652 \cdot 10^{+04}$ &  21 &  $8.7960 \cdot 10^{+01}$    \\
 1/42 &  $1.9834 \cdot 10^{+04}$ & 43 &  $8.4252 \cdot 10^{+01}$ \\
 1/61 &  $1.3656 \cdot 10^{+04}$ & 62 & $8.4102 \cdot 10^{+01}$ \\
 1/114 &  $7.3078 \cdot 10^{+03}$ & 115 & $8.4034 \cdot 10^{+01}$ \\
 1/40000 & $8.5692 \cdot 10^{+01}$ & - & - \\
 \bottomrule 
 \end{tabular}
 \vspace{0.4cm} \caption{Test 1: Value of the cost functional at the POD  
 solution utilizing uniform and adaptive time discretization, respectively, analytical 
 value: $J \approx 8.3988 \cdot 10^{+01} $}
 \label{tab:3}
  \end{table}

  Now, let us discuss the inclusion of step 4 in Algorithm \ref{Alg:OPTPOD}.
  Since we went for an adaptive time grid regarding the adjoint variable, we cannot 
  in general expect that the resulting time grid is a good time grid for the state 
  variable. Table \ref{tab:1} confirms that utilizing a uniform time grid leads to better approximation results in the state variable than using the time adaptive grid.
  In order to improve also the approximation quality in the state variable, we incorporate 
  the error estimation (\ref{post_p}) from \cite{KV02} in a post-processing 
  step after producing the time grid with the strategy of \cite{GHZ12} and 
  before starting the POD solution process. Define\\

  \begin{center}
   $\eta_{\text{POD}_j} := 
    \Delta t_j^2 \left( \int\limits_{I_j} (\| y_{tt}^k \|_H^2 + 
   \|y_t^k\|_V^2) \right) $
  \end{center} 

  where $y_t^k \approx y_t(t_k)$ and $y_{tt}^k \approx y_{tt}(t_k)$ 
  are computed via finite difference approximation. We perform bisection on those 
  time intervals $I_j$, where the quantity $\eta_{\text{POD}_j}$ has its 
  maximum value and repeat this $N_{\text{refine}}$ times. This results in 
  the time grid $\mathcal{T_{\text{new}}}$. The improvement in the 
  approximation quality in the state variable can be seen in Table \ref{tab:4}. 
  The more additional time instances we include according to \eqref{post_p}, 
  the better the approximation results get with respect to the state. Moreover, 
  also the approximation quality in the control and 
  adjoint state is improved.

\begin{table}[htbp]
\centering
 \begin{tabular}{ c | c | c | c  }
 \toprule
 $N_{\text{refine}}$ & $\varepsilon_{\text{abs}}^y$ & $\varepsilon_{\text{abs}}^u$ & $\varepsilon_{\text{abs}}^p$ \\
 \hline
 0 & $5.1874 \cdot 10^{-02}$  & $5.3428 \cdot 10^{-02}$ & $9.6343 \cdot 10^{-03}$\\
 5 & $ 4.0058 \cdot 10^{-02}$ &  $2.1145 \cdot 10^{-02}$ &  $3.6378 \cdot 10^{-03}$   \\
 10 &  $3.0909 \cdot 10^{-02}$ & $1.8396 \cdot 10^{-02}$ &  $3.0895 \cdot 10^{-03}$ \\
 20 &  $2.4759 \cdot 10^{-02}$ & $1.7104 \cdot 10^{-02}$ & $2.8210 \cdot 10^{-03}$\\
 30 &  $2.3028 \cdot 10^{-02}$ & $1.6971 \cdot 10^{-02}$ & $2.7906 \cdot 10^{-03}$ \\
 \bottomrule 
 \end{tabular}
 \vspace{0.4cm} \caption{Test 1: Improvement of approximation quality concerning 
 the state variable. The initial time grid $\mathcal{T}$ is computed with dof=43}
 \label{tab:4}
  \end{table}

  We note that the sum of the neglected eigenvalues $\sum_{i=2}^d \lambda_i$ is approximately 
  zero and the second largest eigenvalue of 
  the correlation matrix is of order $10^{-10}$, which makes 
  the use of additional POD basis functions redundant. Likewise, in this 
  particular example the choice of richer snapshots (even the optimal snapshots) 
  does not bring significant 
  improvements in the approximation quality of the POD solutions. So, this 
  example shows that solely the use of an appropriate adaptive time mesh 
  efficiently improves the accuracy of the POD solution.\\


  \subsection{Test 2: Solution with steep gradient in the middle of the time interval} 

  Let $\Omega = (0,1)$ be the spatial domain and $[0,T]=[0,1]$ be the time interval. 
  We choose $\varepsilon = 10^{-4}$ and $\alpha = 1$. To begin with, we 
  consider an unconstrained optimal control problem and investigate the inclusion 
  of control constraints separately in Test 3. We build the example in 
  such a way that the analytical solution $(\bar{y},\bar{u})$ of \eqref{ocp}
  is given by:
  
  $$ \bar{y} (x,t) = x^3 (x-1)t, \quad \bar{p}(x,t) = \sin 
  (\pi x) \text{atan}\left(\frac{t-0.5}{\varepsilon}\right)(t-1), $$\\
  $$ \bar{u}_1(t) = \bar{u}_2(t) = -\text{atan}
  \left(\frac{t-0.5}{\varepsilon}\right)(t-1)\left(\frac{32}{\pi^3} - 
  \frac{8}{\pi^2}\right), $$\\ 
  $$\bar{\chi}_1(x) = \max(0,1-16(x-0.25)^2), \quad \bar{\chi}_2(x) = 
  \max(0,1-16(x-0.75)^2).$$\\
  
\begin{figure}[htbp]
\includegraphics[scale=0.33]{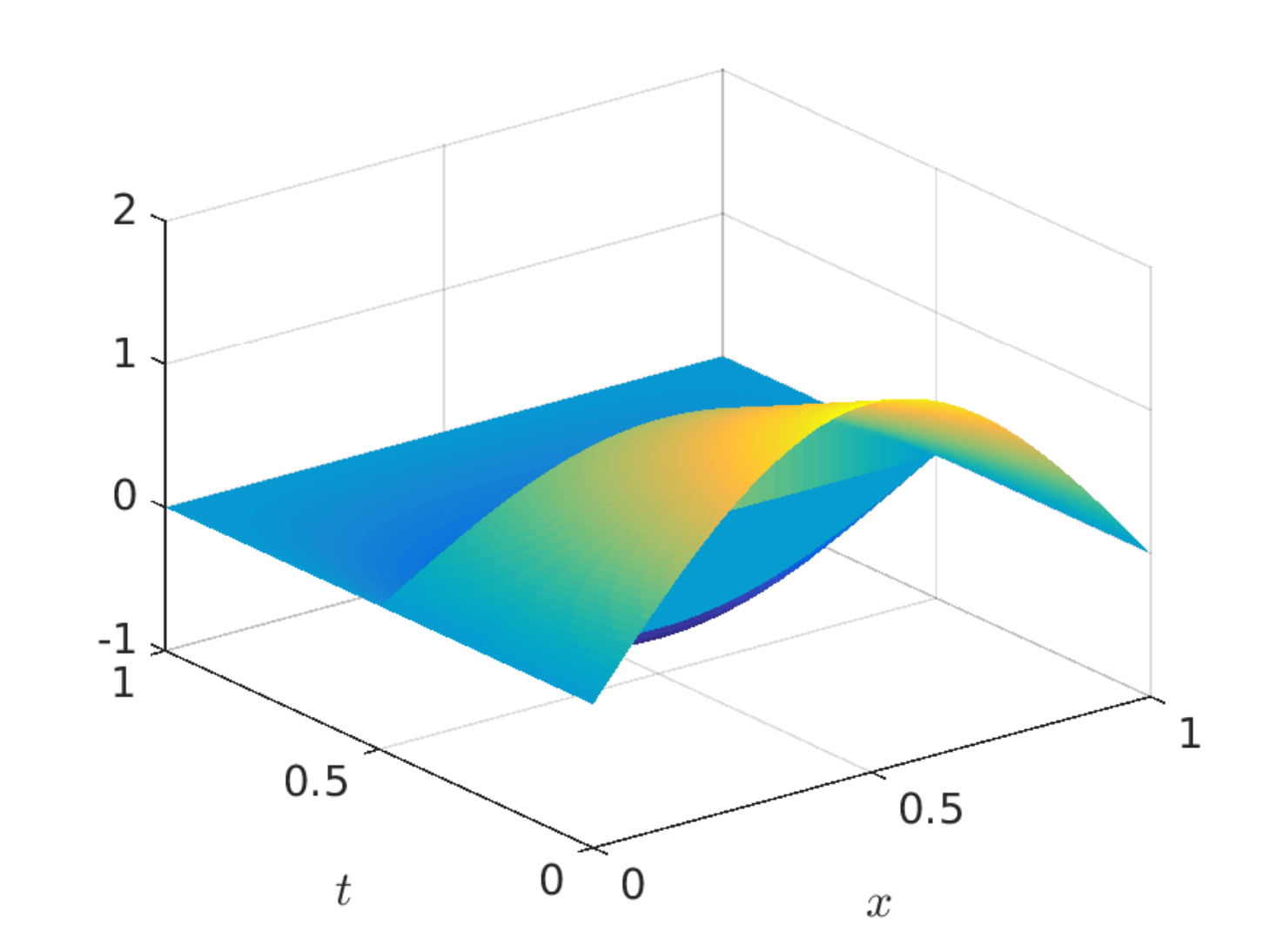} \includegraphics[scale=0.33]{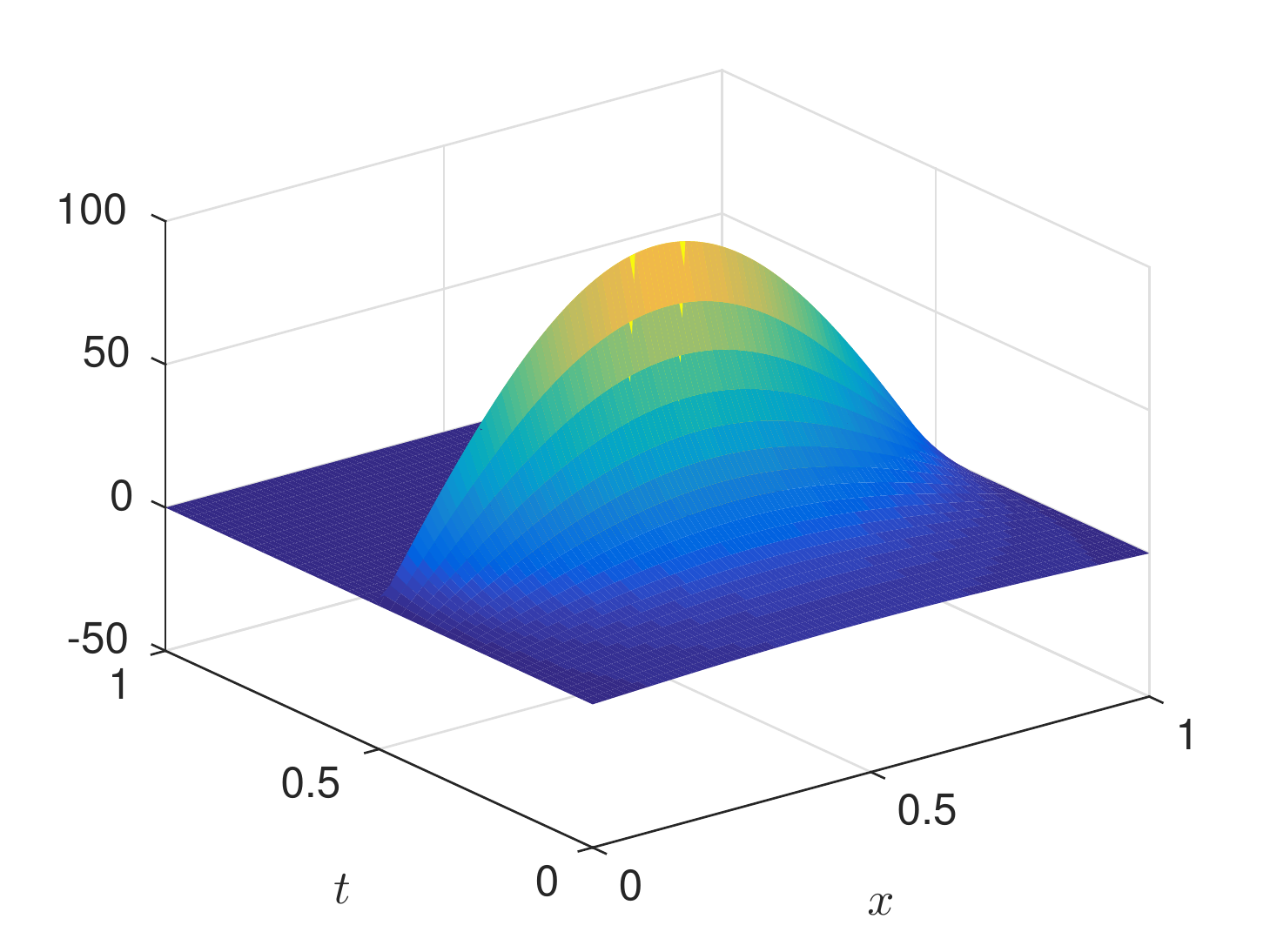} \includegraphics[scale=0.33]{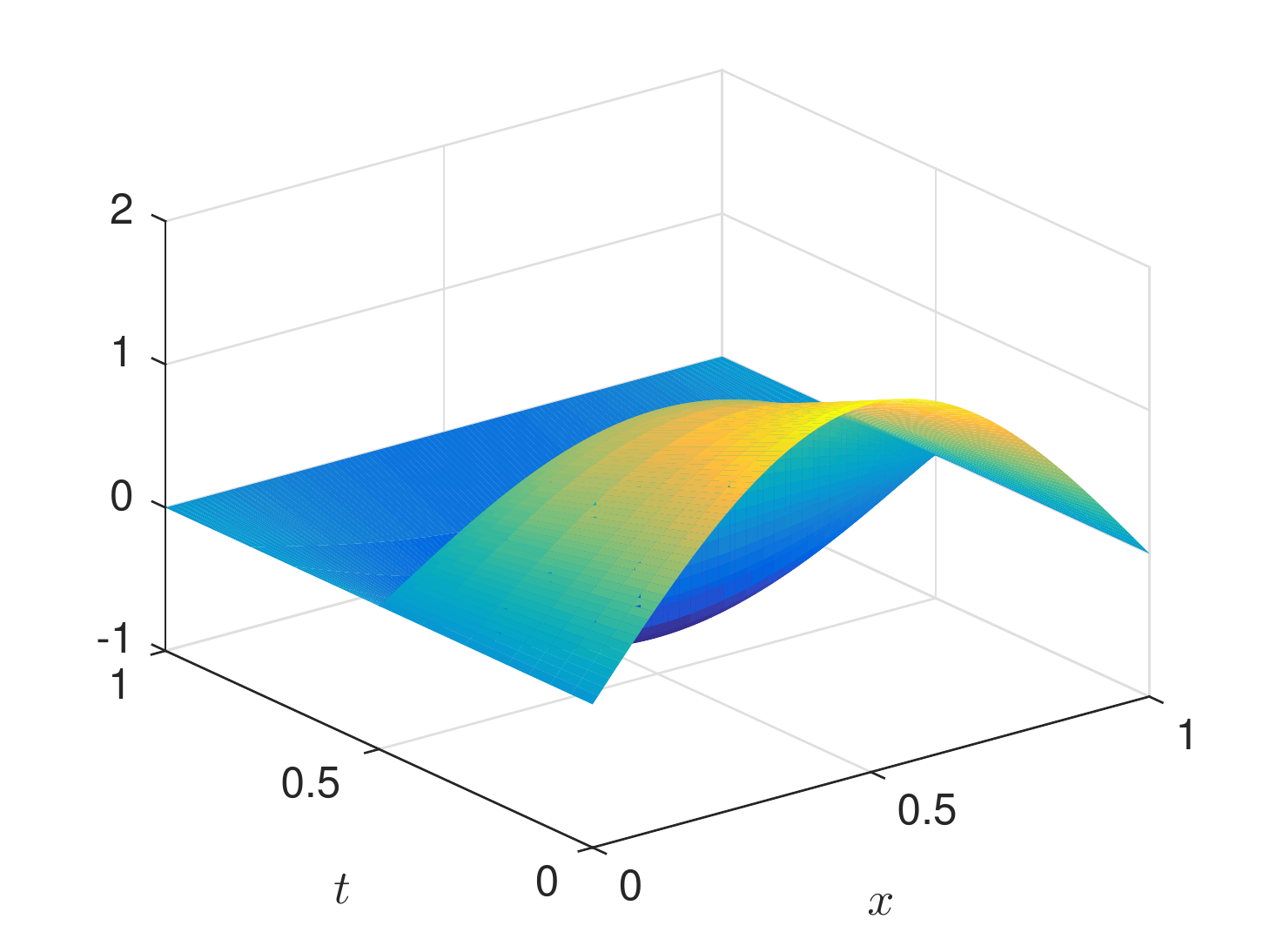}
 \caption{Test 2: Analytical optimal adjoint state $\bar{p}$ (left), POD adjoint 
 solution $p^\ell$ with $\ell = 4$ utilizing an equidistant time grid with 
 $ \Delta t = 1/40$ (middle), POD adjoint solution $p^\ell$ with $\ell = 4$ utilizing 
 an adaptive time grid with dof=41 (right)}
 \label{fig:ex2_p_surf}
 \end{figure}

\begin{figure}[htbp]
 \includegraphics[scale=0.33]{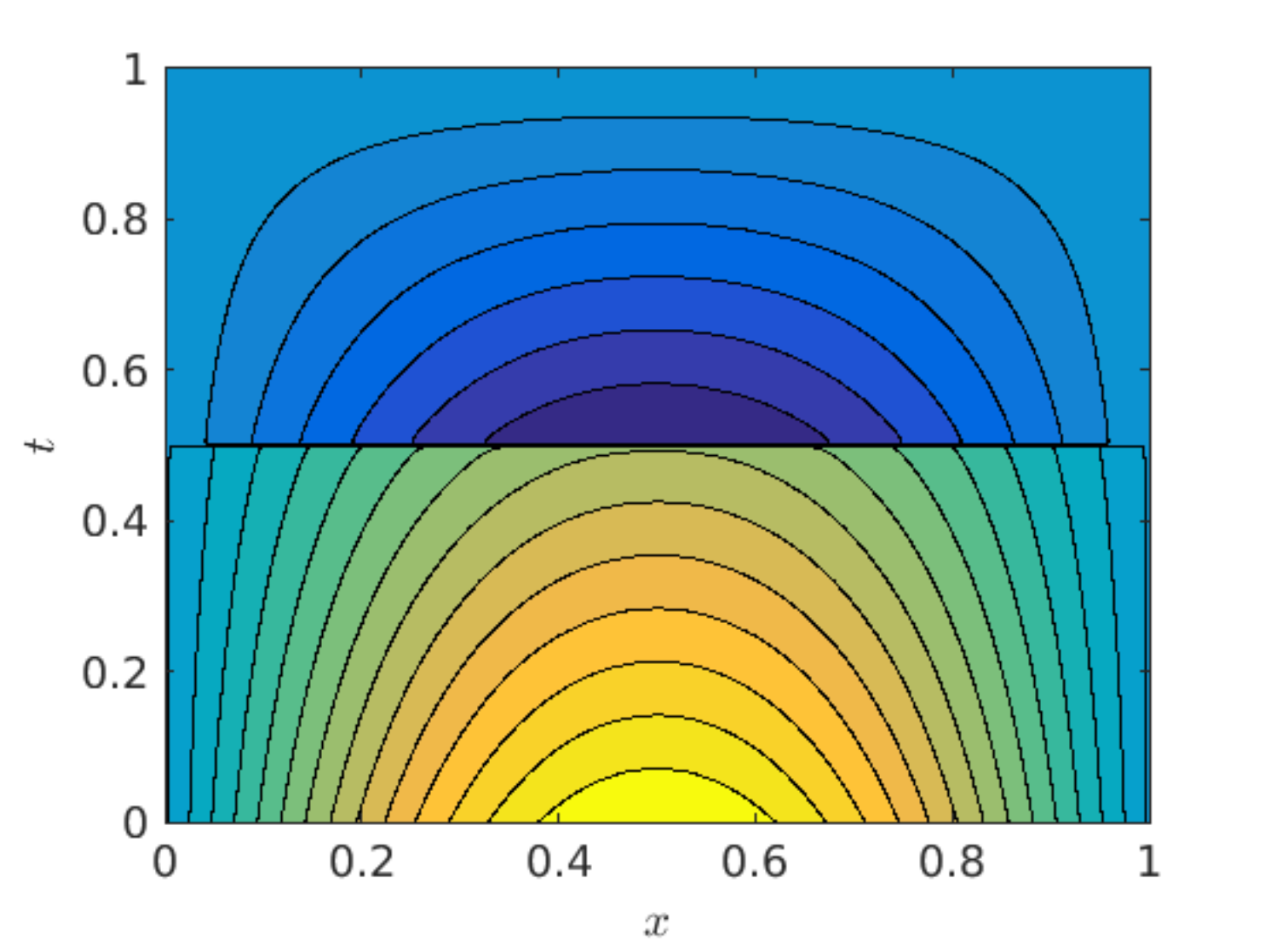}  \includegraphics[scale=0.33]{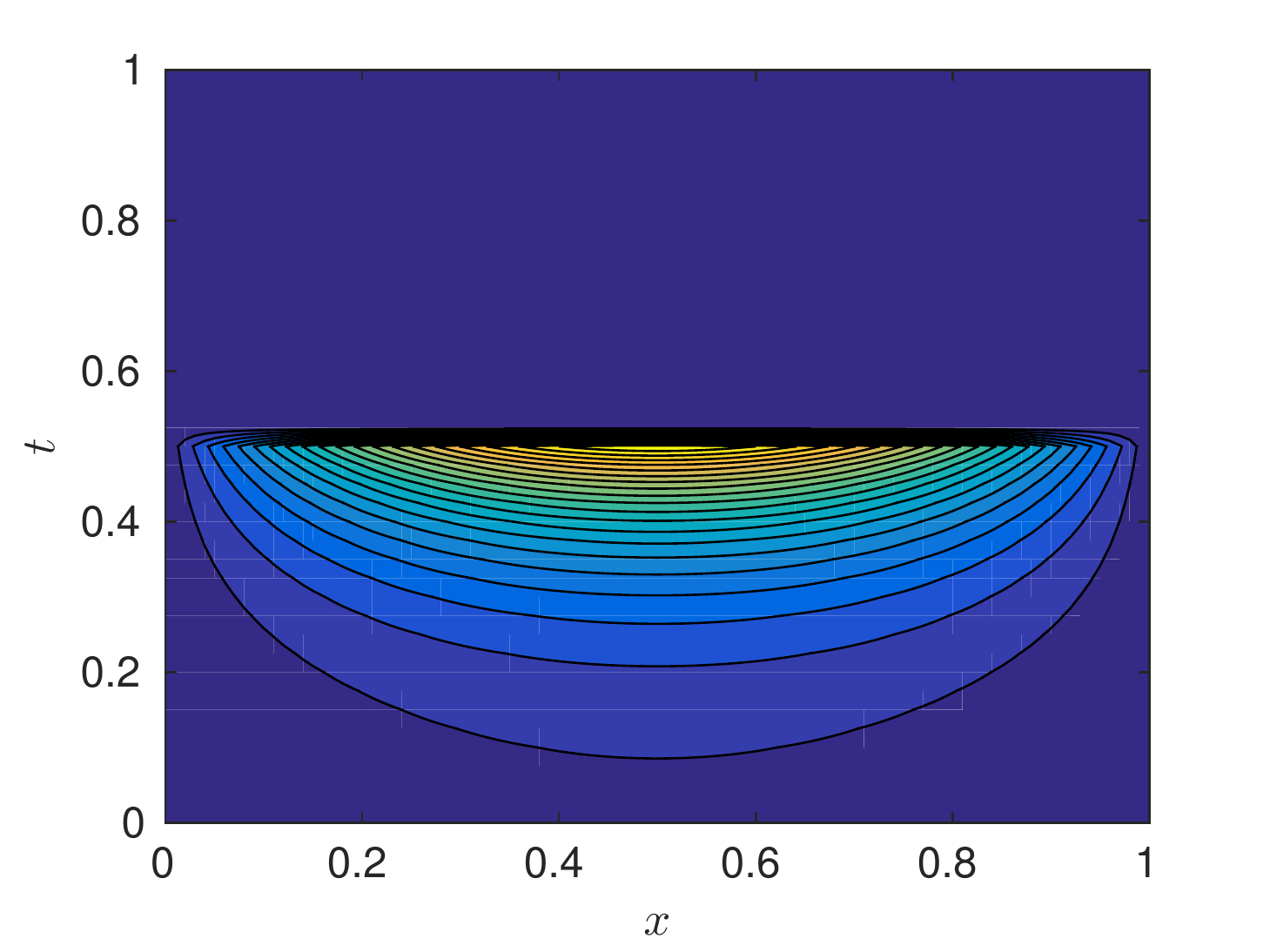}  \includegraphics[scale=0.33]{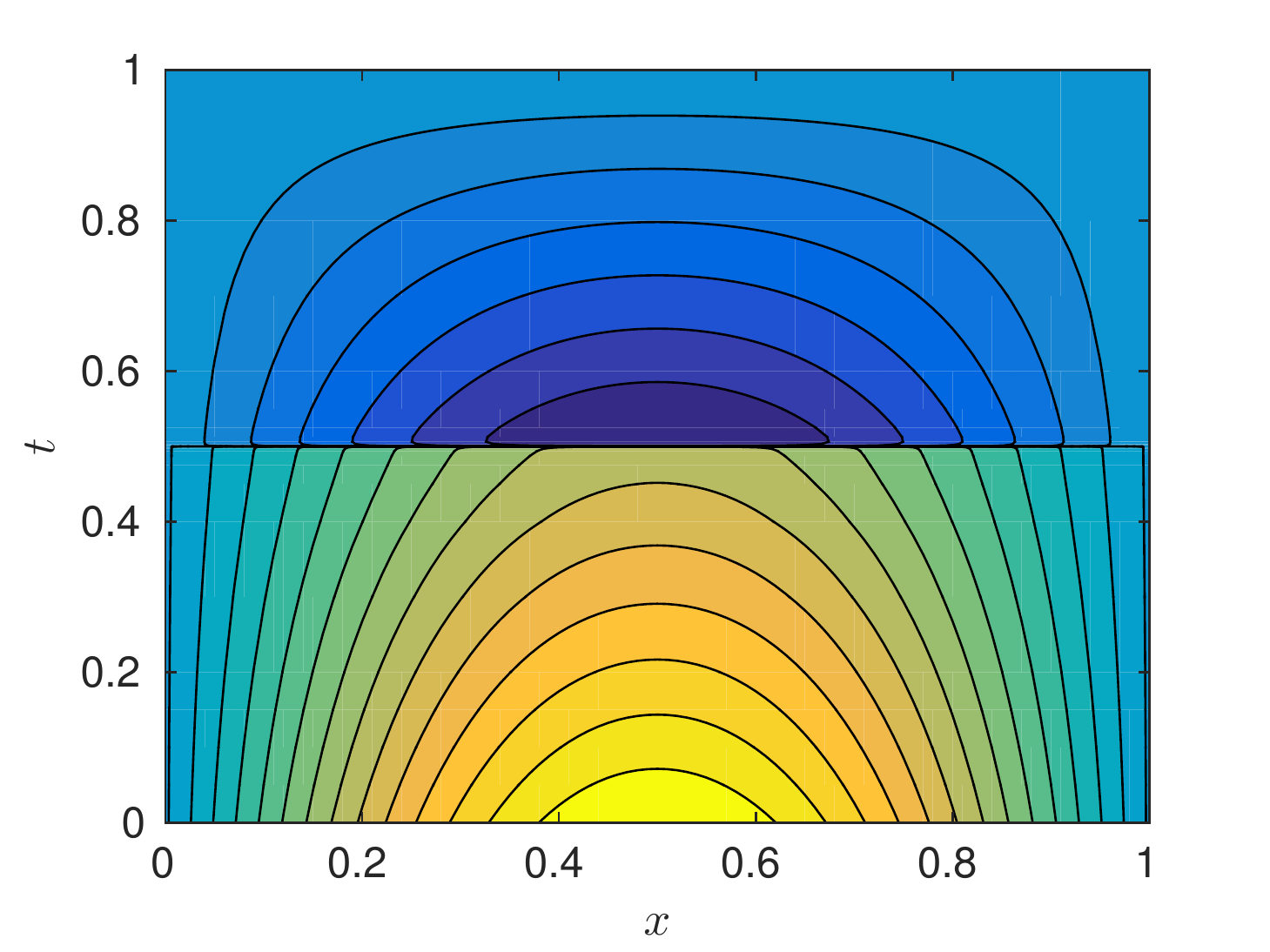}
 \caption{Test 2: Contour lines of the analytical optimal adjoint state 
 $\bar{p}$ (left), POD adjoint solution $p^\ell$ with $\ell = 4$ utilizing an equidistant 
 time grid with $ \Delta t = 1/40$ (middle), POD adjoint solution $p^\ell$ with $\ell = 4$ 
 utilizing an adaptive 
 time grid with dof=41 (right)}
 \label{fig:ex2_p_contour}
 \end{figure}

\begin{figure}[htbp]
\includegraphics[scale=0.33]{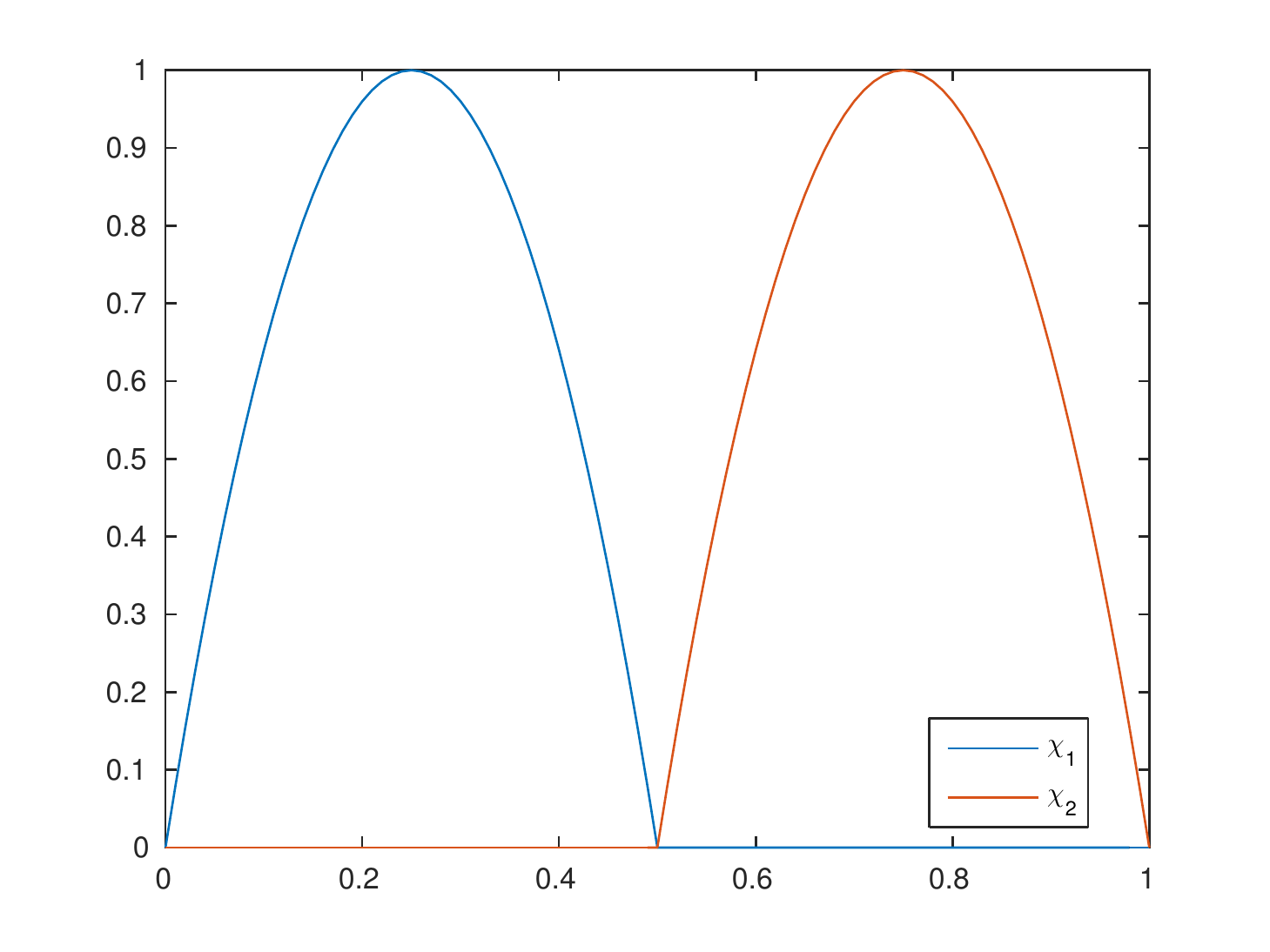} \includegraphics[scale=0.33]{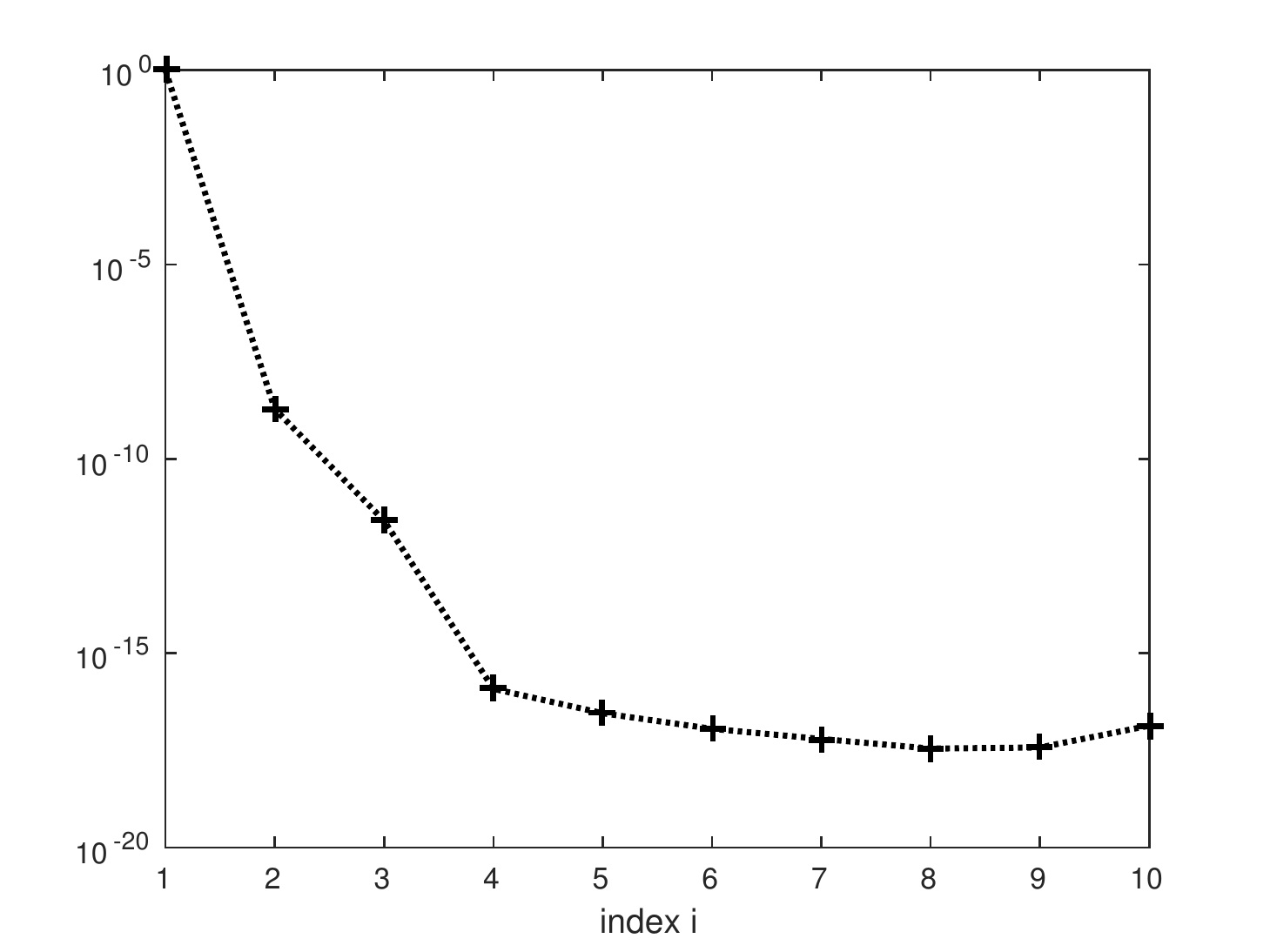} \includegraphics[scale=0.33]{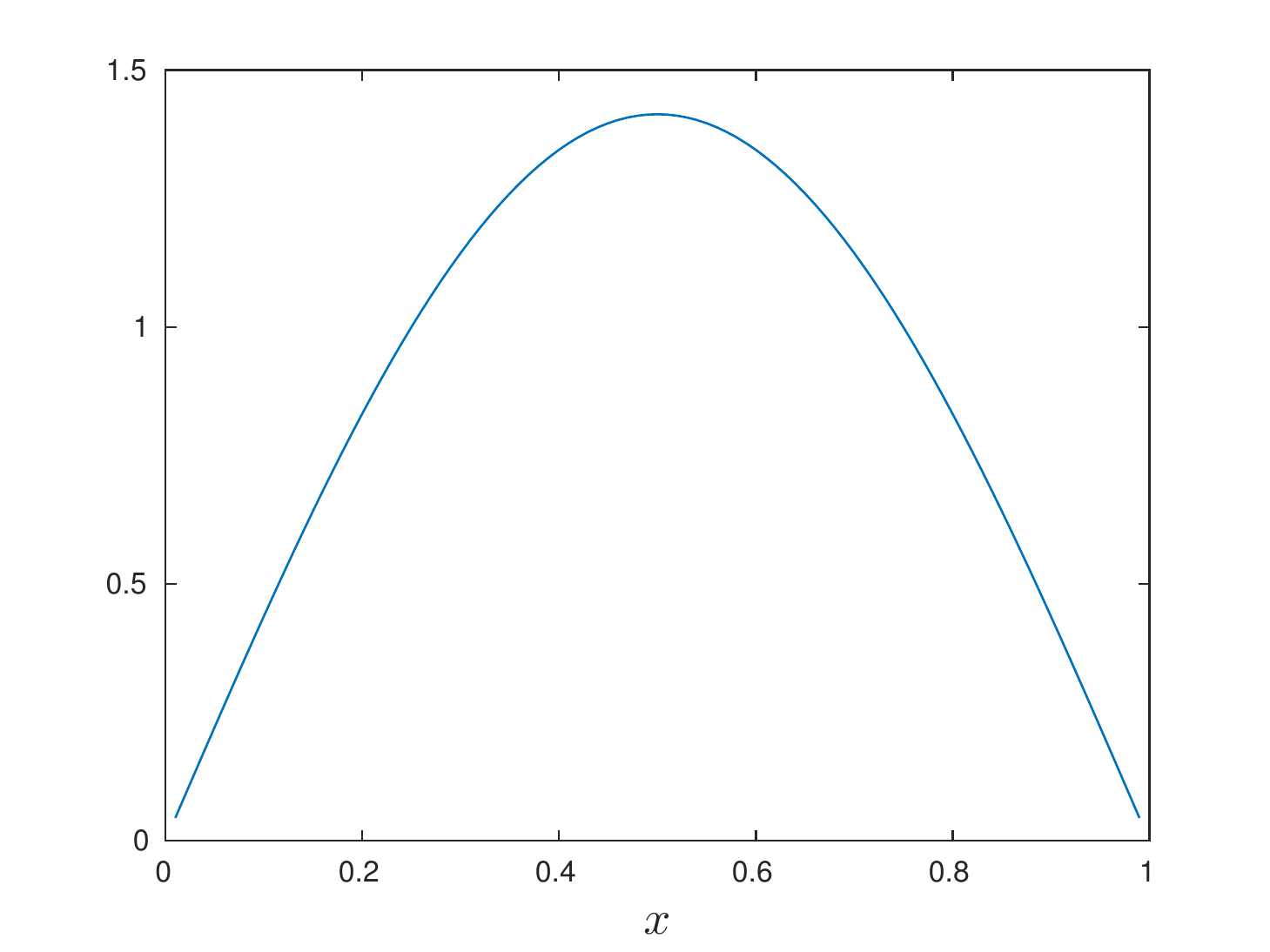}
 \caption{Test 2: Shape functions $\chi_1 (x)$ and $\chi_2(x)$ (left), decay 
 of the eigenvalues on semilog scale (middle) and first POD basis 
 function $\psi_1$ (right) utilizing uniform time grid 
 with $\Delta t = 1/40$ }
 \label{fig:shape_functions_ev}
 \end{figure}

\begin{figure}[htbp]
\includegraphics[scale=0.33]{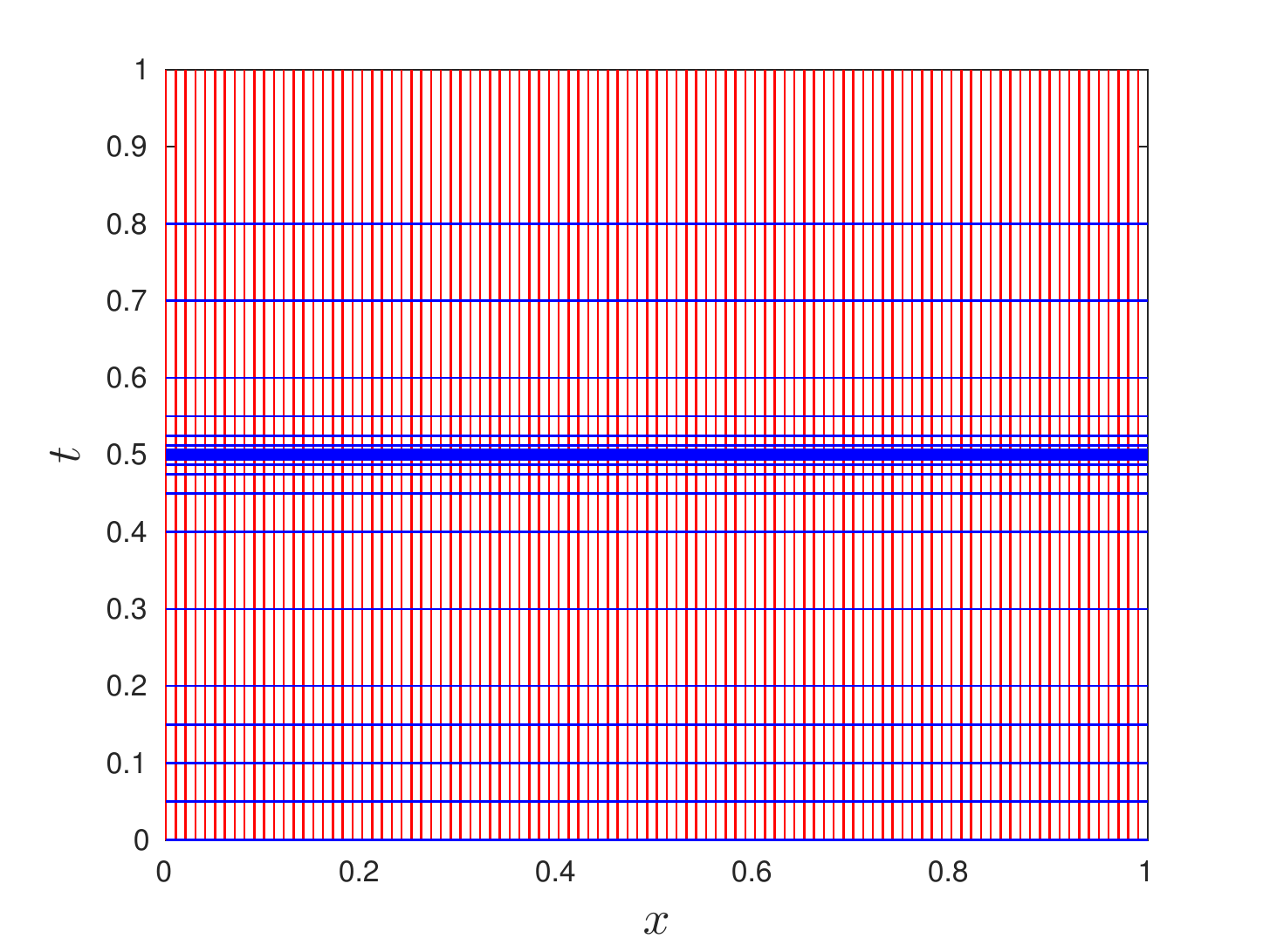}
\includegraphics[scale=0.33]{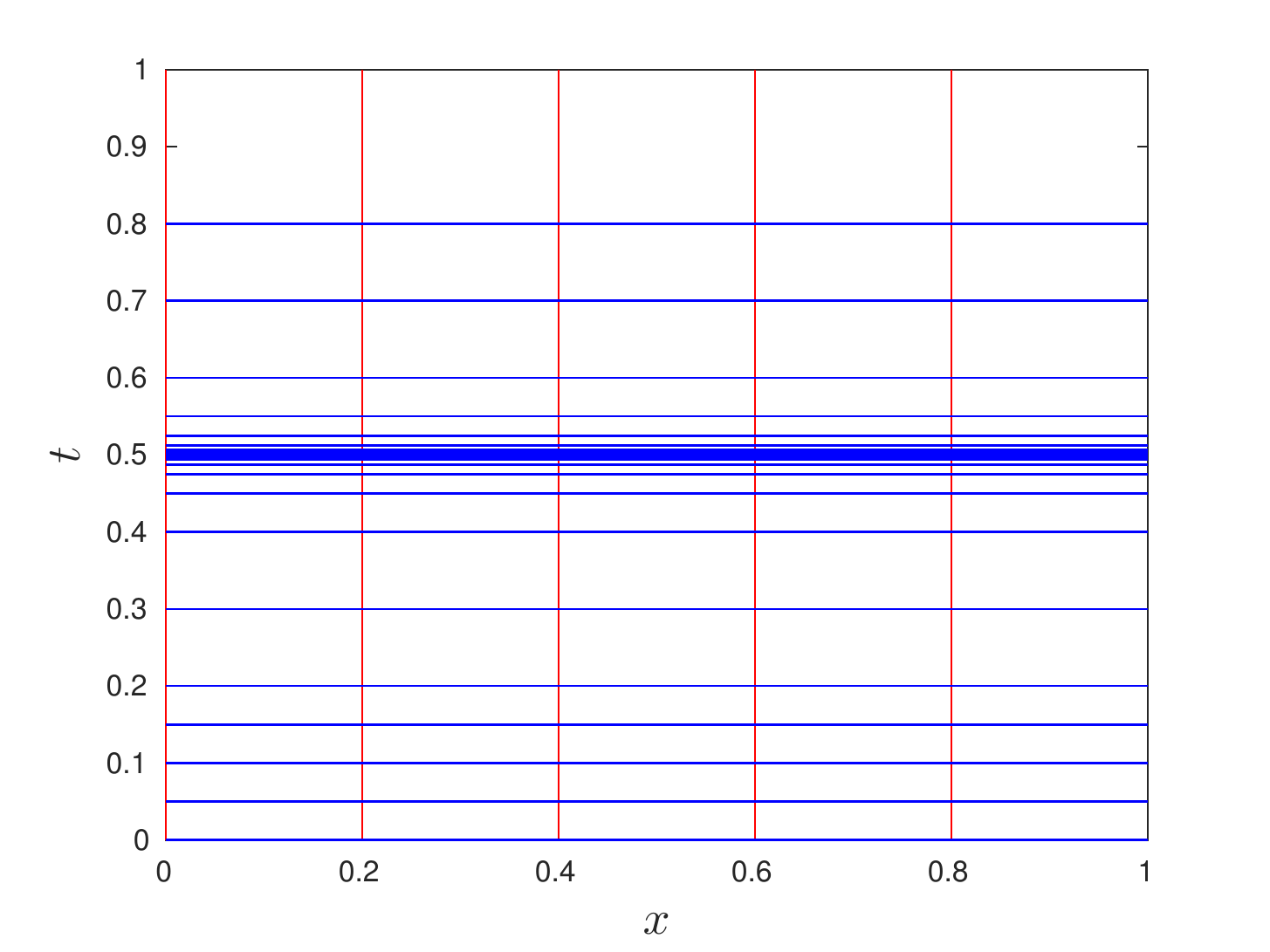}
\includegraphics[scale=0.33]{grid_equi.pdf}
    \caption{Test 2: Adaptive space-time grids with dof $= 41$
    according to the strategy in \cite{GHZ12} and 
    $\Delta x = 1/100$ (left) and $\Delta x = 1/5$ (middle), 
    respectively, and the equidistant grid with $\Delta t = 1/40$ (right)}
    \label{fig:grids_ex2}
\end{figure}

  \noindent The desired state and the forcing term are chosen accordingly. 
  Due to the arcus-tangens term and the small value for $\varepsilon$, 
  the adjoint state exhibits an interior 
  layer with steep gradient at $t=0.5$, which can be seen in the left panel of Figure \ref{fig:ex2_p_surf} and \ref{fig:ex2_p_contour}. The shape functions $\chi_1$ and $\chi_2$ are shown in 
  Figure \ref{fig:shape_functions_ev} on the left side.  
 Like in Test 1, we study the use of two different time grids: 
 an equidistant time discretization and the 
 time adaptive grid computed in step 1 
 of Algorithm \ref{Alg:OPTPOD} (see Figure \ref{fig:grids_ex2}). Once again, we 
 note that spatial and temporal discretization decouple when computing the 
 time adaptive grid utilizing the a-posteriori estimation \eqref{est-thm31}, which enables us to use a large 
 spatial resolution $\Delta x$ for solving the elliptic system and to keep the offline 
 costs low.

 As snaphots, we choose state and adjoint snapshots as well as time derivative adjoint 
 snapshots corresponding to $u_\circ=0$ and we also include the initial condition $y_0$ into 
 our snapshot set. The 
 middle and right plots of Figures \ref{fig:ex2_p_surf} and 
 \ref{fig:ex2_p_contour} show the surface and contour lines of the POD adjoint solution 
 utilizing an equidistant time grid (with $\Delta t = 1/40$) and utilizing 
 the adaptive time grid (with dof = 41), respectively. Clearly, the equidistant 
 time grid fails to capture the interior layer at $t=1/2$ satisfactorily, whereas 
 the POD adjoint state utilizing the adaptive time grid approximates the interior 
 layer well.

\begin{figure}[htbp]
\includegraphics[scale=0.33]{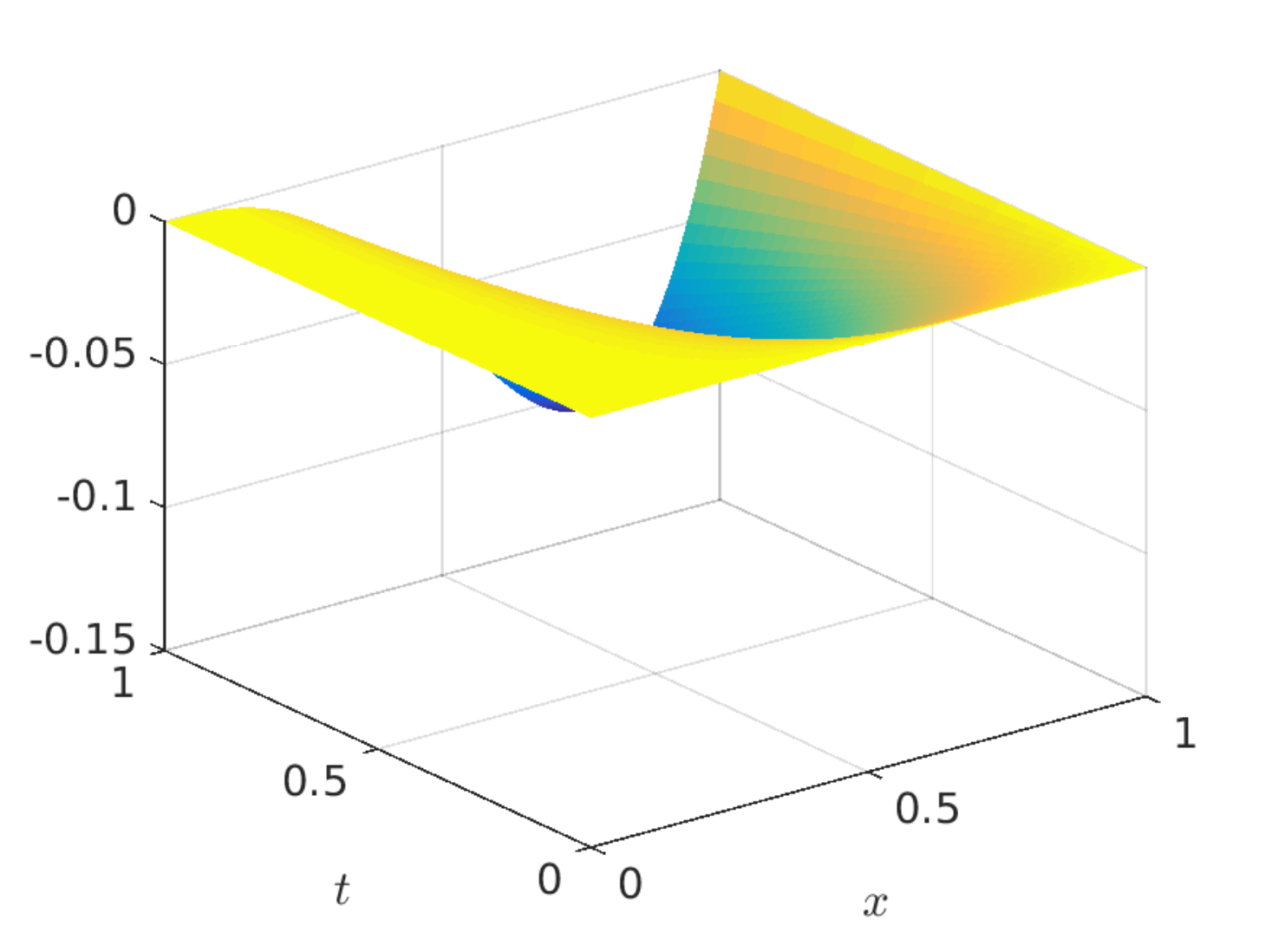} \includegraphics[scale=0.33]{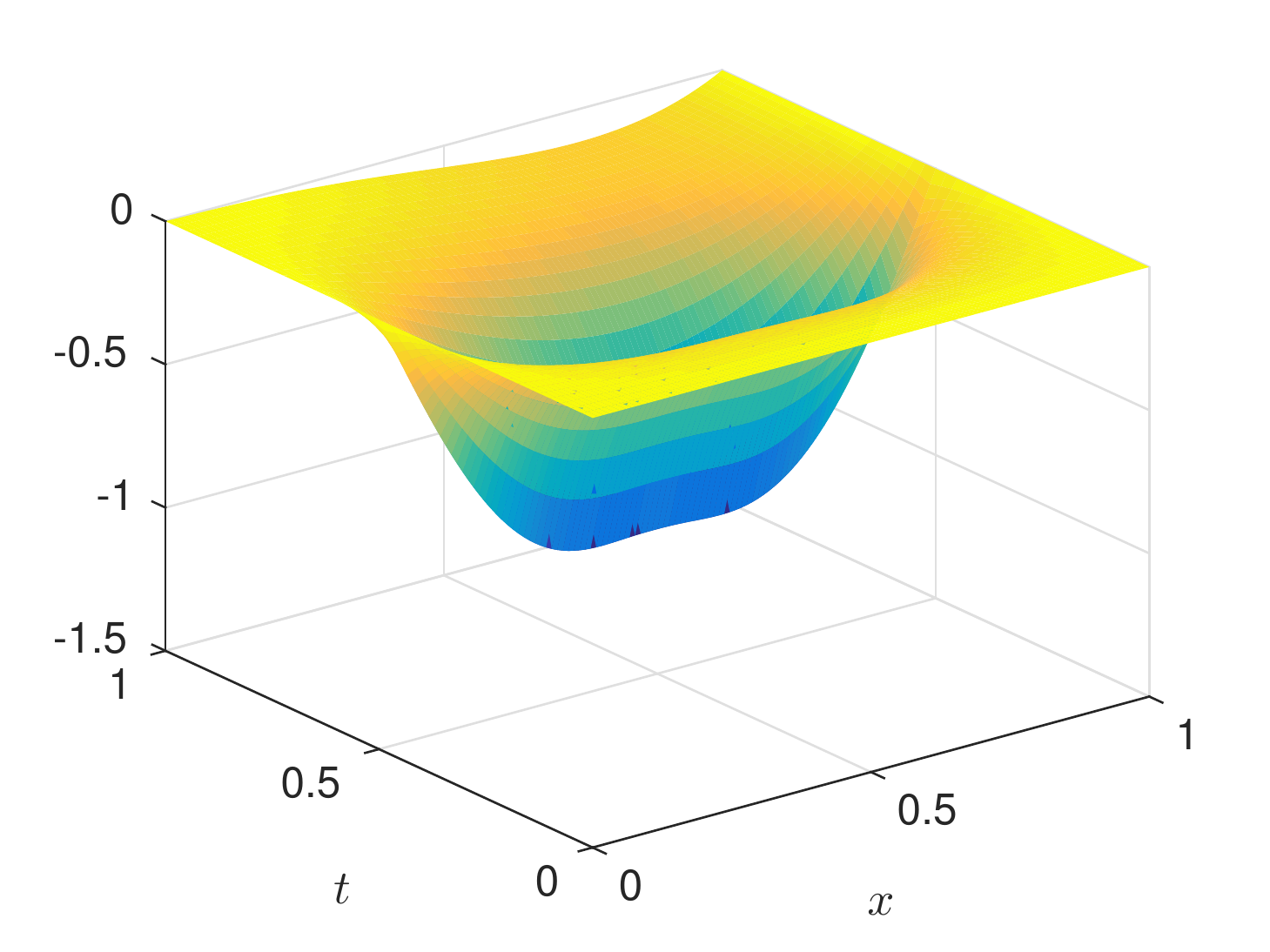} \includegraphics[scale=0.33]{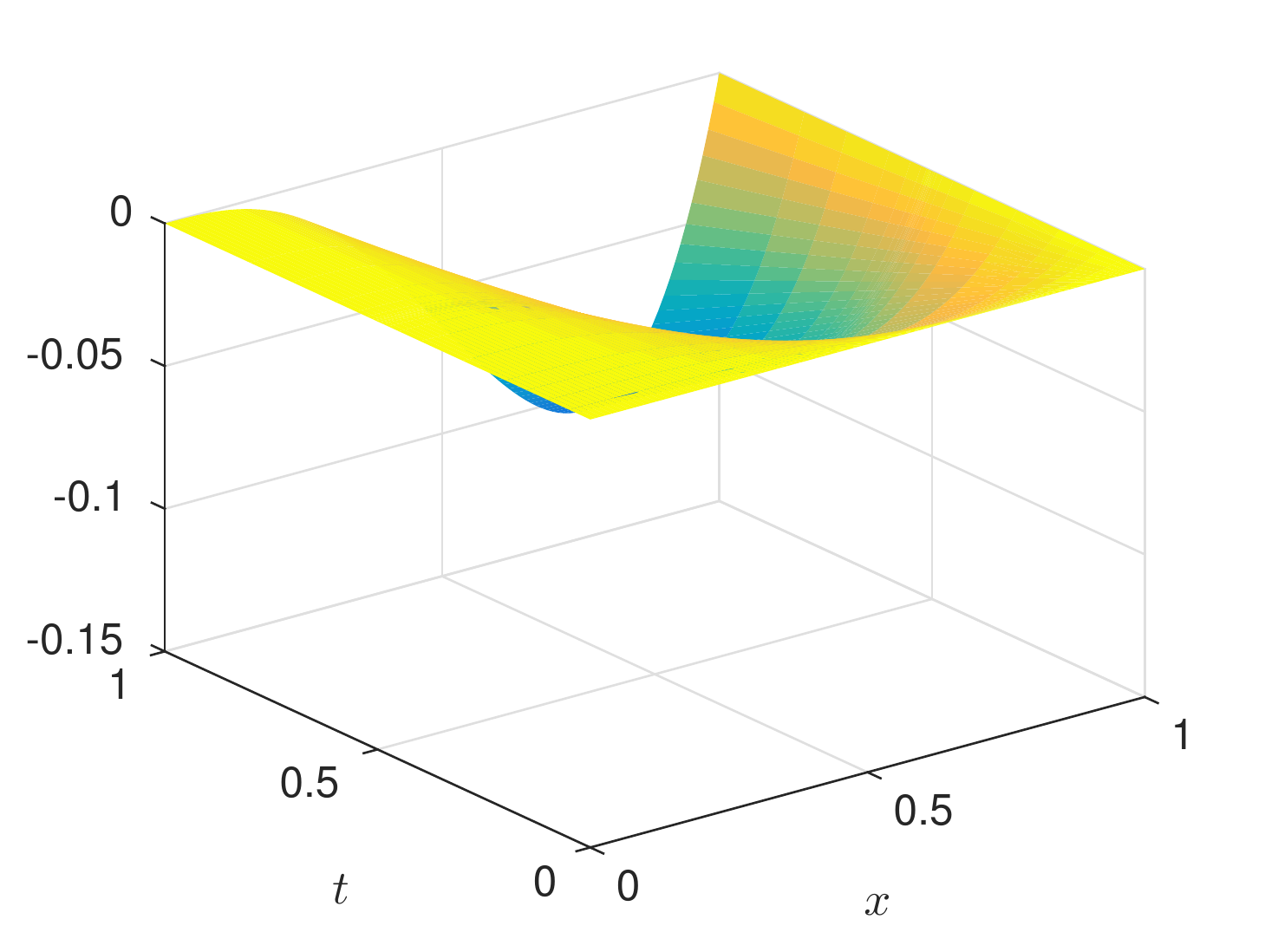}
 \caption{Test 6.2: Analytical optimal state $\bar{y}$ (left), POD  
 solution $y^\ell$ with $\ell = 4$ utilizing an equidistant time grid with
 $ \Delta t = 1/40$ (middle), 
 POD solution $y^\ell$ with $\ell = 4$ utilizing an adaptive time grid with dof=41 (right)}
 \label{fig:ex_2_state_surface} 
 \end{figure}

\begin{figure}[htbp]
\includegraphics[scale=0.33]{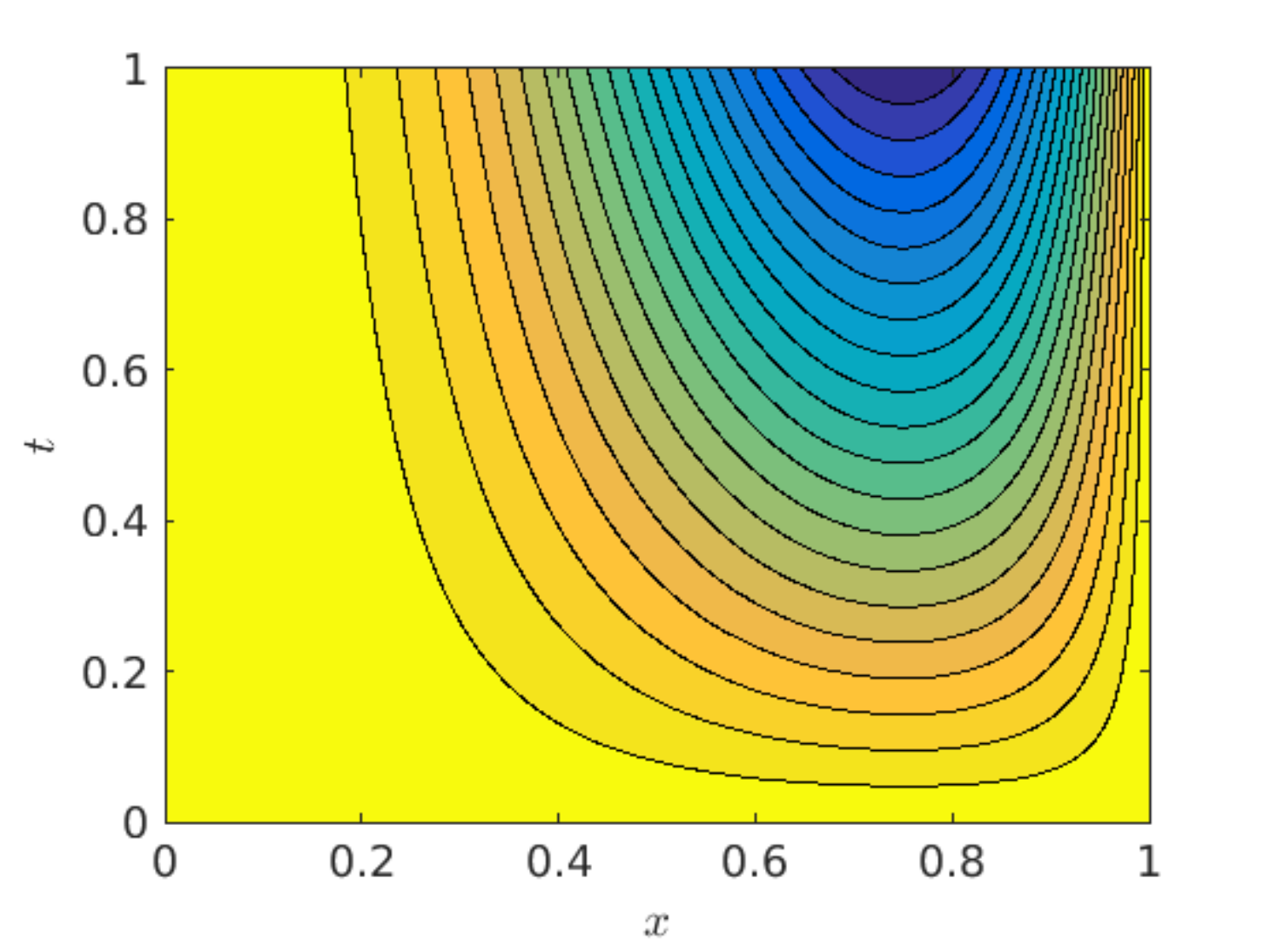} \includegraphics[scale=0.33]{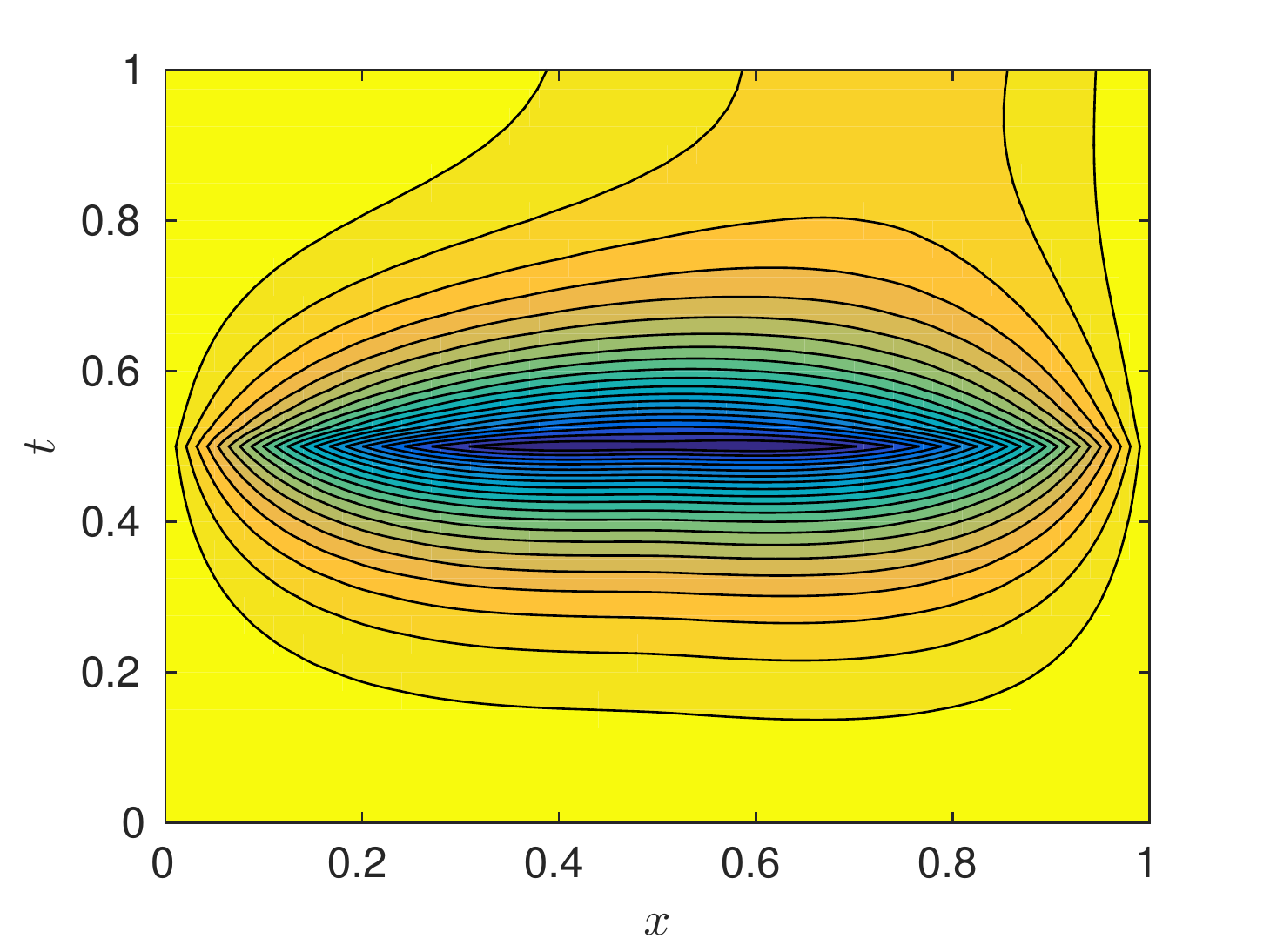} \includegraphics[scale=0.33]{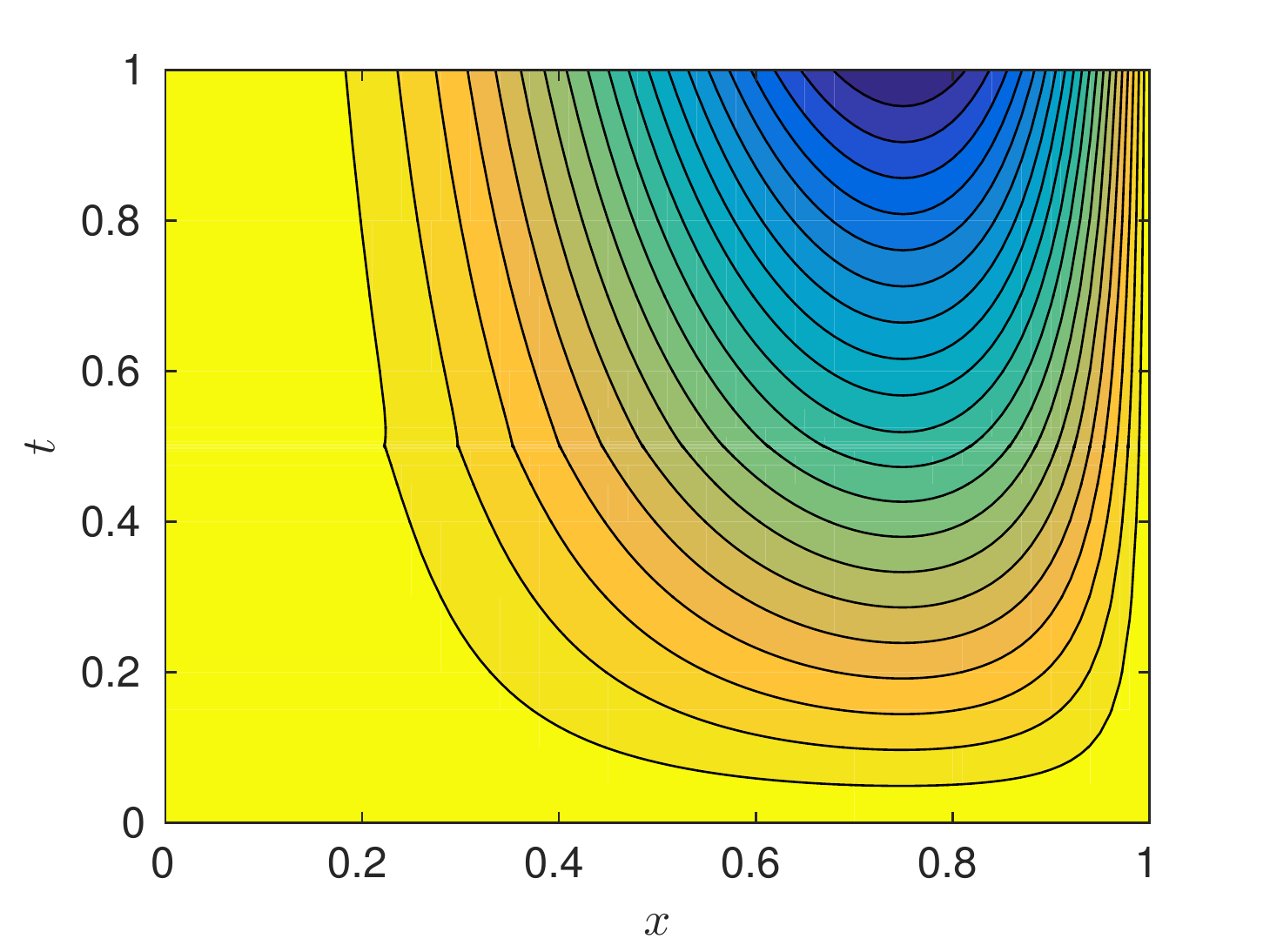}
 \caption{Test 6.2: Contour lines of the analytical optimal state 
 $\bar{y}$ (left), POD solution $y^\ell$ with $\ell = 4$ utilizing an equidistant 
 time grid with $ \Delta t = 1/40$ (middle), POD solution $y^\ell$ with $\ell = 4$ utilizing an adaptive 
 time grid with dof=41 (right)}
 \label{fig:ex_2_state_contour}
 \end{figure}

 Unlike Test Example 6.1, the adaptive time grid is also a suitable time 
grid for the state variable in this numerical test example. This can be seen visually 
when comparing the results for the POD state utilizing uniform discretization and 
utilizing the adaptive time grid with the analytical optimal state, Figures \ref{fig:ex_2_state_surface} 
and \ref{fig:ex_2_state_contour}.

 Table \ref{tab:ex2_absolute_err_l5} summarizes the absolute errors between the 
 analytical optimal solution and the POD solution for the state, control and 
 adjoint state for all test runs with an equidistant 
 and adaptive time grid, respectively. If we compare the results of the numerical approximation, we note that the use of an adaptive time grid heavily improves the quality of the POD solution with respect to an equidistant grid. In fact, we get an 
 improvement up to order four.
 
 \begin{table}[htbp]
\centering
 \begin{tabular}{ c | c | c | c || c | c | c | c}
 \toprule
 $\Delta t$ & $\varepsilon_{\text{abs}}^y$ 
 & $\varepsilon_{\text{abs}}^u$ & $\varepsilon_{\text{abs}}^p$ & dof & 
 $\varepsilon_{\text{abs}}^y$ & $\varepsilon_{\text{abs}}^u$ &  $\varepsilon_{\text{abs}}^p$ \\
 \hline
 1/20 & $5.0767 \cdot 10^{-01}$  & $7.8419 \cdot 10^{+00}$  & $3.5413 \cdot 10^{+01}$  & 21 &  $4.0346 \cdot 10^{-02}$ &  $5.4053 \cdot 10^{-01}$ &  $2.4409 \cdot 10^{+00}$ \\
 1/40 & $2.6242 \cdot 10^{-01}$  & $4.1058 \cdot 10^{+00}$  & $1.8542 \cdot 10^{+01}$  & 41 &  $2.2178 \cdot 10^{-04}$ &  $5.3471 \cdot 10^{-03}$ &  $1.3186 \cdot 10^{-02}$ \\
 1/68 & $1.5603 \cdot 10^{-01}$  & $2.4503  \cdot 10^{+00}$  & $1.1065 \cdot 10^{+01}$  & 69 &  $9.7031 \cdot 10^{-05}$ &  $4.5702 \cdot 10^{-03}$ &  $4.2670 \cdot 10^{-03}$ \\
 1/134 & $7.8741 \cdot 10^{-02}$  & $1.2386 \cdot 10^{+00}$  & $5.5938 \cdot 10^{+00}$  & 135 &  $8.5577 \cdot 10^{-05}$ &  $4.4901 \cdot 10^{-03}$ &   $2.3507 \cdot 10^{-03}$\\
 \bottomrule 
 \end{tabular}
 \vspace{0.4cm} \caption{Test 6.2: Absolute errors between the analytical optimal 
 solution and the POD solution with $\ell = 4$ depending on the time discretization (equidistant: 
 columns 1-4, adaptive: columns 5-8)}
 \label{tab:ex2_absolute_err_l5}
  \end{table}
 
\begin{figure}[htbp]
\includegraphics[scale=0.28]{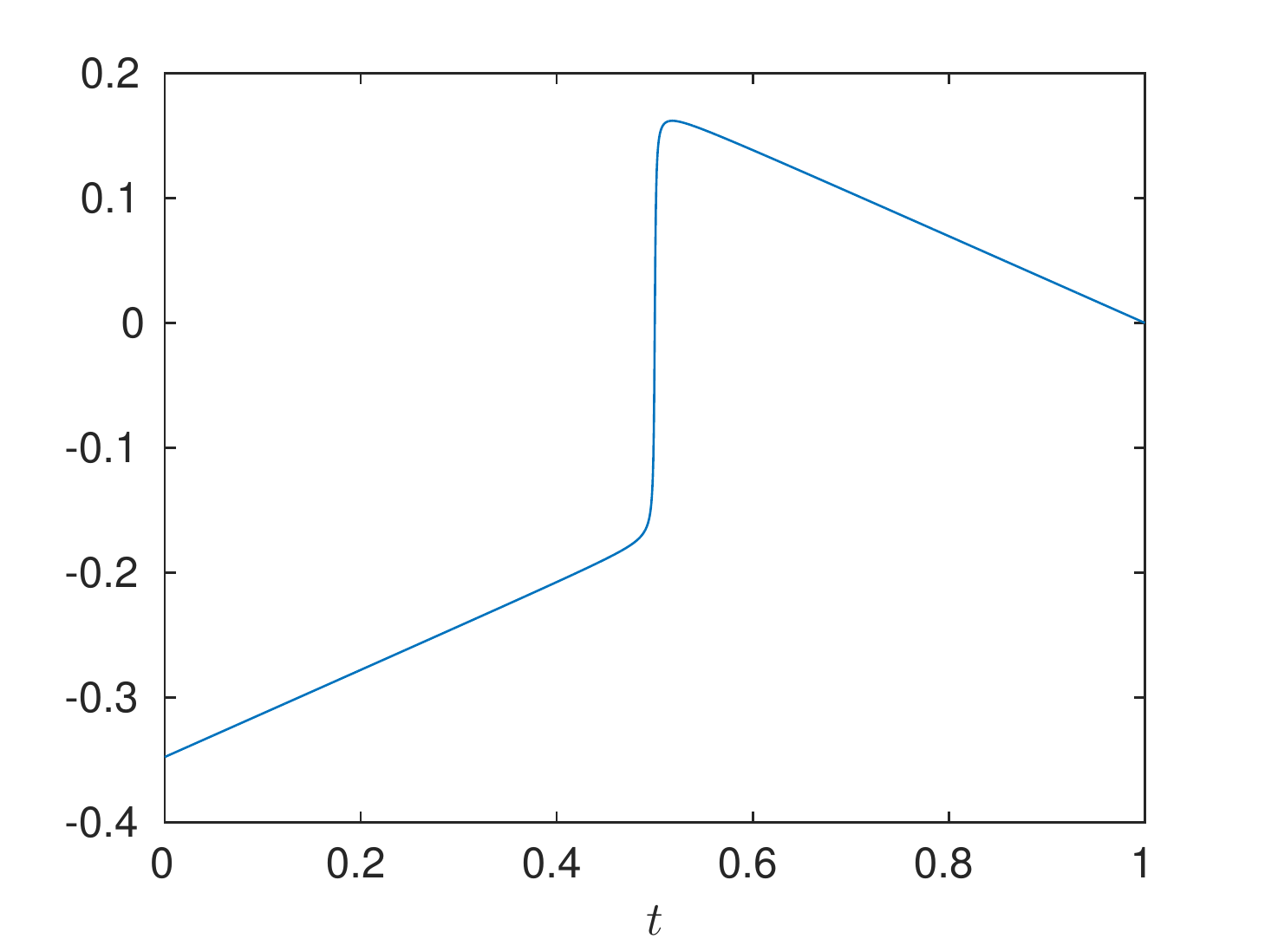}  \includegraphics[scale=0.28]{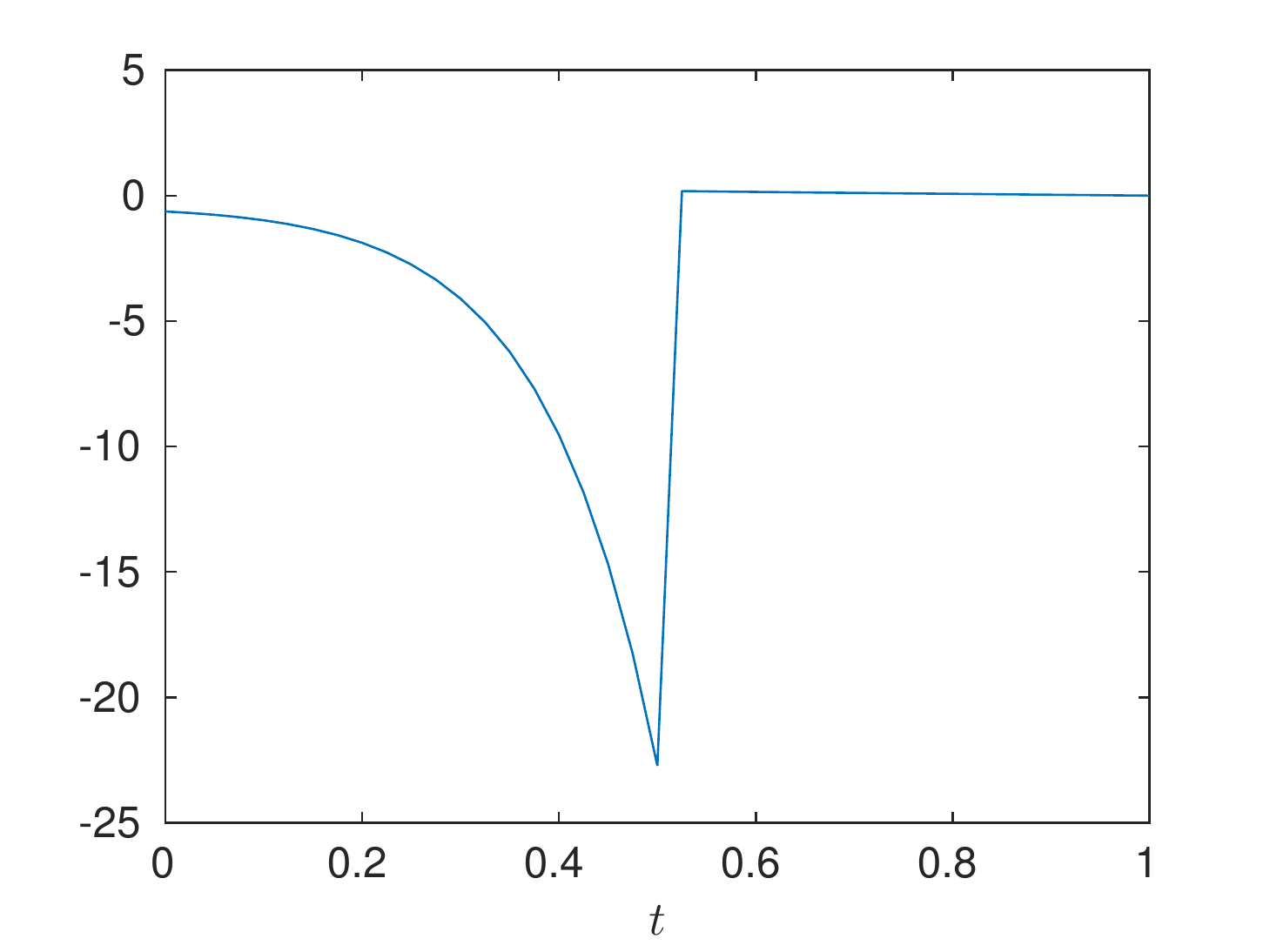}   \includegraphics[scale=0.28]{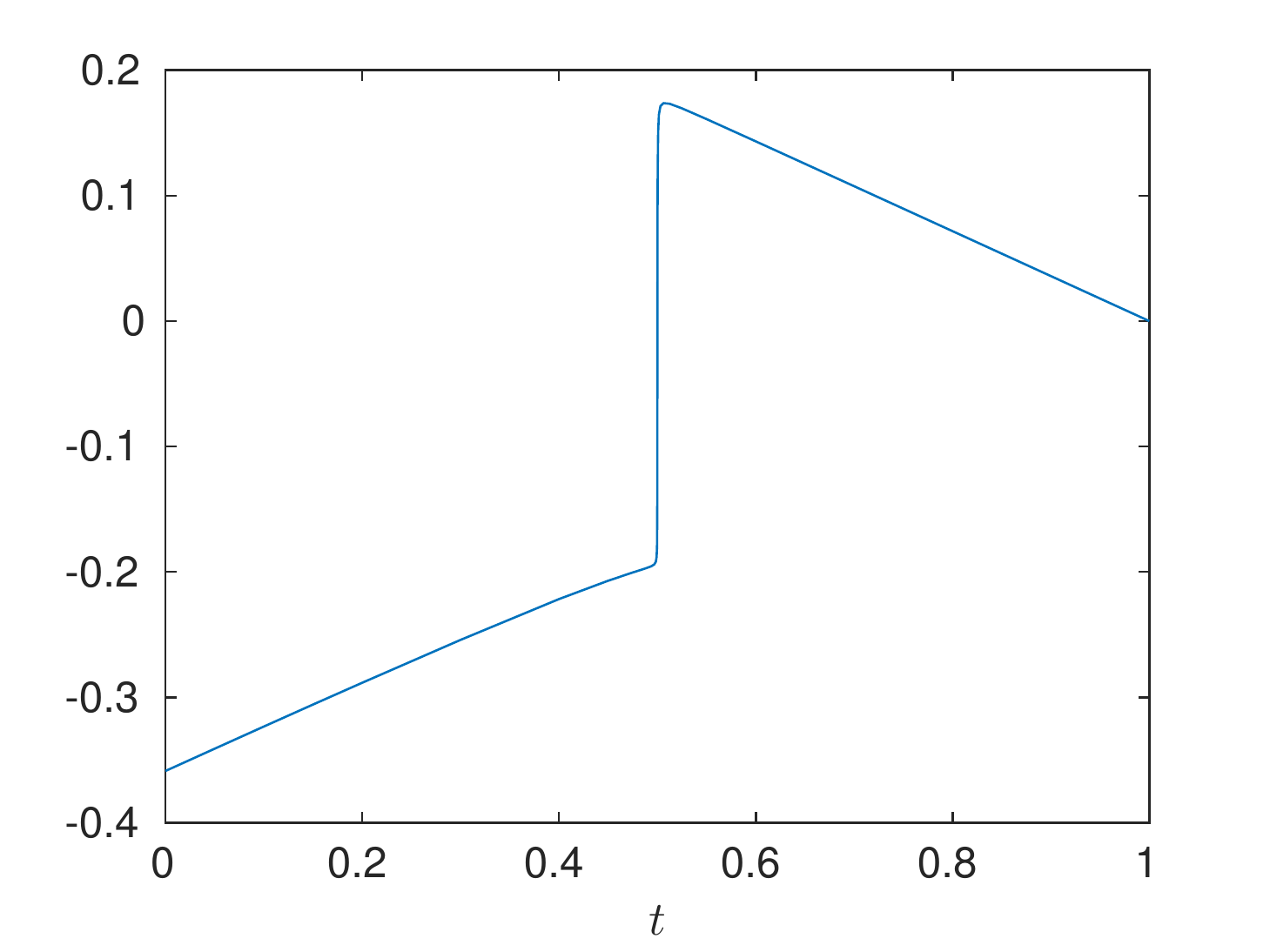}\\
  \includegraphics[scale=0.28]{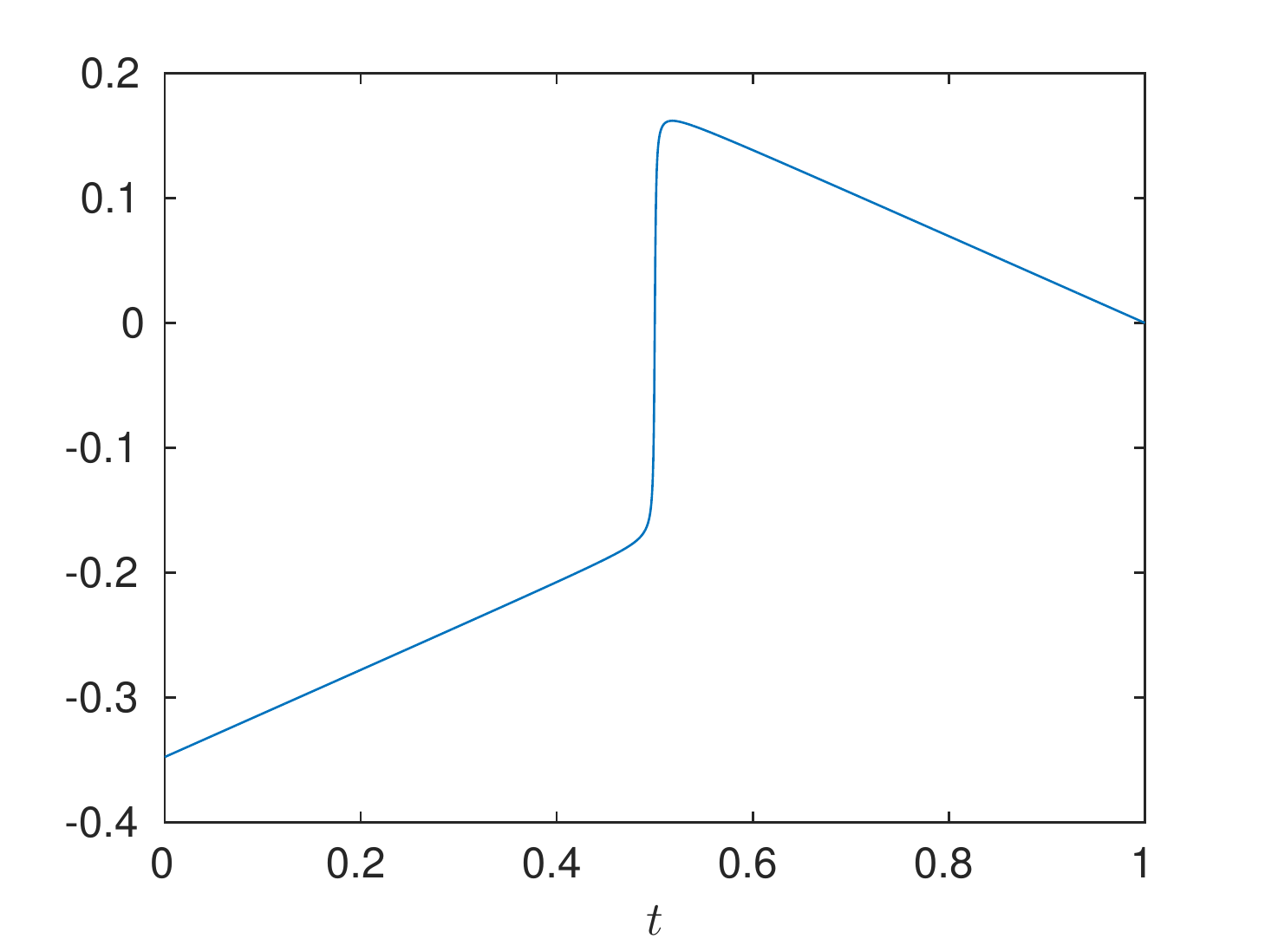} \includegraphics[scale=0.28]{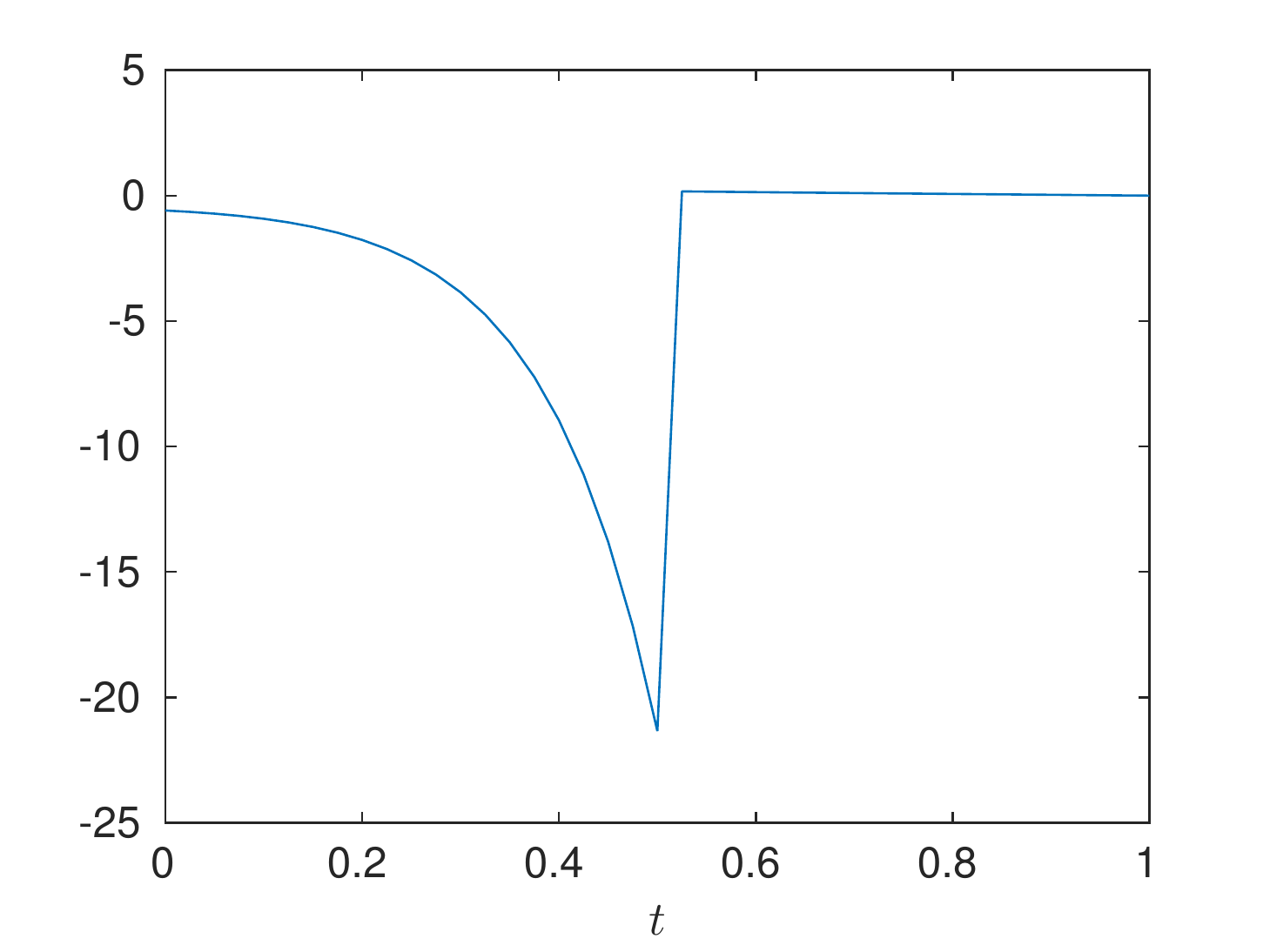}   \includegraphics[scale=0.28]{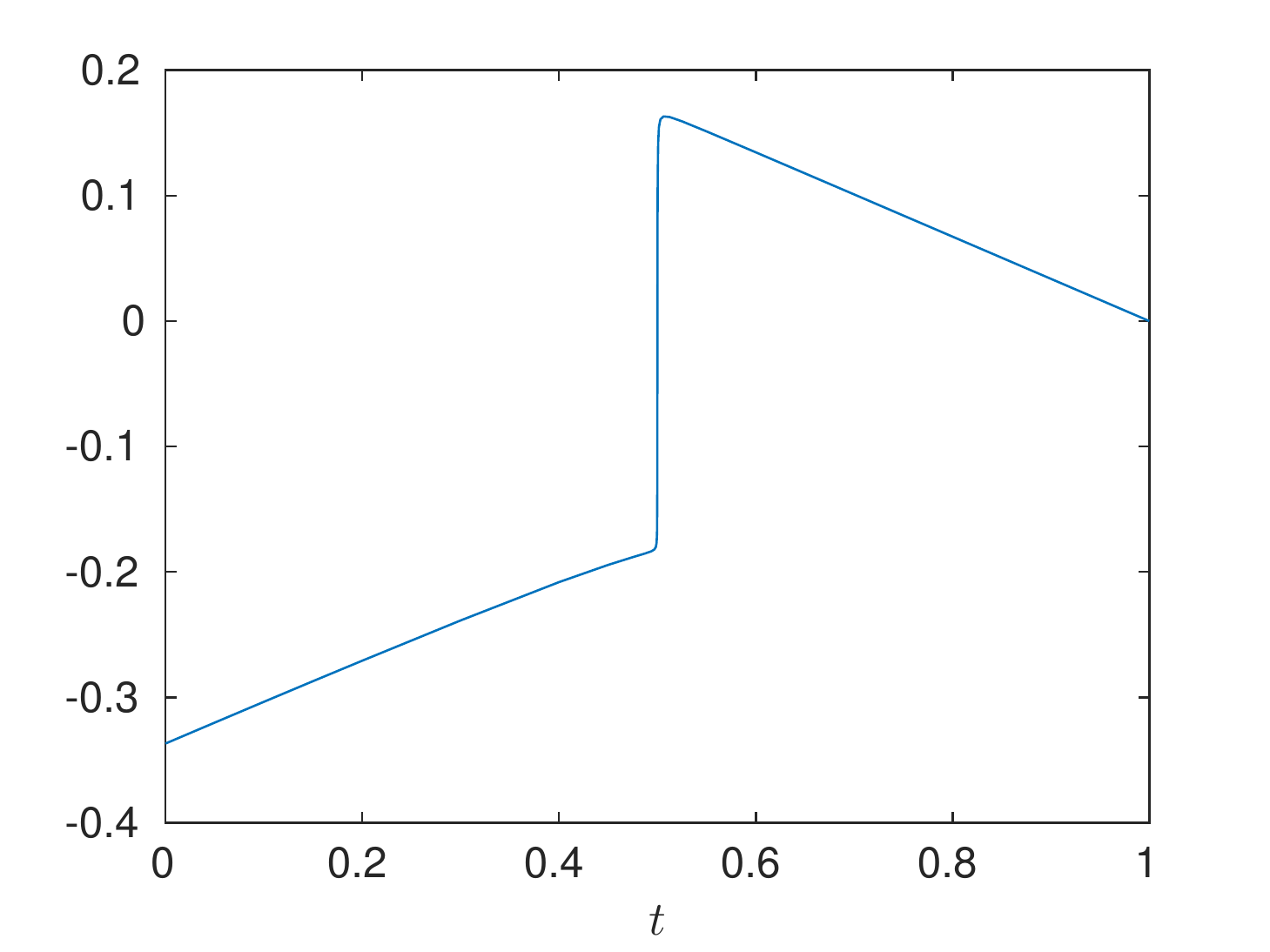}
 \caption{Test 2: Analytical control intensities $\bar{u}_1(t)$ (top left) and $\bar{u}_2(t)$ 
 (bottom left), POD control utilizing an equidistant 
 time grid with $ \Delta t = 1/40$ (middle) and $\ell = 4$, POD control 
 utilizing an adaptive time grid with dof=41 (right) and $\ell = 4$}
 \label{fig:ex_2_controls}
 \end{figure}

 The exact optimal control intensities $\bar{u}_1(t)$ and $\bar{u}_2(t)$ as well 
 as the POD solutions utilizing uniform and adaptive temporal discretization are 
 illustrated in Figure \ref{fig:ex_2_controls}.\\
 Another point of comparison is the evaluation of the cost functional. Since the 
 exact optimal solution is known analytically, we can compute the exact value 
 of the cost functional, which is $J(\bar{y},\bar{u})= 1.0085 \cdot 10^{+03}$. As expected, 
 utilizing an adaptive time grid enables us to approximate this value of the 
 cost functional quite well when using dof=135, see 
 Table \ref{tab:ex_2_cost_value}. In contrast, the use of a very fine 
 temporal discretization with $\Delta t = 1 / 10000$ is still worse than the 
 results with the adaptive time grid and dof $\geq 41$. Again, this emphasizes 
 the importance of a suitable time grid. \\

\begin{table}[htbp]
\centering
 \begin{tabular}{ c | c  || c | c }
 \toprule
 $\Delta t$
 &  $J(y^\ell, u)$   & dof  & $J(y^\ell, u)$  \\
 \hline
 1/20   & $3.1225 \cdot 10^{+05}$    & 21  &  $1.9553 \cdot 10^{+04}$ \\
 1/40   & $1.5619 \cdot 10^{+05}$   & 41  &  $1.0274 \cdot 10^{+03}$  \\
 1/68   & $9.1901  \cdot 10^{+04}$   & 69  &  $1.0065 \cdot 10^{+03}$  \\
 1/134   & $4.6655 \cdot 10^{+04}$   & 135  & $1.0082 \cdot 10^{+03}$ \\
  1/10000   & $1.0350 \cdot 10^{+03}$   & -- & --  \\
 \bottomrule 
 \end{tabular}
 \vspace{0.4cm} \caption{Test 2: Value of the cost functional with $\ell = 4$, true 
 value $J \approx 1.0085 \cdot 10^{+03}$}
 \label{tab:ex_2_cost_value}
 \end{table}

 Now, we like to investigate which influence the number $\ell$ of utilized POD basis functions 
 has on the approximation quality of the POD solution. First, we have a look at the 
 decay of the eigenvalues, which is displayed in Figure \ref{fig:shape_functions_ev}, middle.
 The eigenvalues stagnate nearby the order of machine precision, which is why the 
 use of more than $\ell = 4$ POD basis functions will not lead to better POD
 approximation results. The first POD basis function $\psi_1$ can be seen in the 
 right plot of Figure \ref{fig:shape_functions_ev}. For the use of only 
 $\ell = 1$ POD basis function, the absolute error between the analytical solution 
 and the POD solution in the state, control and adjoint state for uniform as well as for adaptive time 
 discretization are summarized in Table \ref{tab:ex2_absolute_err_l1}. Let us compare the results in this Table \ref{tab:ex2_absolute_err_l1} where $\ell = 1$ POD 
 basis function is used with the results in Table 
 \ref{tab:ex2_absolute_err_l5} where $\ell = 4$ POD basis functions are used. 
 We note that in the case of the uniform temporal discretization, the use of $\ell = 1$ 
 POD basis function leads to similar approximation results like when using $\ell = 4$ 
 POD modes. On the contrary, in the case of the adaptive time discretization, the use 
 $\ell = 4$ POD basis functions leads to better approximation 
 results with respect to the state variable than using $\ell = 1$ POD basis. The approximation 
 results concerning the control and adjoint state differ only slightly when increasing the 
 number of utilized POD basis functions. Nevertheless, also for the use of only 
 $\ell = 1$ POD mode, the use of the time adaptive grid 
 leads to an improvement of the absolute errors of up to four decimal points in comparison to using a uniform 
 time grid.

\begin{table}[htbp]
\centering
 \begin{tabular}{ c | c | c | c || c | c | c | c}
 \toprule
 $\Delta t$ & $\varepsilon_{\text{abs}}^y$ 
 & $\varepsilon_{\text{abs}}^u$ & $\varepsilon_{\text{abs}}^p$ & dof & 
 $\varepsilon_{\text{abs}}^y$ & $\varepsilon_{\text{abs}}^u$ &  $\varepsilon_{\text{abs}}^p$ \\
 \hline
 1/20 & $5.0631 \cdot 10^{-01}$  & $7.8420 \cdot 10^{+00}$  & $3.5413 \cdot 10^{+01}$  & 21 &  $4.5255 \cdot 10^{-02}$ &  $5.4054 \cdot 10^{-01}$ &  $2.4409 \cdot 10^{+00}$ \\
 1/40 & $2.6230 \cdot 10^{-01}$  & $4.1059 \cdot 10^{+00}$  & $1.8542 \cdot 10^{+01}$  & 41 &  $2.0721 \cdot 10^{-02}$ &  $5.3475 \cdot 10^{-03}$ &  $1.3186 \cdot 10^{-02}$ \\
 1/68 & $1.5684 \cdot 10^{-01}$  & $2.4503  \cdot 10^{+00}$  & $1.1065 \cdot 10^{+01}$  & 69 &  $2.0713 \cdot 10^{-02}$ &  $4.5706 \cdot 10^{-03}$ &  $4.2670 \cdot 10^{-03}$ \\
 1/134 & $8.1129 \cdot 10^{-02}$  & $1.2386 \cdot 10^{+00}$  & $5.5938 \cdot 10^{+00}$  & 135 &  $2.0664 \cdot 10^{-02}$ &  $4.4905 \cdot 10^{-03}$ &   $2.3507 \cdot 10^{-03}$\\
 \bottomrule 
 \end{tabular}
 \vspace{0.4cm} \caption{Test 2: Absolute errors between the analytical optimal 
 solution and the POD solution with $\ell = 1$ depending on the time discretization (equidistant: 
 columns 1-4, adaptive: columns 5-8)}
 \label{tab:ex2_absolute_err_l1}
  \end{table}

 \subsection{Test 3: Control constrained problem} 
In this test we add control constraints to the previous example.
We set $u_{1,a}(t) 
 \leq u_1(t) \leq u_{1,b}(t)$ and $u_{2,a}(t) \leq u_2(t) \leq u_{2,b}(t)$ for the time dependent 
 control intensities $u_1(t)$ and $u_2(t)$. The analytical value range 
 for both controls is $u_1(t), u_2(t) \in [-0.3479, 0.1700]$ for $t \in [0,1]$. For each control intensity we choose different upper and lower bounds: we set $u_{1,a}(t) = -100$ (i.e. no 
 restriction), $ u_{1,b} = 0.1$ and $u_{2,a}(t) = -0.2, \; u_{2,b}(t) = 0$. For the 
 solution of problem \eqref{ocp_pod} we use a projected gradient method.

\begin{figure}[htbp]
\includegraphics[scale=0.33]{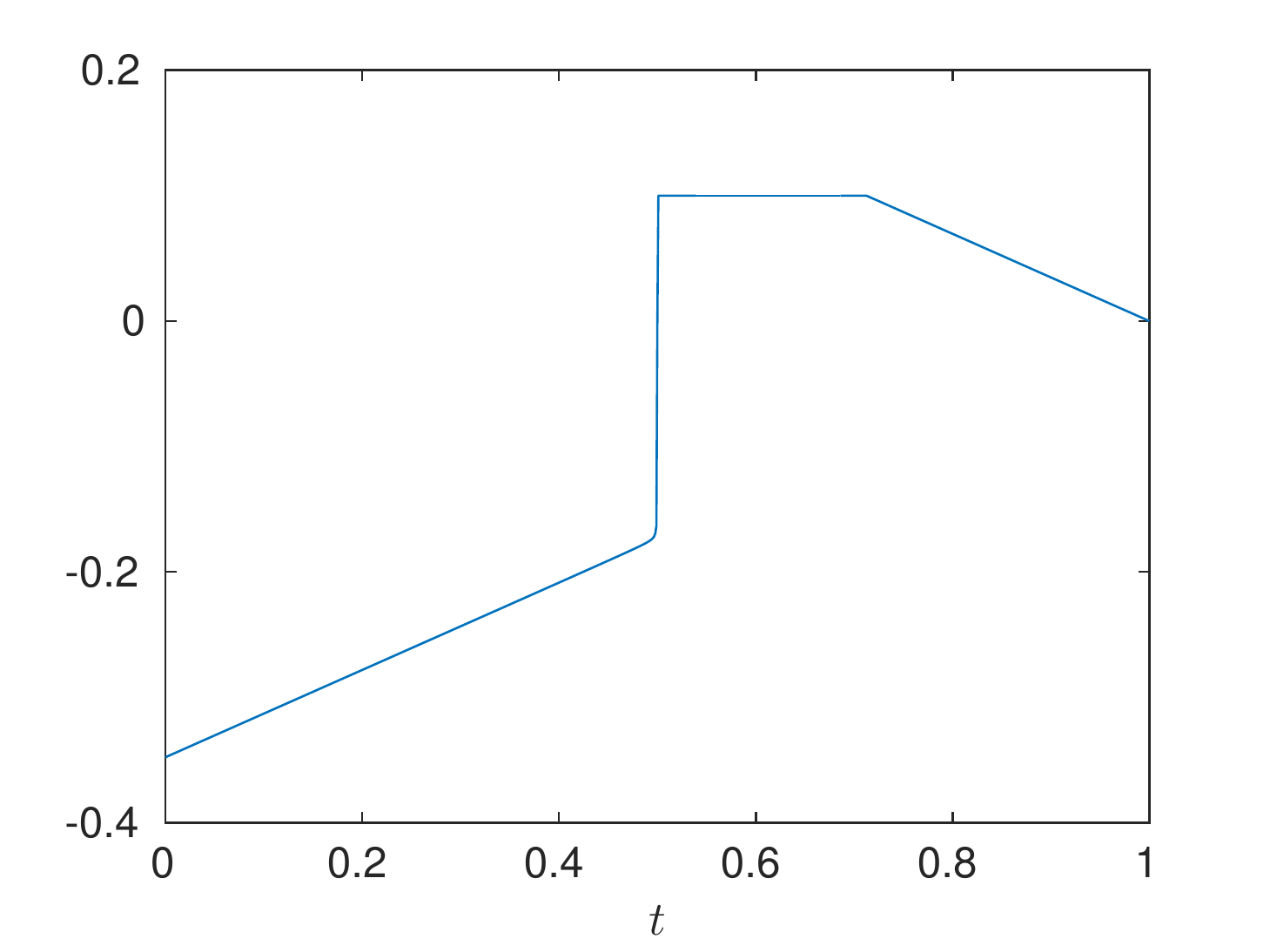}  \includegraphics[scale=0.33]{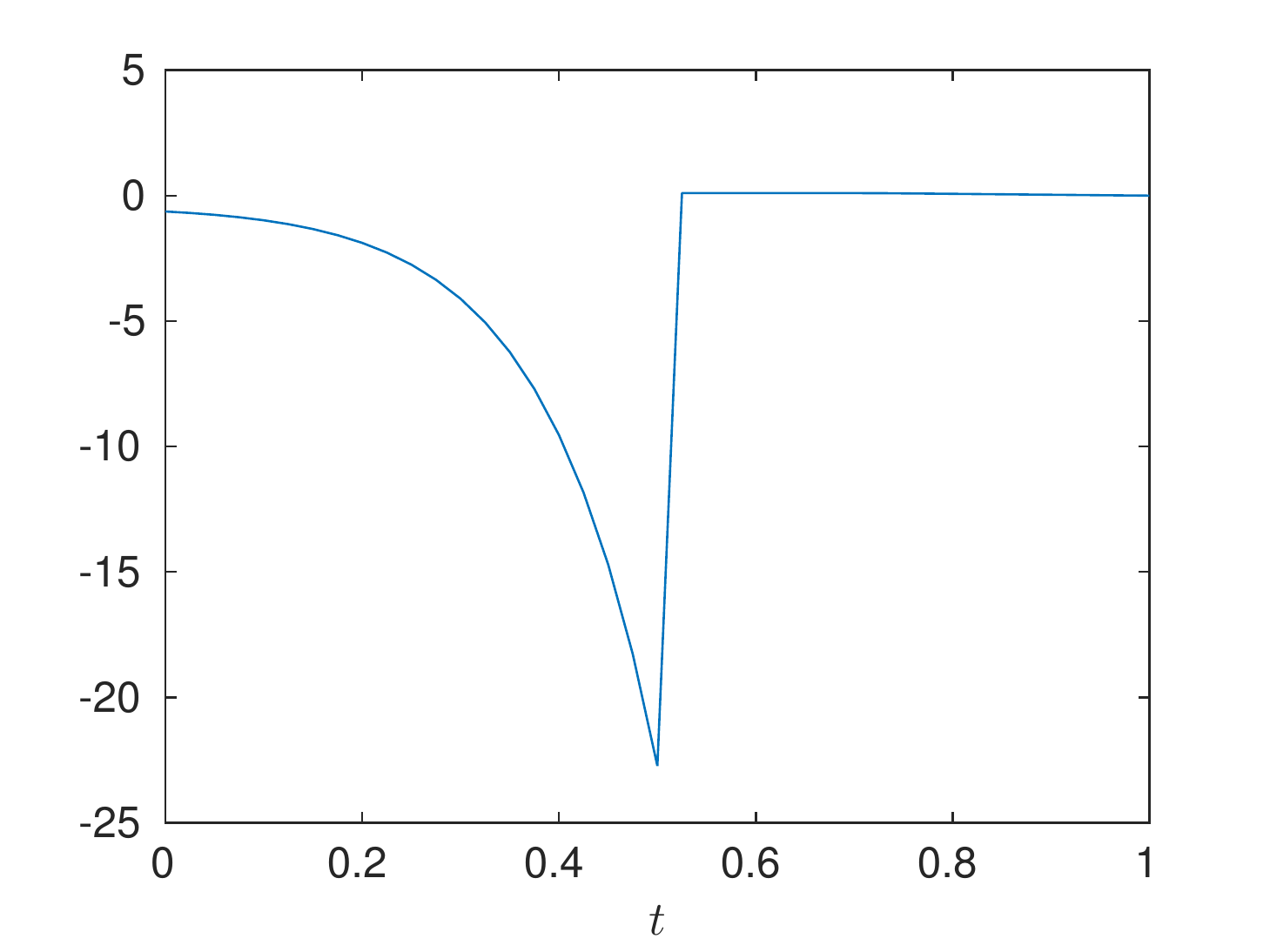}  \includegraphics[scale=0.33]{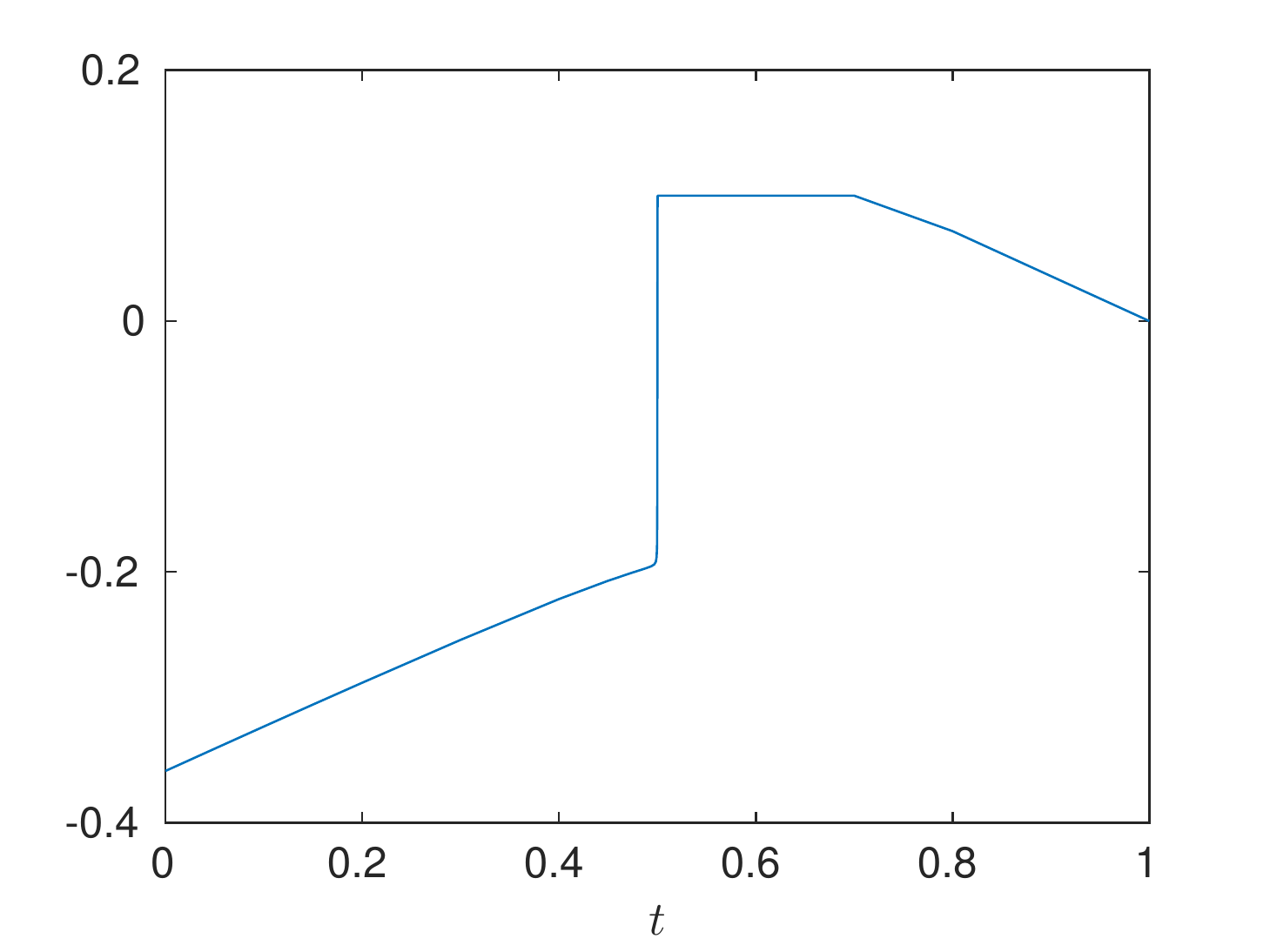}\\
  \includegraphics[scale=0.33]{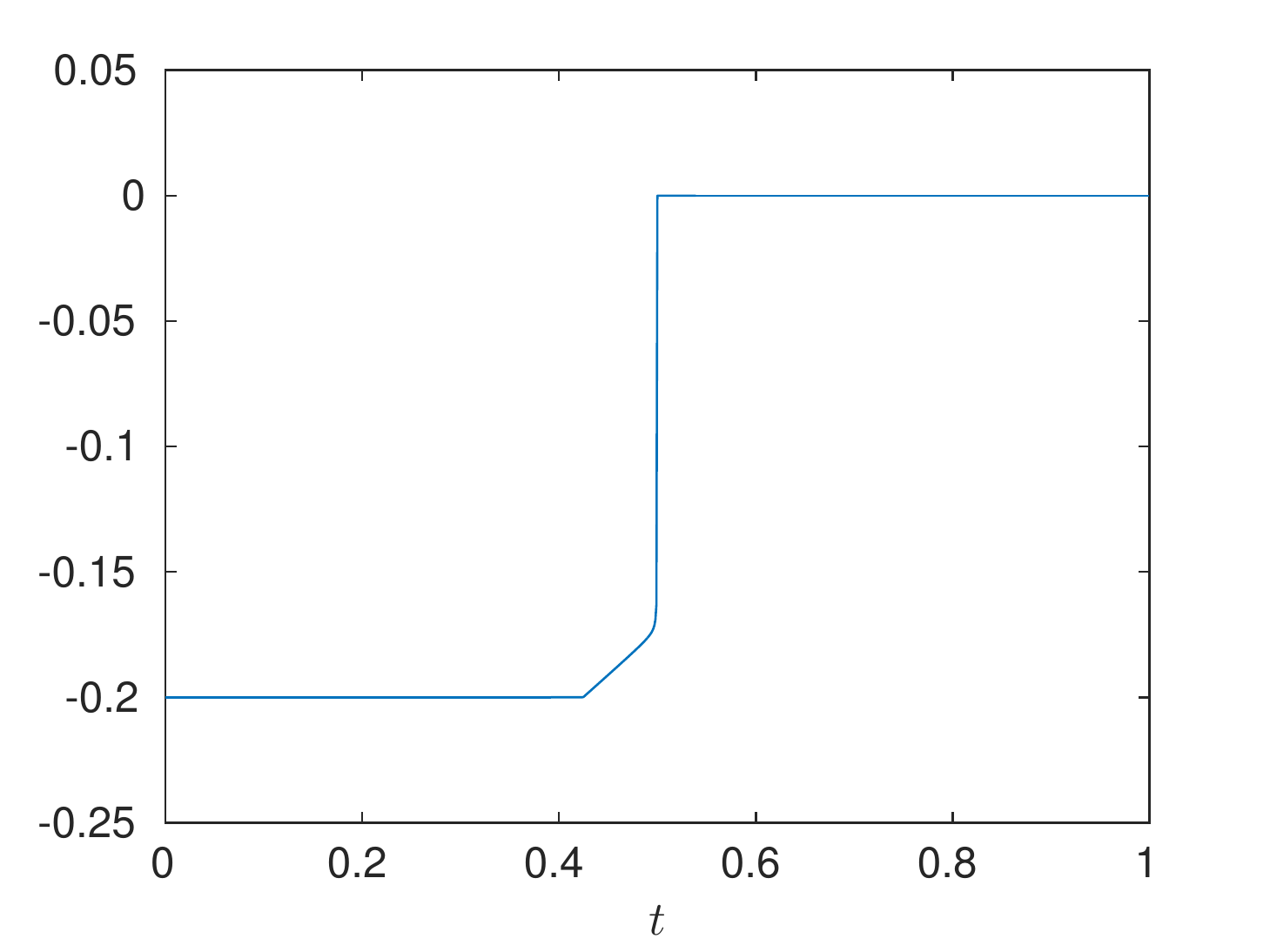} \includegraphics[scale=0.33]{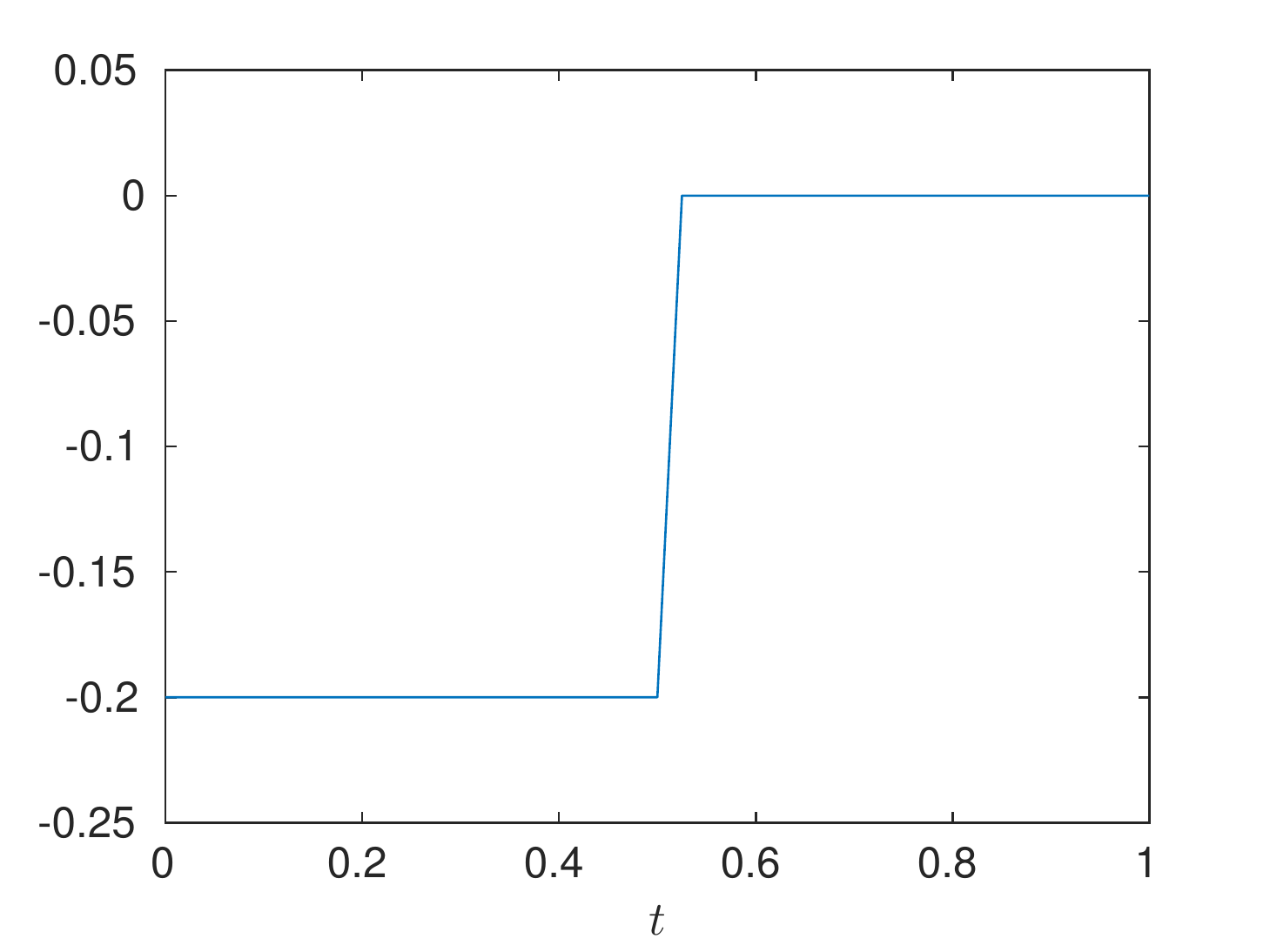}  \includegraphics[scale=0.33]{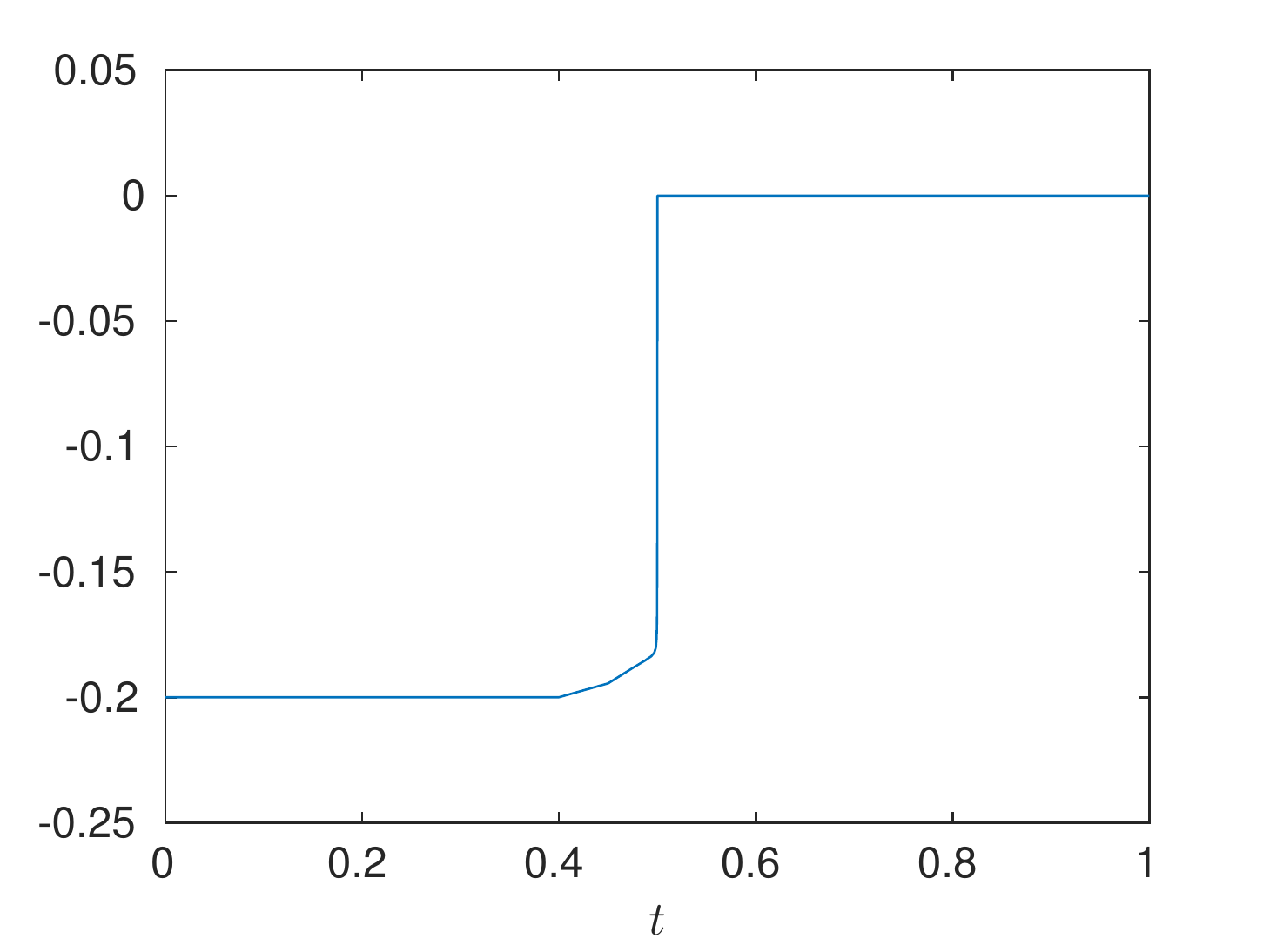}
 \caption{Test 3: Inclusion of box constraints for the control intensities: Analytical control intensities $\bar{u}_1(t)$ (top left) and $\bar{u}_2(t)$ 
 (bottom left), POD control utilizing an equidistant 
 time grid with $ \Delta t = 1/40$ (middle) and $\ell = 4$, POD control utilizing an adaptive 
 time grid with dof=41 (right) and $\ell = 4$}
 \label{fig:ex_2_controls_box}
 \end{figure}

 The solution of the nonlinear, nonsmooth equation \eqref{2ordp} can be done 
 by a semi-smooth Newton method or by a Newton method utilizing a 
 regularization of the projection formula, see \cite{NPS11}. In our numerical 
 tests we compute the approximate solution to \eqref{ocp_pod} with a fixed point 
 iteration and initialize the method with the adjoint state corresponding to the 
 control unconstrained optimal control problem. In this way, only two iterations 
 are needed for convergence. Convergence of the fixed point iteration can be argued for large enough values of $\alpha$, see \cite{HV12}.\\ 
 The analytical optimal solutions $\bar{u}_1$ and $\bar{u}_2$ are shown in the left 
 plots in Figure \ref{fig:ex_2_controls_box}. For POD basis computation, we use state, adjoint and 
 time derivative adjoint snapshots corresponding to the reference control $u_\circ = 0$ and we also include 
 the initial condition $y_0$ into our snapshot set. The plots in the middle and on the 
 right in Figure \ref{fig:ex_2_controls_box} refer to the POD controls using a uniform and an 
 adaptive temporal discretization, respectively. Once again, we note that utilizing 
 an adaptive time grid leads to far better results than using a uniform temporal 
 grid. The numerical results in Table \ref{tab:ex_2_errors_box} confirm this 
 observation. We observe that the inclusion of box constraints on the control 
 functions lead in general to better approximation results, compare Table
 \ref{tab:ex2_absolute_err_l5} with Table \ref{tab:ex_2_errors_box}. This is due to 
 the fact that on the active sets the error between the analytical optimal 
 controls and the POD solutions vanishes.

\begin{table}[htbp]
\centering
 \begin{tabular}{ c | c | c | c || c | c | c | c}
 \toprule
 $\Delta t$ & $\varepsilon_{\text{abs}}^y$ 
 & $\varepsilon_{\text{abs}}^u$ & $\varepsilon_{\text{abs}}^p$ & dof & 
 $\varepsilon_{\text{abs}}^y$ & $\varepsilon_{\text{abs}}^u$ &  $\varepsilon_{\text{abs}}^p$ \\
 \hline
 1/20 & $ 2.8601 \cdot 10^{-01}$  & $5.7201 \cdot 10^{+00}$  & $3.5430 \cdot 10^{+01}$  & 21 &  $2.2714 \cdot 10^{-02}$ &  $3.9586 \cdot 10^{-01}$ &  $2.4423 \cdot 10^{+00}$ \\
 1/40 & $ 1.4802 \cdot 10^{-01}$  & $2.9955 \cdot 10^{+00}$  & $1.8551 \cdot 10^{+01}$  & 41 &  $2.9482 \cdot 10^{-04}$ &  $4.4969 \cdot 10^{-03}$ &  $1.3183 \cdot 10^{-02}$ \\
 1/68 & $8.8124 \cdot 10^{-02}$  & $1.7882  \cdot 10^{+00}$  & $1.1071 \cdot 10^{+01}$  & 69 &  $2.1247 \cdot 10^{-04}$ &  $3.2811 \cdot 10^{-03}$ &  $4.2629 \cdot 10^{-03}$ \\
 1/134 & $4.4570 \cdot 10^{-02}$  & $9.0470 \cdot 10^{-01}$  & $5.5965 \cdot 10^{+00}$  & 135 &  $2.1330 \cdot 10^{-04}$ &  $3.1321 \cdot 10^{-03}$ &  $2.3474 \cdot 10^{-03}$\\
 \bottomrule 
 \end{tabular}
 \vspace{0.4cm} \caption{Test 3: Inclusion of box constraints for the control 
 intensities: Absolute errors between the analytical optimal 
 solution and the POD solution with $\ell = 4$ depending on the time discretization (equidistant: 
 columns 1-4, adaptive: columns 5-8)}
 \label{tab:ex_2_errors_box}
  \end{table}

 \section{Conclusion}
In this paper we investigated the problem of snapshot location in optimal control problems. We showed that the numerical POD solution is 
much more accurate if we use an adaptive time grid, especially when the solution 
of the problem presents steep gradients. The time grid was computed by means of an a-posteriori error estimation 
strategy of space-time approximation of a second order in time and fourth order in space elliptic equation which 
describes the optimal control problem and has the advantage that it is independent of an input control function. 
Furthermore, a coarse approximation with respect to space of the latter equation gives information on the snapshots one 
can use to build the surrogate model. Finally, we provided a certification of our surrogate model by 
means of an a-posteriori error estimation for the error between the optimal solution and the POD 
solution.\\
For future work, we are interested in transferring our approach to optimal control problems subject to 
nonlinear parabolic equations.



\end{document}